# AUTOMORPHISMS OF GENERALIZED THOMPSON GROUPS

## MATTHEW G. BRIN AND FERNANDO GUZMÁN

Department of Mathematical Sciences

State University of New York at Binghamton

Binghamton, NY 13902-6000

USA

August 28, 1997

## CONTENTS

## 0. INTRODUCTION

**0.1. Results.** We study the automorphisms of some generalizations of Thompson's groups and their underlying structures. The automorphism groups of two of Thompson's original groups were analyzed in [2] and were shown to be "small" and "unexotic." Our results differ sharply from [2] in that we show that the automorphism groups of the generalizations are "large" and have "exotic" elements. The term exotic is from [1] and is explained later in this introduction.

Richard J. Thompson introduced the triple of infinite groups $F \subseteq T \subseteq G$ in the 1960s and showed that they have several interesting properties. These groups were later generalized to infinite groups $F_{n,\infty} \subseteq F_n \subseteq T_{n,r} \subseteq G_{n,r}$ for $n \geq 2$ and $r \geq 1$ for which $F = F_{2,\infty} = F_2$, $T = T_{2,1}$ and $G = G_{2,1}$. Further background is given in the second part of this introduction.

The automorphism groups $\mathrm{Aut}(F)$ and $\mathrm{Aut}(T)$ were analyzed in [2]. In this paper, we study $\mathrm{Aut}(F_{n,\infty})$, $\mathrm{Aut}(F_n)$ and $\mathrm{Aut}(T_{n,n-1})$ for $n > 2$. We find that all three contain exotic elements, and that the complexity increases as $n$ increases. We also find that the structures of $\mathrm{Aut}(F_n)$ and $\mathrm{Aut}(T_{n,n-1})$ are closely related, while differing significantly from the structure of $\mathrm{Aut}(F_{n,\infty})$. This difference increases as $n$ increases. Our analysis is only partial in that we discover large subgroups of these groups containing exotic elements, but we do not know whether the subgroups are proper. Automorphism groups of a class of groups that include the Thompson groups were studied in the unpublished notes [1] where the structures were analyzed modulo the unanswered question of whether exotic automorphisms could exist.

We now give more details. Throughout the paper, functions will be written to the right of their arguments and composition will proceed from left to right.

We give two fundamental facts that underly our analyses. The first is that $F_{n,\infty}$ and $F_n$ act faithfully by homeomorphisms on the real line $\mathbf{R}$, and $T_{n,r}$ acts faithfully by homeomorphisms on the circle $S_r = \mathbf{R}/r\mathbf{Z}$ of length $r$. (The groups $G_{n,r}$, not studied, act on the Cantor set.) Further, all the homeomorphisms used in the representations satisfy certain "local structure conditions." That is, for each $n \geq 2$, the representations use only homeomorphisms that

1. preserve orientation,
2. are piecewise linear (PL) ,
3. have slopes restricted to $\langle n \rangle$, the multiplicative, cyclic subgroup of the positive reals generated by $n$,
4. have "breaks" (discontinuities of slope) restricted to a discrete subset of $\mathbf{Z}[\frac{1}{n}]$, the ring of fractions with denominators in $\langle n \rangle$, and
5. preserve $\mathbf{Z}[\frac{1}{n}]$ as a set.

In fact, $T_{n,r}$ is exactly the set of self homeomorphisms of $S_r$ that satisfy 1–5.

For $F_{n,\infty}$, there is another restriction. There is a surjective ring homomorphism $\phi_n : \mathbf{Z}[\frac{1}{n}] \to \mathbf{Z}_{n-1}$ defined by $(a/n^b)\phi_n \equiv a \pmod{n-1}$ (see Section 1) that extends modding out by $n - 1$. If we let $\Delta_n$ denote the kernel of $\phi_n$, then the homeomorphisms used in the representation of $F_{n,\infty}$ satisfy 1–4 above and

5′. preserve $\Delta_n$ as a set.

See Section 2 for full definitions of the generalized Thompson groups.

The second fundamental fact is that the actions above are highly transitive. Because of this, a theorem of McCleary-Rubin [20] applies which allows us to identify the automorphism group of a generalized Thompson group $H$ with its normalizer $N(H)$ in the full homeomorphism group of the real line or circle. Thus we can identify an automorphism with a homeomorphism that normalizes the group in question and discuss properties of the homeomorphism (e.g., satisfies properties 1–5 above) as if they were properties of the automorphism. This practice starts with the next paragraph.

It is very easy to analyze the automorphisms of $T_{n,n-1}$ and $F_n$ that satisfy 1–5 and the automorphisms of $F_{n,\infty}$ that satisfy 1–4 and 5′. The automorphisms that are PL can also be fairly completely analyzed. Automorphisms that are not PL are labeled *exotic* in [1].

The structure of exotic automorphisms depends only on the local structures of the Thompson groups and the transitivity properties that bring in the McCleary-Rubin Theorem. Thus it suffices to study just one group with each of these sets of properties, and we choose the largest available. We thus invent "containers" where the containing group for groups with the properties of $F_n$ is $A_n(\mathbf{R})$ consisting of all self homeomorphisms of $\mathbf{R}$ or $S_r$ that satisfy 1–5, and the containing group for groups with the properties of $F_{n,\infty}$ is $B_n(\mathbf{R})$ consisting of all self homeomorphisms of $\mathbf{R}$ that satisfy 1–4 and 5′. Note that $A_2(\mathbf{R}) = B_2(\mathbf{R})$. The group $T_{n,n-1}$ can serve as its own containing group.

The generalized Thompson groups are countable while $A_n(\mathbf{R})$ and $B_n(\mathbf{R})$ are uncountable. However, the outer automorphism groups $\mathrm{Out}(A_n(\mathbf{R}))$ and $\mathrm{Out}(B_n(\mathbf{R}))$ are countable, and each of $\mathrm{Aut}(F_n)$ and $\mathrm{Aut}(F_{n,\infty})$ can be largely determined from $\mathrm{Out}(A_n(\mathbf{R}))$ and $\mathrm{Out}(B_n(\mathbf{R}))$ respectively. Conversely, we use our knowledge of the Thompson groups to help understand the outer automorphisms of the containing groups. It is a consequence of [2, Lemma 3.7] that $\mathrm{Out}(A_n(\mathbf{R}))$ and $\mathrm{Out}(T_{n,n-1})$ are isomorphic, so the role of $T_{n,n-1}$ as a containing group will be brief.

In [2] it is shown that $\mathrm{Out}(T) = \mathrm{Out}(A_2(\mathbf{R}))$ is the cyclic group of order 2 (generated by a reflection). From there it is shown that the automorphism groups of $F$ and $T$ are "small" and easy to describe. In particular, an index 2 subgroup in each of $\mathrm{Aut}(F)$ and $\mathrm{Aut}(T)$ consists entirely of elements that satisfy 1–5 above.

For $n > 2$, we show that $B_n(\mathbf{R})$ is a characteristic subgroup of index $n - 1$ in $A_n(\mathbf{R})$ and that the automorphism groups (normalizers) satsify $\mathrm{Aut}(A_n(\mathbf{R})) \subseteq \mathrm{Aut}(B_n(\mathbf{R}))$. Thus there is a subgroup $\mathrm{Out}_*(A_n(\mathbf{R}))$ of $\mathrm{Out}(B_n(\mathbf{R}))$ that surjects to $\mathrm{Out}(A_n(\mathbf{R}))$ with kernel of order $n - 1$.

For $n$ not a prime, there are PL automorphisms of $A_n(\mathbf{R})$ and $B_n(\mathbf{R})$ that are not in $A_n(\mathbf{R})$ and $B_n(\mathbf{R})$ in that slopes outside of $\langle n \rangle$ are used. We can define $\mathrm{Aut}_{PL}(A_n(\mathbf{R}))$ and $\mathrm{Aut}_{PL}(B_n(\mathbf{R}))$ as those automorphisms that are PL, with images $\mathrm{Out}_{PL}(A_n(\mathbf{R}))$ and $\mathrm{Out}_{PL}(B_n(\mathbf{R}))$ in the outer automorphism groups.

We can now indicate some of our results.

I. $\mathrm{Aut}_{PL}(A_n(\mathbf{R})) = \mathrm{Aut}_{PL}(B_n(\mathbf{R}))$ and the structure is describable modulo difficult questions in number theory (Section 3 and Theorem 6.1.4).

II. $\mathrm{Aut}(B_n(\mathbf{R}))$ is complete (has trivial center and trivial outer automorphism group so that $\mathrm{Aut}^2(B_n(\mathbf{R})) = \mathrm{Aut}(B_n(\mathbf{R}))$), but we know nothing about $\mathrm{Aut}^2(A_n(\mathbf{R}))$ (to the extent of not knowing if it is contained in $\mathrm{Aut}(B_n(\mathbf{R}))$) (Theorems 2.4.3 and 6.1.9).

III. For $n > 2$, $\mathrm{Out}_{PL}(A_n(\mathbf{R}))$ is of infinite index in $\mathrm{Out}(A_n(\mathbf{R}))$ and in particular intersects trivially an isomorphic copy of $F_2$ (Theorem 6.1.5).

IV. For $n > 2$, $\mathrm{Out}_*(A_n(\mathbf{R}))$ is of infinite index in $\mathrm{Out}(B_n(\mathbf{R}))$ and in particular intersects trivially an isomorphic copy of $F_{m,\infty}$ for each $2 \leq m < n$ (Theorems 6.1.6(d) and 6.3.2).

V. For $n \geq 4$, there is an element of order $n-2$ in $\mathrm{Out}(B_n(\mathbf{R}))$ (Theorem 6.1.10).

VI. For all $n > 2$, $F = F_2 = F_{2,\infty}$ contains isomorphic copies of $F_n$ and $F_{n,\infty}$ (Proposition 2.5.2).

Item VI is completely different from the others and its relevance(?) is touched on in the second part of the introduction. It is elementary but seems not to be known.

In IV, $\mathrm{Out}_*(A_n(\mathbf{R}))$ actually has trivial intersection with large parts of copies of $\mathrm{Aut}(F_{m,\infty})$ for $2 \leq m < n-1$.

0.2. **Background and remarks.** Thompson invented the groups $F \subseteq T \subseteq G$ in connection with studies in logic. Later, Thompson established $T$ and $G$ as the first known examples of infinite, finitely presented, simple groups [25] (see [7] for published details). Higman [18] defined an infinite family of finitely presented, infinite groups $G_{n,r}$ with $G = G_{2,1}$ and showed that $G_{n,r}$ is simple when $n$ is even, and that it contains a simple subgroup of index 2 when $n$ is odd.

Brown and Geoghegan [6] established $F$ as the first known example of a torsion free group of type $FP_\infty$ and not of type $FP$. Later Brown [5] defined infinite families of groups $F_{n,\infty} \subseteq F_n \subseteq T_{n,r}$ contained in $G_{n,r}$ with $F = F_2 = F_{2,\infty}$ and $T = T_{2,1}$ with the $F_n$ and $F_{n,\infty}$ torsion free. Brown showed that all $F_n$, $F_{n,\infty}$, $T_{n,r}$ and $G_{n,r}$ are of type $FP_\infty$ and not of type $FP$.

Thompson's groups have also appeared in homotopy theory [12], dynamical systems [14], the algebra of string rewriting [8], and diagram groups over semi-group presentations [16]. In abstract measure theory, Geoghegan has conjectured (see [4] or [13, Question 13]) that $F$ might be a counterexample to the conjecture that any finitely presented group with no non-cyclic free subgroup is amenable (admits a bounded, non-trivial, finitely additive measure on all subsets that is invariant under left multiplication). The contents and bibliographies of [2, 3, 5, 7, 8, 16, 19] give access to all that the authors know concerning Thompson's groups and their connections to other areas of mathematics.

Our results give no information on the amenability question, but they do raise interesting parallels between Thompson's groups and free groups. The free group on $n > 1$ generators has a complete automorphism group [10, 11], and the corresponding outer automorphism group is large and contains copies of the automorphism groups of the free groups on $m < n$ generators. Also, the free group on 2 generators contains copies of all countable free groups. These observations can be compared with II, IV and VI above.

The authors are rather surprised with V. We have no idea if there is an upper bound on the order of torsion elements in the outer automorphism groups.

Item III contains the main divergence with the results of [2]. It implies that for $n > 2$ there are large groups of automorphisms of the Thompson groups that do not share the local structure of the elements of the Thompson groups and in fact are not even PL. This item proves the existence of exotic automorphisms of Thompson's groups. In [1] exotic endomorphisms are found: non-PL homeomorphisms of $\mathbf{R}$ that conjugate the Thompson groups into themselves, but not onto.

Item IV shows that very closely related groups (one a finite index, characteristic subgroup of the other) can have very different automorphism groups (one is of infinite index in the other). This phenonemon is not new, but examples require some work to construct. See Section 6.

We have small things to say about other subgroups of $A_n(\mathbf{R})$ and $B_n(\mathbf{R})$. We point out that certain subgroups that capture enough of the local structure of the elements will have their automorphism groups contained in $\mathrm{Aut}(A_n(\mathbf{R}))$ and $\mathrm{Aut}(B_n(\mathbf{R}))$. Thus we discuss a class of groups associated to each of $A_n(\mathbf{R})$ and $B_n(\mathbf{R})$. We have nothing to say about $T_{n,r}$ when $r$ is not congruent to 0 modulo $n - 1$ since in this case, $\mathrm{Out}(A_n(\mathbf{R}))$ and $\mathrm{Out}(T_{n,r})$ are hard to compare.

When $n - 1$ has non-trivial divisors, there are containing groups between $A_n(\mathbf{R})$ and $B_n(\mathbf{R})$. It is clear that some of our techniques could be applied to these groups, but we do not do so here.

The results in this paper were obtained by exploiting various algebraic and geometric relationships between the various Thompson groups. For example, the algebraic Proposition 4.1.2 and geometric Lemma 4.1.5 combine to yield Theorem 4.1.6. A more purely geometric point of view is taken in [2].

The unpublished notes [1] of Bieri and Strebel contain much material about the automorphisms of groups of piecewise linear homeomorphisms of the real line. There is overlap of varying amounts between the material in [1] and material below in Sections 1, 2.6 and 3 and some material in Section 6 that ties our various results together. The parts of the current paper that do not overlap with [1] are Sections 4 and 5. The contents of the various sections are discussed briefly below. See [2, Section 1.3] for a discussion of material in [1] that complements our results. There are also some details from a very small part of [1] given in an appendix to [24].

0.3. **The structure of the paper.** We combine algebraic and geometric analysis. The algebra is easier to work with on $\mathbf{R}$ and the geometry is easier to work with on the circle. The results of [2] allow us to move back and forth between the two settings. Most of the paper takes place on the line with a few short trips to the circle to establish some key geometric facts. By and large the algebra furnishes examples and the geometry shows that the examples are non-trivial.

The point of view of [2] is adopted here that the groups we work with consist of locally defined homeomorphisms and it is the local properties that determine the automorphisms. It is for this reason that we define two classes of groups with the local properties of $F_n$ and $F_{n,\infty}$ and work with $A_n(\mathbf{R})$ and $B_n(\mathbf{R})$, the largest members of each class. The classes are introduced in Section 1 and the Thompson groups are introduced in Section 2. The piecewise linear automorphisms are analyzed in Section 3. Techniques for constructing automorphisms are developed in Section 4. The emphasis to this point is the existence of automorphisms. Section 5 develops restrictions on the automorphisms. From the results in Section 5 all "negative" results are obtained: certain automorphisms are not PL, certain automorphisms of $B_n$ are not automorhisms of $A_n$, etc. Results are gathered in Section 6, and some questions are raised in Section 7.

## 1. The classes

1.1. **The class "containers".** We define the groups that contain all the groups in their class.

For $n \geq 2$ an integer, let $\mathbf{Z}[\frac{1}{n}] = \{an^b | a, b \in \mathbf{Z}\}$ and let $\langle n \rangle = \{n^k | k \in \mathbf{Z}\}$. For a PL function $f$ from $\mathbf{R}$ to itself, a *break* of $f$ is an $x$ in the interior of the domain of $f$ at which $f'$ is discontinuous.

**Definition 1.1.1.** Let $n \geq 2$ be an integer. Let $A_n(\mathbf{R})$ be the set of all homeomorphisms $h : \mathbf{R} \to \mathbf{R}$ so that

1. $h$ is orientation preserving,
2. $h$ is PL,
3. the slopes of $h$ are all in $\langle n \rangle$,
4. the breaks of $h$ are in a discrete subset of of $\mathbf{Z}[\frac{1}{n}]$, and
5. $(\mathbf{Z}[\frac{1}{n}])h = \mathbf{Z}[\frac{1}{n}]$.

Item 1 follows from item 3, but is included for emphasis. The fact that Property 5 can be weakened to $(\mathbf{Z}[\frac{1}{n}])h \subseteq \mathbf{Z}[\frac{1}{n}]$ in the presence of Property 3 will be of no use in this paper.

**Lemma 1.1.2.** *The set $A_n(\mathbf{R})$ is a group, and so is the set of all homeomorphisms $h : \mathbf{R} \to \mathbf{R}$ that satisfy Property 1, 2, 4 and 5 of Definition 1.1.1.*

*Proof.* That Properties 1–2 are preserved under inversion and composition is standard. That Property 5 is preserved by inversion and composition is immediate. That Property 4 is preserved by inversion follows from Property 5, and that Property 4 is preserved by composition follows from the chain rule and Property 5. That Property 3 is preserved by composition and inversion follows from the chain rule. □

For $n \geq 2$ an integer, define $\phi_n : \mathbf{Z}[\frac{1}{n}] \to \mathbf{Z}_{n-1}$ by $(an^b)\phi_n = a \pmod{n-1}$, and for $f \in A_n(\mathbf{R})$, let $f\rho_n = (xf)\phi_n - (x)\phi_n$ for some $x \in \mathbf{Z}[\frac{1}{n}]$.

**Lemma 1.1.3.** *For $n \geq 2$ an integer,*

(a) *$\phi_n$ is well defined and a surjective ring homomorphism,*
(b) *if $f \in A_n(\mathbf{R})$ and $x, x' \in \mathbf{Z}[\frac{1}{n}]$, then $(xf)\phi_n - x\phi_n = (x'f)\phi_n - x'\phi_n$,*
(c) *$\rho_n : A_n(\mathbf{R}) \to \mathbf{Z}_{n-1}$ is a well defined, surjective homomorphism.*

*Proof.* (a) The function $\phi_n$ is well defined since $an^b = a'n^{b'}$ with all of $a$, $b$, $a'$, $b'$ integers implies that $a = a' = 0$ or $a/a'$ is an integral power of $n$ and we know $n \equiv 1 \pmod{n-1}$. It is a ring homomorphism for the same reason. That it is surjective is trivial.

(b) Let $x = x_0 < x_1 < x_2 < \cdots < x_n < x_{n+1} = x'$ where $x_1, \ldots, x_n$ are the breaks of $f$ between $x$ and $x'$. The lemma follows if it is true that $(x_i f)\phi_n - x_i \phi_n = (x_{i+1} f)\phi_n - x_{i+1}\phi_n$ so we assume that there are no breaks between $x$ and $x'$. We have $x'f - xf = n^k(x' - x)$, and since $\phi_n$ is a ring homomorphism taking $n$ to 1, we get $(x'f)\phi_n - (xf)\phi_n = x'\phi_n - x\phi_n$ which implies the result.

(c) That $\rho_n$ is well defined follows from (b). It is surjective since all translations by integers are in $A_n(\mathbf{R})$. Now

$$(gf)\rho_n = (xgf)\phi_n - x\phi_n$$
$$= (xgf)\phi_n - (xg)\phi_n + (xg)\phi_n - x\phi_n$$
$$= f\rho_n + g\rho_n$$

$\square$

.

For $n \geq 2$ an integer, let $\Delta_n$ be the kernel of $\phi_n : \mathbf{Z}[\frac{1}{n}] \to \mathbf{Z}_{n-1}$. For an $f \in A_n(\mathbf{R})$, we think of $f\rho_n$ as the amount that $f$ rotates the cosets of $\Delta_n$ in $\mathbf{Z}[\frac{1}{n}]$.

Let $B_n(\mathbf{R})$ be the kernel of $\rho_n : A_n(\mathbf{R}) \to \mathbf{Z}_{n-1}$. The following is immediate from Lemma 1.1.3(c).

**Lemma 1.1.4.** *For $n \geq 2$ an integer, $|A_n(\mathbf{R}) : B_n(\mathbf{R})| = n - 1$.* $\square$

For a real $r > 0$, we let $S_r = \mathbf{R}/r\mathbf{Z}$ be the circle of length $r$. If $r \in \mathbf{Z}[\frac{1}{n}]$, then $r\mathbf{Z} \subseteq \mathbf{Z}[\frac{1}{n}]$, and we will denote the quotient $\mathbf{Z}[\frac{1}{n}]/r\mathbf{Z} \subseteq S_r$ by $\mathbf{Z}[\frac{1}{n}]$. If $r \in \Delta_n$, then $r\mathbf{Z} \subseteq \Delta_n$ and $\phi_n$ induces a surjective homomorphism $\phi_n : \mathbf{Z}[\frac{1}{n}]/r\mathbf{Z} \to \mathbf{Z}_{n-1}$. From now on we will let context determine whether we are in $S_r$ or $\mathbf{R}$.

For a real $r > 0$, $r \in \mathbf{Z}[\frac{1}{n}]$, we let $A_n(S_r)$ be those self homeomorphisms of $S_r$ that satisfy the 5 properties in Definition 1.1.1. An argument similar to the proof of Lemma 1.1.3(c) shows the following.

**Lemma 1.1.5.** *For a real $r > 0$, $r \in \Delta_n$, there is a well defined homomorphism $\rho_n : A_n(S_r) \to \mathbf{Z}_{n-1}$ given by $f\rho_n = (xf)\phi_n - x\phi_n$ for any $x \in \mathbf{Z}[\frac{1}{n}]$.* $\square$

Note that in the proof of Lemmas 1.1.3(b) and 1.1.5, we do not need $f$ to be a homomorphism, and that $f\rho_n$ is well defined for any continuous self-map of $\mathbf{R}$ or $S_r$ that satisfies 1–5 in Definition 1.1.1.

For a real $r > 0$, $r \in \Delta_n$, we let $B_n(S_r)$ be the kernel of $\rho_n : A_n(S_r) \to \mathbf{Z}_{n-1}$.

## 1.2. Transitivity and the equality of automorphism and normalizer.

**Lemma 1.2.1.** *For $n \geq 2$ an integer, let $x_1 < x_2 < \cdots < x_n$ and $y_1 < y_2 < \cdots < y_n$ all be in $\mathbf{Z}[\frac{1}{n}]$ such that $x_i \phi_n = y_i \phi_n$ for all $i$ with $1 \leq i \leq n$. Then there is an $f \in B_n(\mathbf{R})$ so that $x_i f = y_i$ for each $i$ with $1 \leq i \leq n$.*

*Proof.* By putting together pieces of a function with infinite portions of slope 1, it is clear that the lemma follows if we show that there is a function taking $[x_1, x_2]$ onto $[y_1, y_2]$ satisfying all the local properties of elements of $B_n(\mathbf{R})$. We can divide the interval $[x_1, x_2]$ into subintervals of length $n^k$ for some integer $k$ and we can do the same to $[y_1, y_2]$ with intervals of length $n^l$ for some integer $l$. Since $(y_1 - x_1)\phi_n = 0 = (y_2 - x_2)\phi_n$, we have $(y_2 - y_1)\phi_n = (x_2 - x_1)\phi_n$ and thus the number of intervals used in the subdivision of $[x_1, x_2]$ is equivalent to the number of intervals in the subdivision of $[y_1, y_2]$ modulo $n - 1$. We can increase the number of intervals in a subdivision by any multiple of $n - 1$ by taking intervals in the subdivision and subdividing them (repeatedly if necessary) into $n$ intervals of equal length. This preserves the property that each interval has length an integral power of $n$. In this way we can get subdivisions of $[x_1, x_2]$ and $[y_1, y_2]$ that have the same number of intervals. Now a function can be defined from $[x_1, x_2]$ onto $[y_1, y_2]$ that takes the intervals in the subdivision of $[x_1, x_2]$ in order to those of the subdivision of $[y_1, y_2]$ so as to be affine on each interval. The properties of elements of $A_n(\mathbf{R})$ are automatically satisfied, and the extra property required to be in $B_n(\mathbf{R})$ is guaranteed by the fact that $x_1\phi_n = y_1\phi_n$. $\quad\square$

The techniques in the proof of previous lemma yield the next lemma.

**Proposition 1.2.2.** *Let* $-\infty \leq a < b \leq \infty$ *be extended integers and let* $\{x_i\}$ *and* $\{y_i\}$, $a < i < b$, *be discrete sets in* $\mathbf{Z}[\frac{1}{n}]$ *which satisfy* $x_i < x_{i+1}$ $y_i < y_{i+1}$, *and* $(y_{i+1}-x_i)\phi_n = (y_{i+1}-x_{i+1})\phi_n$ *for each* $i$ *with* $a < i < i+1 < b$. *Assume that for each* $j$ *for which* $a < 2j < 2j+1 < b$ *that there is a homeomorphism* $f_j$ *from* $[x_{2j}, x_{2j+1}]$ *onto* $[y_{2j}, y_{2j+1}]$ *satisfying the 5 properties in Definition 1.1.1. Then there is a function* $f \in A_n(\mathbf{R})$ *that agrees with* $f_j$ *for each* $j$ *with* $a < 2j < 2j + 1 < b$. $\quad\square$

Let $D$ be a subset of $\mathbf{R}$ and let a group $G$ act on $\mathbf{R}$ by order preserving homeomorphisms. We say that $G$ acts *o-k-transitively* on $D$ if for every $x_1 < \cdots < x_k$ and $y_1 < \cdots < y_k$ in $D$, there is an $f \in G$ so that $x_i f = y_i$ for each $i$ with $1 \leq i \leq k$. For a subset $D$ of $S_r$, we say that a group $G$ that acts on $S_r$ by orientation preserving homeomorphisms acts *o-k-transitively* on $D$ if for every $x_1, \ldots, x_k$ and $y_1, \ldots, y_k$ in $D$ with $x_1 < \cdots < x_k$ and $y_1 < \cdots < y_k$ in some closed interval in $S_r$ with the order inherited from the projection $\mathbf{R} \to S_r = \mathbf{R}/r\mathbf{Z}$, there is an $f \in G$ so that $x_i f = y_i$ for each $i$ with $1 \leq i \leq k$.

For an $f \in \text{Homeo}(\mathbf{R})$, we let the *support* of $f$ be the set of all $x \in \mathbf{R}$ for which $xf \neq x$. If $G$ acts on $\mathbf{R}$ by order preserving homeomorphisms, then we say that $G$ *has a non-trivial element of bounded support* if for some $f \in G$ not the identity, $f$ is fixed off some compact subset of $\mathbf{R}$. If $G$ acts on $S_r$ by orientation preserving homeomorphisms, then we say that $G$ *has a non-trivial element of bounded support* if for some $f \in G$ not the identity, $f$ is fixed on some open set.

**Lemma 1.2.3.** *The groups* $A_n(\mathbf{R})$, $A_n(S_r)$, $(r \in \mathbf{Z}[\frac{1}{n}])$, $B_n(\mathbf{R})$ *and* $B_n(S_r)$, $(r \in \Delta_n)$, *all have non-trivial elements of bounded support and act o-k-transitively on* $\mathbf{Z}[\frac{1}{n}]$ *or* $\Delta_n$ *(in the appropriate domain) for all* $k$.

*Proof.* This follows from Lemma 1.2.1 and an easily proven corresponding lemma for $S_r$. $\quad\square$

**Theorem 1.2.4** (McCleary-Rubin [20]). *Let* $G$ *act on* $\mathbf{R}$ *or* $S_r$ *by orientation preserving homeomorphisms. Assume that* $G$ *has a non-trivial element of bounded support. If* $G$ *acts on* $\mathbf{R}$, *assume that* $G$ *acts o-2-transitively on a dense subset of* $\mathbf{R}$. *If* $G$ *acts on* $S_r$, *assume that* $G$ *acts o-3-transitively on a dense subset of* $S_r$. *Then for each automorphism* $\alpha$ *of* $G$, *there is a unique self homeomorphism* $h$ *of* $\mathbf{R}$ *or* $S_r$ *(whichever is appropriate) so that* $f\alpha = h^{-1}fh$ *for every* $f \in G$.

**1.3. The two classes and their automorphisms.** Let the class $[A_n]$ consist of all subgroups $G$ of $A_n(\mathbf{R})$ or $A_n(S_r)$, $r \in \mathbf{Z}[\frac{1}{n}]$, that have a non-trivial element of bounded support, that act o-2-transitively on $\mathbf{R}$ or o-3-transitively on $S_r$, and that use all the germs of $A_n(X)$, $X = \mathbf{R}$ or $X = S_r$, in that for each $x \in X$ and each $f \in A_n(X)$, there is a $g \in G$ so that $f = g$ on some open neighborhood of $x$.

Let the class $[B_n]$ consist of all subgroups $G$ of $B_n(\mathbf{R})$ or $B_n(S_r)$, $r \in \Delta_n$, that have a non-trivial element of bounded support, that act o-2-transitively on $\mathbf{R}$ or o-3-transitively on $S_r$, and that use all the germs of $B_n(X)$, $X = \mathbf{R}$ or $X = S_r$, in that for each $x \in X$ and each $f \in B_n(X)$, there is a $g \in G$ so that $f = g$ on some open neighborhood of $x$.

We immediately get the following.

**Lemma 1.3.1.** *The groups* $A_n(\mathbf{R})$, $A_n(S_r)$, $(r \in \mathbf{Z}[\frac{1}{n}])$, *are in the class* $[A_n]$, *and the groups* $B_n(\mathbf{R})$ *and* $B_n(S_r)$, $(r \in \Delta_n)$, *are in the class* $[B_n]$. $\qquad\square$

Let Homeo$(X)$ be the group of homeomorphisms of the topological space $X$. If $G$ is a subgroup of Homeo$(X)$, then let $N(G)$ denote the normalizer of $G$ in Homeo$(X)$.

**Lemma 1.3.2.** *Let* $X = \mathbf{R}$ *or* $X = S_r$, *and let* $C = A$ *or* $C = B$. *Then the following hold.*

  (a) *For* $G$ *in the class* $[C_n]$, *taking* $h \in N(G)$ *to conjugation of* $G$ *by* $h$ *is an isomorphism from* $N(G)$ *to* Aut$(G)$.

  (b) *If* $G$ *is in* $[C_n]$, *we have* $N(G) \subseteq N(C_n(X))$.

  (c) $G \in [C_n]$ *is characteristic in* $C_n(X)$ *if and only if* $N(G) = N(C_n(X))$.

*Proof.* (a) This follows immediately from the McCleary-Rubin theorem.

(b) For $h \in N(G)$, $f \in C_n(X)$, and $x \in X$, some $g \in G$ agrees with $f$ on a neighborhood of $x$ so that $h^{-1}gh$ and $h^{-1}fh$ agree on a neighborhood of $xh$. Since $h \in N(G)$ and $G \subseteq C_n(X)$, the behavior of $h^{-1}gh$ on a neighborhood of $xh$ satisfies the requirements in Definition 1.1.1. Thus $h^{-1}fh$ is in $C_n(X)$.

(c) If $G$ is characteristic in $C_n(X)$, then restriction gives a homomorphism from $N(C_n(X))$ to $N(G)$ which is seen to be the bijective inverse of the inclusion in (b) by the uniqueness part of the McCleary-Rubin theorem. The converse follows from (a). $\qquad\square$

For $X = \mathbf{R}$ or $X = S_r$, the isomorphism between Aut$(A_n(X))$ and $N(A_n(X))$ allows us to identify Out$(A_n(X))$ with $N(A_n(X))/A_n(X)$. Similarly for $B_n(X)$. We let Out$_o(A_n(X))$ and Out$_o(B_n(X))$ be the images in Out$(A_n(X))$ and Out$(B_n(X))$ of the orientation preserving elements of $N(A_n(X))$ and $N(B_n(X))$. They are index 2 subgroups since all elements of $A_n(X)$ and $B_n(X)$ preserve orientation and $x \mapsto -x$ is a normalizer on both $\mathbf{R}$ and $S_r$. For $r \in \mathbf{R}$, let $t_r$ be the translation $xt_r = x + r$.

The next lemma follows from [2, Lemma 3.7] and its proof.

**Proposition 1.3.3.** *(a) The projection* $\mathbf{R} \to S_r = \mathbf{R}/r\mathbf{Z}$ *induces an isomorphism* Out$(B_n(\mathbf{R})) \to $ Out$(B_n(S_r))$ *for* $r \in \Delta_n$.

*(b) Each element of* Out$_o(B_n(\mathbf{R}))$ *has a representative that fixes* $t_{n-1}$. $\qquad\square$

**1.4. $B_n$ is characteristic in $A_n$.** We start with the geometric property that makes $B_n$ characteristic.

**Lemma 1.4.1.** *With* $X$ *either* $\mathbf{R}$ *or* $S_{n-1}$, $B_n(X)$ *is the subgroup of* $A_n(X)$ *generated by elements with non-empty open sets of fixed points.*

*Proof.* Elements fixed on non-empty open sets are all in the kernel of $\rho_n$ giving one containment. When $X = \mathbf{R}$, and an $\alpha \in B_n(X)$ is given, then Proposition 1.2.2 gives a $\beta \in A_n(X)$ that is the identity off some finite interval that agrees with $\alpha$ on an open set. Now $\alpha = ((\alpha\beta^{-1})\beta$ and both $\beta$ and $\alpha\beta^{-1}$ have non-empty open sets of fixed points. When $X = S_{n-1}$, an open set $U$ is chosen so small that the closure of $U \cup U\alpha$ is contained in an open interval $I$ in $S_{n-1}$ whose complement contains a non-empty open set. Now $\beta$ can be chosen with support in $I$ to agree with $\alpha$ on $U$. $\square$

**Proposition 1.4.2.** *With $X$ either $\mathbf{R}$ or $S_{n-1}$, the normalizer of $A_n(X)$ is contained in the normalizer of $B_n(X)$. This implies that $B_n(X)$ is characteristic in $A_n(X)$.*

*Proof.* Conjugation preserves the property of being fixed on an open set. The last statement comes from the fact that the normalizers are naturally identified with the automorphism groups. $\square$

1.5. **The permutation induced by a normalizer.** In this section we visit the circle to obtain information that is then lifted to the real line.

Fix an integer $n \geq 2$. Let $\Pi_n$ be the permutation group on the set $\{0, 1, \ldots, n-1\}$. Let $\mu_n : S_{n-1} \to S_{n-1}$ be the map defined by $x\mu_n = nx$. We say that a map is a *local homeomorphism* if there is an open cover of the domain so that, for each $U$ in the cover, the map restricted to $U$ is a homeomorphism onto the image of $U$.

**Lemma 1.5.1.** *The map $\mu_n$ is a local homeomorphism, is not a homeomorphism, and otherwise satisfies the 5 properties in Definition 1.1.1. There are $n-1$ fixed points of $\mu_n$ and these are the $n-1$ integer points in $S_{n-1}$. The germs of $\mu_n$ are contained in the germs of $B_n(S_{n-1})$. For $x \in S_{n-1}$ in $\mathbf{Z}[\frac{1}{n}]$, the forward orbit of $x$ under iterates of $\mu_n$ ends at the fixed point $x\phi_n$.*

*Proof.* The first sentence is immediate, and the second sentence is shown by noting that the solutions to $nx = x + k(n-1)$ in $\mathbf{R}$ are the integers. We already know that the germs of $\mu_n$ are contained in those of $A_n(S_{n-1})$. The third sentence follows because $\mu_n$ has fixed points, so by the remarks following Lemma 1.1.5, $\mu_n\rho_n = 0$. For the last sentence, the action of $\mu_n$ on the base $n$ expansion of $x$ is to shift the $n$-ary point one position to the right and then reduce the integer part modulo $(n-1)$. This preserves the value of $\phi_n$. Eventually, the characters in the fractional part are used up and an integer remains. $\square$

**Lemma 1.5.2.** *If $h$ is a normalizer of $B_n(S_{n-1})$, then $\mathbf{Z}[\frac{1}{n}]h = \mathbf{Z}[\frac{1}{n}]$.*

*Proof.* The germs at an $x \in S_{n-1}$ of elements of $B_n(S_{n-1})$ that have a fixed point at $x$ form a group. If $x \in \mathbf{Z}[\frac{1}{n}]$, then this group is isomorphic to $\mathbf{Z} \times \mathbf{Z}$ with generators being "slope $n$ to the right of $x$" and "slope $n$ to the left of $x$." For any $x \notin \mathbf{Z}[\frac{1}{n}]$, this group is either isomorphic to $\mathbf{Z}$ or is trivial. Conjugation by $h$ induces an isomorphism from the group of germs at $x$ of elements with fixed point $x$ to the group of germs at $xh$ of elements with fixed point $xh$. The result follows. $\square$

**Lemma 1.5.3.** *Let $h$ be a normalizer of $B_n(S_{n-1})$. Then the germs of both both $h\mu_n h^{-1}$ and $h^{-1}\mu_n h$ are germs of $B_n(S_{n-1})$.*

*Proof.* This follows immediately from Lemma 1.5.1. $\square$

**Proposition 1.5.4.** *There are homomorphims that we ambiguously refer to as $\pi_n$ from each of $N(B_n(S_{n-1}))$, $N(A_n(S_{n-1}))$, $\mathrm{Out}(B_n(S_{n-1}))$, $N(B_n(\mathbf{R}))$, $N(A_n(\mathbf{R}))$ and $\mathrm{Out}(B_n(\mathbf{R}))$, to $\Pi_{n-1}$. For $h$ in one of these groups, the equation $(x\phi_n)(h\pi_n) =$* $(xh)\phi_n$ holds for any $x \in \mathbf{Z}[\frac{1}{n}]$ and gives a well defined permutation on the residues $x\phi_n$, $x \in \mathbf{Z}[\frac{1}{n}]$, in that $(xh)\phi_n = (yh)\phi_n$ if and only if $x\phi_n = y\phi_n$ for any $x$ and $y$ in $\mathbf{Z}[\frac{1}{n}]$.

*Proof.* Any result for $N(B_n)$ will apply to $N(A_n)$ since normalizers of $A_n$ are normalizers of $B_n$.

We start with an $h$ in $N(B_n(S_{n-1}))$. Since the forward orbit of each $x \in \mathbf{Z}[\frac{1}{n}]$ ends in the fixed point $i = x\phi_n$, we have that the $(n-1)$ sets of pre-images of the fixed points under iterates of $\mu_n$ are exactly the $(n-1)$ classes of points in $\mathbf{Z}[\frac{1}{n}]$ with the same values of $\phi_n$. These are carried to the $(n-1)$ sets of pre-images of the fixed points under iterates of $h^{-1}\mu_n h$. However, since the germs of $h^{-1}\mu_n h$ are in the germs of $B_n(S_{n-1})$, each of these sets consists entirely of elements of $\mathbf{Z}[\frac{1}{n}]$ with a common value of $\phi_n$. Since $\mathbf{Z}[\frac{1}{n}]h = \mathbf{Z}[\frac{1}{n}]$, these sets cover all of $\mathbf{Z}[\frac{1}{n}]$ and therefore all of the $(n-1)$ values of $\phi_n$ on $\mathbf{Z}[\frac{1}{n}]$. Thus the $(n-1)$ sets of pre-images of the fixed points under iterates of $h^{-1}\mu_n h$ are exactly the $(n-1)$ classes of points in $\mathbf{Z}[\frac{1}{n}]$ with the same values of $\phi_n$. The fact that $\pi_n$ is a homomorphism follows from the fact that $i(h\pi_n) = (ih)\phi_n$. We get an induced homomorphism on $\mathrm{Out}(B_n(S_{n-1}))$ since $B_n(S_{n-1})$ is contained in the kernel of $\pi_n$.

To get $\pi_n$ defined on $N(B_n(\mathbf{R}))$, we note that Proposition 1.3.3(b) implies that any $h$ in $N(B_n(\mathbf{R}))$ induces a permutation on the values of $\phi_n$ modulo an element of $B_n(\mathbf{R})$. However, elements of $B_n(\mathbf{R})$ act trivially on the values of of $\phi_n$. Thus we get a well defined $h\pi_n$. That this is a homomorphism and that we get an induced homomorphism on $\mathrm{Out}(B_n(\mathbf{R}))$ follows as on $S_{n-1}$. $\qquad\square$

## 2. The Thompson groups

### 2.1. Definitions and presentations.

We define the groups referred to as generalized Thompson groups in the introduction, and we prove basic facts about them. The material here (through Proposition 2.4.2) can be found in [5, pp. 51–63] where the Thompson groups are defined as subgroups of automorphism groups of certain algebras and it is proven that they are isomorphic to the groups defined below. Since the basic setting, the arrangement of material and some of the proofs are different from what is found in [5], and since we will need some of the details of arguments, we give full proofs. See the references in [2] for other papers that treat all or parts of this material when $n = 2$. Our notation is largely, but not entirely consistent with that of [5].

Fix an integer $n \geq 2$. For a real $i$, let $g_i$ be the self homeomorphism of $\mathbf{R}$ defined by

$$xg_i = \begin{cases} x & x < i, \\ n(x-i)+i & i \leq x \leq i+1, \\ x+n-1 & x > i+1. \end{cases}$$

Typically, $i$ will be an integer. Recall that $xt_r = x + r$. Note that all the $g_i$ for $i \in \mathbf{Z}[\frac{1}{n}]$ are in $B_n(\mathbf{R})$ and that $t_r$ is in $A_n(\mathbf{R})$ if $r \in \mathbf{Z}[\frac{1}{n}]$ and in $B_n(\mathbf{R})$ if $r \in \Delta_n$.

The Thompson groups will be defined by fixing an $n$ and chosing various of the homeomorphisms $g_i$ and $t_r$ as generators. Parts of this paper make use of the relationships that exist between the Thompson groups for different values of $n$. When such discussions take place, it will be necessary to distinguish between the $g_i$ for various values of $n$. In those places, we will use $g_{n,i}$ to refer to the homeomorphism $g_i$ that is defined in the previous paragraph.

**Definition 2.1.1.** Fix an integer $n \geq 2$. If $Y$ is a set of elements in a group $G$, then we use $\langle Y \rangle$ to denote the subgroup of $G$ generated by the elements in $Y$. Using the elements of $A_n(\mathbf{R})$ and $B_n(\mathbf{R})$ described above we establish notation for certain subgroups of $A_n(\mathbf{R})$ and $B_n(\mathbf{R})$.

$$
\begin{aligned}
F_{n,\infty} &= \langle t_{n-1}, \ldots, g_{-1}, g_0, g_1, \ldots \rangle \\
F_n &= \langle t_1, \ldots, g_{-1}, g_0, g_1, \ldots \rangle \\
F_{n,i} &= \langle g_i, g_{i+1}, \ldots \rangle \\
F_{n,-\infty} &= \langle \ldots, g_{-1}, g_0, g_1, \ldots \rangle
\end{aligned}
$$

The group $F_{n,i}$ is not the same as the group $F_{n,r}$ in [5]. Our group $F_{n,\infty}$ corresponds more closely to the embedding of $F_{n,\infty}$ in $F_n$ in [5, Prop. 4.3] than the definition of $F_{n,\infty}$ in [5, Prop. 4.1].

The next item is established by direct checking. Shortly, we will see that the relations given there are sufficient to give presentations of the groups defined above.

**Lemma 2.1.2.** *Each of the following relations*

$$
\begin{aligned}
g_i^{-1} g_j g_i &= g_{j+n-1} && \text{for all } i < j \\
t_{n-1}^{-1} g_j t_{n-1} &= g_{j+n-1} && \text{for all } j \\
t_1^{-1} g_j t_1 &= g_{j+1} && \text{for all } j
\end{aligned}
$$

*is satisfied by the generators of $F_{n,\infty}$, $F_n$, $F_{n,i}$ and $F_{n,-\infty}$ whenever the generators in the relation are in the group in question.* $\qquad\square$

Thus the groups $F_{n,\infty}$, $F_n$ and $F_{n,i}$ are finitely generated and:

$$
\begin{aligned}
F_{n,\infty} &= \langle t_{n-1}, g_0, g_1, \ldots, g_{n-2} \rangle \\
F_n &= \langle t_1, g_0 \rangle \\
F_{n,i} &= \langle g_i, g_{i+1}, \ldots, g_{i+n-1} \rangle.
\end{aligned}
$$

The group $F_{n,-\infty}$ is not finitely generated since it is the union of

$$
F_{n,0} \subsetneq F_{n,-1} \subsetneq F_{n,-2} \subsetneq \cdots .
$$

**Definition 2.1.3.** Let $F_{n,\infty}^*$, $F_n^*$, $F_{n,i}^*$ and $F_{n,-\infty}^*$ be the groups presented with generators as given in Definition 2.1.1 of the corresponding unstarred group and all relations in Lemma 2.1.2 that apply to the generators of the group.

**Lemma 2.1.4** (Semi-normal form). *An element of $F_{n,\infty}^*$, $F_n^*$, $F_{n,i}^*$ or $F_{n,-\infty}^*$ can be written as $PN^{-1}$ where both $P$ and $N$ have the form*

$$
t_r^\epsilon g_{i_1}^{\delta_1} \cdots g_{i_k}^{\delta_k}
$$

*where $\epsilon \geq 0$ and is positive only if appropriate to the group, $r$ is $n-1$ or $1$ as appropriate to the group, each $\delta_j \geq 1$, and $i_1 < i_2 < \cdots < i_k$.*

*Proof.* The relations of Lemma 2.1.2 imply the following:

$$
\begin{aligned}
g_j^\epsilon g_i &= g_i g_{j+n-1}^\epsilon && \text{for all } i < j \\
g_i^{-1} g_j^\epsilon &= g_{j+n-1}^\epsilon g_i^{-1} && \text{for all } i < j \\
g_j^\epsilon t_{n-1} &= t_{n-1} g_{j+n-1}^\epsilon && \text{for all } j \\
t_{n-1}^{-1} g_j^\epsilon &= g_{j+n-1}^\epsilon t_{n-1}^{-1} && \text{for all } j \\
g_j^\epsilon t_1 &= t_1 g_{j+1}^\epsilon && \text{for all } j \\
t_1^{-1} g_j^\epsilon &= g_{j+1}^\epsilon t_1^{-1} && \text{for all } j
\end{aligned}
$$

where $\epsilon$ is $+1$ or $-1$ in each relation. In words, the first equation above says that positive powers of $g_i$ can move to the left of a $g_j^\epsilon$ with $j > i$ at the expense of raising $j$. There are similar wordings of the remaining equations. The achievement of the semi-normal form now follows by sliding generators to their appropriate positions. ∎

The above is not a normal form. For example, when $n = 2$, we have $g_0^2 g_2 g_0^{-2} = g_0 g_1 g_0^{-1}$. We can achieve a normal form by making suitable restrictions, but we have no need for this. A discussion for $n = 2$ is in [6] and [7].

**Proposition 2.1.5.** *In $F_{n,\infty}^*$, $F_{n,i}^*$ and $F_{n,-\infty}^*$ any non-trivial normal subgroup contains the commutator subgroup. In $F_n^*$, any non-trivial normal subgroup contains the commutators of all pairs in $\{t_{n-1}, \dots, g_{-1}, g_0, g_1, \dots\}$.*

*Proof.* Let a non-trivial element in a normal subgroup $K$ have semi-normal form $PN^{-1}$. By conjugating (if necessary), we can assume that $P$ and $N$ start with different generators. By inverting (if necessary), we can assume that $P$ starts with $t_r$ (with an appropropriate $r$) or with a $g_i$ of lower subscript than the generator at the start of $N$. Let the generator at the start of $P$ be $\gamma$ and let $P = \gamma^k P_1$ where $P_1$ does not start with a non-zero power of $\gamma$.

From this point to the end of the proof, all statements of equality, triviality or commutativity are modulo the normal subgroup $K$. Now either $\gamma^k$ is trivial, or $\gamma^k$ equals a word $w = P_1 N^{-1}$, in generators of "higher subscript" than that of $\gamma$ (all $g_i$ have "higher subscript" than that of $t_r$). We claim:

(*) These facts imply that $\gamma^k$ commutes with all $g_i$ and with $t_{n-1}$ if applicable.

We can assume that $\gamma^k$ is not trivial. For an integer $q$, let $w_q$ be $w$ with the subscripts of all its generators raised uniformly by $q(n-1)$. Conjugating $\gamma^k w^{-1} \in K$ by powers of $\gamma$, gives $\gamma^k = w_q$ for all $q \geq 0$. Any $g_i$ or $t_{n-1}$ conjugates some $w_q$ to $w_{q+1}$ for some sufficiently large $q$, and so commutes with $\gamma^k$. This proves (*).

(**) If $g_j$ has subscript higher than that of $\gamma$, then $\gamma^k$ conjugates $g_j$ to some $g_m$ with $m > j$ and we have $g_j = g_m$ by the conclusion of (*). This says that $g_j$ satisfies the hypotheses about $\gamma^k$ in (*), so $g_j$ commutes with all $g_i$ and with $t_{n-1}$ if applicable.

If $g_j$ is arbitrary and $\alpha = g_i$ for some $i < j$ or $\alpha = t_{n-1}$ if applicable, then $\alpha^q$ conjugates $g_j$ to $(g_j)_q$ (notation as above), and $\alpha^{q+1}$ conjugates $g_j$ to $(g_j)_{q+1}$. For $q$ sufficiently large, $(g_j)_q = (g_j)_{q+1}$ by (**). Thus $\alpha = \alpha^{q+1}(\alpha^q)^{-1}$ commutes with $g_j$ and the proof is complete. ∎

**Corollary 2.1.5.1.** *Each of the groups $F_{n,\infty}$, $F_n$, $F_{n,i}$ and $F_{n,-\infty}$ is presented by the "starred" presentation given in Definition 2.1.3.*

*Proof.* The homomorphism taking each formal generator to the homeomorphism it denotes yields an injection since no two of the homeomorphisms $g_i$ and $g_j$ with $i \neq j$ commute. ∎

Note that conjugating by either $t_1$ or $t_{n-1}$ induces an automorphism on each of $F_{n,\infty}$, $F_n$ and $F_{n,-\infty}$ and an endomorphism on each $F_{n,i}$. Conjugating by $t_1$ fixes $t_1$ and $t_{n-1}$ and takes each $g_i$ to $g_{i+1}$ and conjugating by $t_{n-1}$ fixes $t_1$ and $t_{n-1}$ and takes each $g_i$ to $g_{i+n-1}$. We use $\sigma_1$ to denote the automorphisms and endomorphisms induced by conjugation by $t_1$ and use $\sigma_{n-1}$ to denote the automorphisms and endomorphisms induced by conjugation by $t_{n-1}$.

The next lemma follows easily from the lemmas thus far established.

**Lemma 2.1.6.** *The groups $F_{n,\infty}$ and all $F_{n,i}$, $i \in \mathbf{Z}$, are isomorphic. The group $F_{n,\infty}$ is a semidirect product of $F_{n,-\infty}$ and $\mathbf{Z}$ where the action of one generator of $\mathbf{Z}$ on $F_{n,-\infty}$ is that of $\sigma_{n-1}$. The group $F_n$ is a semidirect product of $F_{n,-\infty}$ and $\mathbf{Z}$ where the action of one generator of $\mathbf{Z}$ on $F_{n,-\infty}$ is that of $\sigma_1$. The groups $F_{n,\infty}$ and $F_n$ are also ascending HNN extensions of each $F_{n,i}$ where the defining endomorphisms are $\sigma_{n-1}$ and $\sigma_1$ respectivley.* $\square$

2.2. **Characterizing the elements of Thompson's groups.** For each of the groups defined above, we characterize the self homeomorphisms of $\mathbf{R}$ that are its elements. To do this we will build new machinery for describing self homeomorphisms of $\mathbf{R}$. For the rest of this section, we fix $n$ at some integer greater than 1.

A *subdivision* of $\mathbf{R}$ into closed intervals is a locally finite collection of bounded, closed intervals whose union is all of $\mathbf{R}$ and whose interiors are non-empty and pairwise disjoint. Note that a subdivision of $\mathbf{R}$ can be described completely by its set of endpoints, and that the sets of endpoints of subdivisions of $\mathbf{R}$ are precisely those subsets of $\mathbf{R}$ that are closed, discrete and that have neither an upper nor a lower bound.

**Definition 2.2.1.** Let the *standard subdivision* $\mathbf{S}$ of $\mathbf{R}$ refer to the collection of closed unit intervals $[i, i+1]$, $i \in \mathbf{Z}$. The intervals in the standard subdivision $\mathbf{S}$ will be called the *standard intervals*. Given a subdivision of $\mathbf{R}$ into closed intervals, an *allowable modification* of the subdivision is to replace one of the closed intervals $[a, b]$ in the subdivision by its subdivision into the $n$ intervals $[a + (b-a)\frac{i}{n}, a + (b-a)\frac{i+1}{n}]$, $0 \le i < n$. An *allowable subdivision* of $\mathbf{R}$ is one obtained from the standard subdivision $\mathbf{S}$ of $\mathbf{R}$ by a finite number of allowable modifications. Note that all but finitely many of the intervals in an allowable subdivision of $\mathbf{R}$ will be standard intervals. An *allowable interval* is an interval in some allowable subdivision of $\mathbf{R}$. It is immediate that the allowable intervals are those intervals of the form $[in^k, (i+1)n^k]$ for some integer $i$ and non-positive integer $k$. We say that an allowable subdivision of $\mathbf{R}$ is *supported in a subset $A$ of $\mathbf{R}$* if $A$ contains all the non-standard intervals of the subdivision.

Let $D$ and $R$ be two allowable subdivisions of $\mathbf{R}$. An *isomorphism from $D$ to $R$* is an orientation preserving homeomorphism from $\mathbf{R}$ to itself that is the identity near $-\infty$ and that takes each interval in $D$ affinely onto an interval of $R$. There is only one isomorphism from a given $D$ to given $R$. The definition of isomorphism could be made more flexible, the next lemma could cover more groups than $F_{n,0}$, and Definition 2.2.3 extended to cover these changes, but we will have no need of this.

**Lemma 2.2.2.** *The group $F_{n,0}$ consists of all isomorphisms between pairs of allowable subdivisions of $\mathbf{R}$ that are supported in $[0, \infty)$.*

We will give a constructive proof since we will have later need of the construction. We first give a lemma that will go into the proof and whose statement we will need later.

**Definition 2.2.3.** We call a word $P$ in the generators $g_i$, $i \ge 0$, of $F_{n,0}$ as given in Definition 2.1.1 a *positive word* if all generators appear with positive exponent. Note that each of the parts $P$ and $N$ of a semi-normal form $PN^{-1}$ of a word in the generators of $F_{n,0}$ is a positive word. We will define the subdvision of $\mathbf{R}$ that we associate with such a positive word $P$. Let $P = g_{i_1} g_{i_2} \cdots g_{i_k}$. Let the prefixes of $P$ be $P_0 = 1, P_1 = g_{i_1}, \ldots, P_{k-1}, P_k = P$ where $P_j = P_{j-1} g_{i_j}$. We inductively associate subdivisions $D_0, \ldots, D_k$ to the prefixes where $D_0 = \mathbf{S}$. To do so we need a numbering scheme. Note that any allowable subdivision of $\mathbf{R}$ has 0 as an endpoint of an interval. We number the intervals of an allowable subdivision of $\mathbf{R}$ that is supported in $[0, \infty)$ so that interval 0 has 0 as its left endpoint and the intervals are numbered by the integers consecutively from left to right in $\mathbf{R}$. We now define $D_j$ as that subdivision obtained from $D_{j-1}$ by subdividing interval numbered $i_j$ in $D_{j-1}$ into $n$ subintervals of equal size (here $i_j$ is the subscript of the ending generator in $P_j$). The subdivision associated with each $P_j$ is $D_j$. In $F_{n,0}$, the subscripts of all the $g_i$ are non-negative, so all of the subdivisions $D_j$ will be supported in $[0, \infty)$.

**Lemma 2.2.4.** *(a) Every allowable subdivision of $\mathbf{R}$ that is supported in $[0, \infty)$ is associated with some positive word in semi-normal form in the generators of $F_{n,0}$.*
*(b) If $P$ is a positive word in the generators of $F_{n,0}$, then $P$ represents an isomorphism from the subdivision $D$ associated with $P$ to the standard subdivision $\mathbf{S}$ of $\mathbf{R}$.*
*(c) If $PN^{-1}$ is a semi-normal form for an element of $F_{n,0}$, then it represents an isomorphism from the allowable subdivision $D$ of $\mathbf{R}$ supported on $[0, \infty)$ and associated with $P$ to the allowable subdivision $R$ of $\mathbf{R}$ supported on $[0, \infty)$ and associated with $N$.*

*Proof.* We get (a) from the definitions by always choosing the leftmost interval to subdivide when there is an ambiguity as to which interval to subdivide next. We get (b) by induction on the length of $P$ and the observation that (b) is true when $P$ is a single generator $g_i$, and we get (c) from (b). $\quad\blacksquare$

*Proof of Lemma 2.2.2.* Both containments of the claimed equality follow from the previous lemma. $\quad\blacksquare$

**Proposition 2.2.5.** *The group $F_{n,0}$ consists of all elements $f$ of $A_n(\mathbf{R})$ that are the identity to the left of 0 and are translations by integral multiples of $n-1$ near $\infty$ in that there are integers $i$ and $M$ so that $xf = x + i(n-1)$ for all $x > M$.*

*Proof.* The generators of $F_{n,0}$ all fit the description in the statement and the description is preserved under inverses and composition. Thus it remains to show that all elements as described are in $F_{n,0}$.

Let $f$ be as described in the statement. By the previous lemma, we must show that $f$ is an isomorphism between two allowable subdivisions $P$ and $Q$ of $\mathbf{R}$ that are supported in $[0, \infty)$. Near $\pm\infty$, the standard subdivision supports $f$ as an isomorphism. We must arrange things on a compact set. There is an allowable subdivision $P$ of $S$ so that $f$ is linear on each interval in $P$. This is done by subdividing intervals until the breaks of $f$ are contained in the endpoints of the subdivision. Since $f$ is the identity to the left of 0, only intervals in $[0, \infty)$ will have to be subdivided. The image subdivision $Pf$ has its endpoints in $\mathbf{Z}[\frac{1}{n}]$, and it can be subdivided to a subdivision $Q'$ by allowable modifications so that $Q'$ is equal to $S$ except on some interval $[a, b]$, $0 \le a < b$, $a, b \in \mathbf{Z}$, where $Q'$ consists of intervals of one constant length $n^{-k}$. Let $P'$ be the subdivision $Q'f^{-1}$ of $P$. Since $f$ is linear with slope an integral power of $n$ on each interval of $P$, this subdivision of $P$ uniformly subdivides each interval of $P$ into subintervals of length $n^{-j}$ for an integer $j$ that depends on the interval. It follows that $P'$ is an allowable subdivision of $P$ and $f$ is an isomorphism carried by the allowable pair $(P', Q')$. $\quad\blacksquare$

In the proof above, we did not use the fact that the translation near $\infty$ was by an integral multiple of $n-1$ because we did not have to. Any element of $A_n(\mathbf{R})$ that fixes 0 and is a translation by an integer near $\infty$ is a translation by an integral multiple of $n-1$ near $\infty$.

The next proposition follows easily from the proposition just proven.

**Proposition 2.2.6.** *The group $F_{n,-\infty}$ consists of all elements $f$ of $A_n(\mathbf{R})$ that are the identity near $-\infty$ and that are translations by integral multiples of $n-1$ near $\infty$. The group $F_{n,\infty}$ consists of all elements $f$ of $A_n(\mathbf{R})$ that are translations by integral multiples of $n-1$ near $\pm\infty$ in that there are integers $i$, $j$ and $M$ so that $xf = x + i(n-1)$ for all $x < -|M|$ and $xf = x + j(n-1)$ for all $x > |M|$. The group $F_{n,i}$ consists of all elements $f$ of $F_{n,\infty}$ that are the identity on $(-\infty, i)$. The group $F_n$ consists of all elements of $A_n(\mathbf{R})$ that are translations by integers near $\pm\infty$ in that there are integers $i$, $j$ and $M$ so that $xf = x + i$ for all $x < -|M|$ and $xf = x + j$ for all $x > |M|$.* ∎

Since $F_n \subseteq A_n(\mathbf{R})$, each $f \in F_n$ rotates the cosets of $\Delta_n$ by $\rho_n(f)$. Thus the values of $i$ and $j$ in the description of elements of $F_n$ in the previous statement must be congruent modulo $n-1$. This observation and the preceeding proposition yield:

**Corollary 2.2.6.1.** $|F_n : F_{n,\infty}| = n-1$. ∎

**Lemma 2.2.7.** *The group $F_n$ is in $[A_n]$, and the groups $F_{n,\infty}$, $F_{n,-\infty}$ and $F_{n,i}$ are in $[B_n]$.*

*Proof.* If $C$ is any compact subset of $\mathbf{R}$ and $f \in A_n(\mathbf{R})$, then Proposition 1.2.2 implies that there is a $g$ in $A_n(\mathbf{R})$ that agrees with $f$ on $C$ and is a translation by $f\phi_n$ near $\pm\infty$. This suffices to put $g$ in $F_n$, and suffices to prove that $F_n$ is in $[A_n]$. Similar arguments give the rest of the statement. ∎

**2.3. Thompson's group on a closed interval.** Let $F_{n,[0,n-1]}$ be the set of elements of $A_n(\mathbf{R})$ that have support in $[0, n-1]$. By Proposition 2.2.6, this is a subgroup of $F_{n,0}$ and thus of $B_n(\mathbf{R})$.

**Lemma 2.3.1.** *There is a homeomorphism from $[0, \infty)$ to $[0, n-1)$ that induces, by conjugation, an isomorphism from $F_{n,0}$ to $F_{n,[0,n-1]}$.*

*Proof.* Let $D_0$ be the subdivision of $[0, n-1]$ that consists of the $n-1$ standard unit intervals $[i, i+1]$, $0 \le i < n-1$, in $[0, n-1]$. Let $D_i$ be obtained from the subdivision $D_{i-1}$ by subdividing the rightmost interval of $D_{i-1}$ into $n$ equal intervals. Let $\mathbf{S}'$ be the limit of this process. This is a subdivision of the interval $[0, n-1)$ into countably many allowable intervals. If we define allowable subdivisions of $\mathbf{S}'$ as the result of finite sequences of allowable modifications of $\mathbf{S}'$, then we can define isomorphisms between pairs of allowable subdivisions of $\mathbf{S}'$.

We note that an allowable modification raises the number of intervals in given closed interval by $n-1$. Thus an isomorphism between allowable subdivisions of $\mathbf{S}'$ results in a shift of the intervals of $\mathbf{S}'$ near $n-1$ by some integral multiple $k$ of $n-1$. It is easy to show that shuch a shift is an affine function in a neighborhood of $n-1$ with fixed point $n-1$ and slope $n^{-k}$. Thus every such isomorphism represents a function in $F_{n,[0,n-1]}$. An argument identical to that of Proposition 2.2.5 shows that every element of $F_{n,[0,n-1]}$ is such an isomorphism. The claimed conjugating homeomorphism in the statement of the lemma is the isomorphism taking the intervals in the standard subdivision **S** that are in $[0, \infty)$ to the intervals in **S'**. $\qquad \square$

**Corollary 2.3.1.1.** *The subgroup $B_n^0(S_{n-1})$ of $B_n(S_{n-1})$ that consists of all elements that fix 0 is isomorphic to $F_{n,[0,n-1]}$, to $F_{n,0}$ and also to $F_{n,\infty}$.*

*Proof.* The subgroup is clearly isomorphic to $F_{n,[0,n-1]}$, and we use Lemmas 2.3.1 and 2.1.6. $\qquad \square$

2.4. **A simple subgroup.** We define two more subgroups by properties of their elements. The subgroups are not conveniently given as presentations.

Let $F_n^0$ consist of those elements of $F_n$ that are the identity near $\pm\infty$.

**Lemma 2.4.1.** *We have $F_n^0 \subseteq F_{n,-\infty} \subseteq F_{n,\infty}$. The group $F_n^0$ consists of all $f \in F_{n,-\infty}$ whose every expression as a word $w$ in the generators $g_j$ has the exponent sum of $w$ over all $g_i$ equal to zero.*

*Proof.* The containments are clear. The second claim follows because each generator $g_j$ behaves as translation by $n-1$ near $+\infty$. $\qquad \square$

Let $F_n^s$ consist of those elements $f$ of $F_{n,-\infty}$ so that an expression of $f$ as a word $w$ in the generators $g_j$ satisfies that for each $i$ with $0 \le i < n-1$, the exponent sum of $w$ over all generators $g_j$ with $j \equiv i \pmod{n-1}$ is zero. Note that this is well defined since the relations in $F_{n,-\infty}$ preserve such sums.

**Proposition 2.4.2.** *The following are true:*

(a) *$F_n^s \subseteq F_n^0 \subseteq F_{n,-\infty} \subseteq F_{n,\infty}$.*
(b) *If $g \in B_n(\mathbf{R})$ and a compact $C \subseteq \mathbf{R}$ are given, then there is an $f \in F_n^s$ that agrees with $g$ on $C$.*
(c) *$F_n^s$ and $F_n^0$ are in $[B_n]$.*
(d) *$F_n^0$ is characteristic in $F_{n,-\infty}$ and $B_n(\mathbf{R})$.*
(e) *$F_n^s$ is the commutator subgroup of $F_{n,-\infty}$ and $F_n^0$.*
(f) *$F_n^s$ is characteristic in $F_n^0$, $F_{n,-\infty}$ and $B_n(\mathbf{R})$.*
(g) *$F_n^s$ is simple.*

*Proof.* (a) This is immediate.

(b) As mentioned in the proof of Lemma 2.2.7, it is easy to find an $f \in F_{n,-\infty}$ that equals $g$ on $C$. To get $f$ into $F_n^s$ one composes $f$ on the right by powers of various generators $g_j$. Since one can choose the subscripts $j$ up to integral multiples of $n-1$, the subscripts can be chosen so large that each new $g_j$ is the identity on $C$.

(c) This follows from (a) and (b).

(d) From the McCleary-Rubin Theorem, every automorphism of $B_n(\mathbf{R})$ and every automorphism of $F_{n,-\infty}$ is realized as conjugation by a self homeomorpism of $\mathbf{R}$. Such a conjugation preserves the property of being fixed off some compact set.

(e) The relations of $F_{n,-\infty}$ establish that $g_i$ is conjugate to $g_j$ if and only if $i \equiv j \pmod{n-1}$. Modulo the commutator subgroup, there are no other relations. Thus the abelianization of $F_{n,-\infty}$ is the sum of $n-1$ copies of $\mathbf{Z}$ and $F_n^s$ is the kernel of the abelianization homomorphism. Thus $F_n^s$ is the commutator subgroup of $F_{n,-\infty}$.

Since $F_n^0$ is characteristic in $F_{n,-\infty}$, the commutator subgroup of $F_n^0$ is normal in $F_{n,-\infty}$ and it lies in $F_n^s$ the commutator subgroup of $F_{n,-\infty}$. By Proposition 2.1.5, the commutator subgroup of $F_n^0$ is all of $F_n^s$ if the commutator subgroup of $F_n^0$ is non-trivial. However this follows from (b) since $F_n^s \subseteq F_n^0$.

(f) This follows from (e) and (d).

(g) We can use a simplification of an argument of Higman [17]. We must show that every non-trivial normal subgroup of $F_n^s$ contains all of $F_n^s$ the commutator subgroup of $F_n^0$. That is, we must show that the commutator of any two elements of $F_n^0$ lies in the normal closure in $F_n^s$ of any non-trivial element of $F_n^s$.

Let $f$ be non-trivial in $F_n^s$ and let $g$ and $h$ be in $F_n^0$. There is an open set $U \subseteq \mathbf{R}$ so that $Uf \cap U = \emptyset$. Since $g$ and $h$ are in $F_n^0$, the union of their supports lies in some compact set $C \subseteq \mathbf{R}$. Because $F_n^s$ is in $[B_n]$, it has enough transitivity properties so that there is a $j \in F_n^s$ so that $C \subseteq Uj$. Now $j^{-1}fj$ carries $Uj$ to $Ufj$ which is disjoint from $Uj$. Thus $j^{-1}fj$ carries $C \subseteq Uj$ to a set disjoint from $C$. The support of $h$ conjugated by $j^{-1}fj$ lies in $Cj^{-1}fj$ and is disjoint from the support of $g$. Thus $h$ conjugated by $j^{-1}fj$ commutes with $g$, and modulo the normal closure in $F_n^s$ of $f$, we have that $g$ and $h$ commute. $\qquad\square$

**Theorem 2.4.3.** *The automorphism group of* $\mathrm{Aut}(B_n(\mathbf{R}))$ *equals* $\mathrm{Aut}(B_n(\mathbf{R}))$.

*Proof.* By Lemmas 1.3.2(c) and 2.4.2(f), $\mathrm{Aut}(B_n(\mathbf{R})) = \mathrm{Aut}(F_n^s)$, and since $F_n^s$ is simple, [23, Theorem 13.5.10] or [22, Stmt. 6.2.1.5] give that $\mathrm{Aut}(F_n^s)$ is its own automorphism group. $\qquad\square$

2.5. **Monomorphisms from universal properties.** This short section recalls a universal property of the Thompson groups and uses it to show the existence of monomorphisms between the Thompson groups defined for different values of $n$, including monomorphisms from any Thompson group acting on $\mathbf{R}$ to $F_{2,0}$. In Section 4, we will derive different monomorphisms between Thompson groups that allow for extensions of automorphisms. The monomorphisms in this section do not seem to have such extension properties.

The endomorphism $\sigma_1 : F_{n,0} \to F_{n,0}$ taking each $g_i$ to $g_{i+1}$ has the property that $\sigma_1^n \neq \sigma_1$ but $\sigma_1^n = \sigma_1 T_0$ where $T_0 : F_{n,0} \to F_{n,0}$ is the inner automorphism gotten by conjugating by the generator $g_{n,0}$. The next lemma shows that $F_{n,0}$ is universal in this behavior. This is a trivial extension of the $n = 2$ case [12, p. 95] and [6, Prop. 1.1].

**Lemma 2.5.1.** *Let the group $G$ admit an endomorphism* $\sigma : G \to G$ *with the property that* $\sigma^n \neq \sigma$ *but* $\sigma^n = \sigma T_a$ *for some non-identity* $a \in G$ *where* $T_a : G \to G$ *is the inner automorphism gotten by conjugating by* $a$. *Then there is a homomorphism* $\kappa : F_{n,0} \to G$ *with* $g_{n,i}\kappa = a\sigma^i$, *for each* $i \geq 0$. *Further, if $a$ fails to commute with* $a\sigma$, *then $\kappa$ is a monomorphism.*

*Proof.* Note that for $0 \leq j < k$ we have

$$
\begin{aligned}
\left(a\sigma^j\right)^{-1}\left(a\sigma^k\right)\left(a\sigma^j\right) &= \left(a\sigma^j\right)^{-1}\left(a\sigma^{k-j}\sigma^j\right)\left(a\sigma^j\right) \\
&= \left[a^{-1}(a\sigma^{k-j})a\right]\sigma^j \\
&= \left[a^{-1}\left((a\sigma^{k-j-1})\sigma\right)a\right]\sigma^j \\
&= \left[(a\sigma^{k-j-1})\sigma^n\right]\sigma^j \\
&= a\sigma^{k+n-1}.
\end{aligned}
$$

Thus the images of $a$ under non-negative powers of $\sigma$ satisfy the relations of $F_{n,0}$. If $a$ does not commute with $a\sigma$, then $\kappa$ is a monomorphism by Proposition 2.1.5. $\qquad\square$

**Proposition 2.5.2.** *Let* $2 \leq m < n$ *satisfy* $(m-1)|(n-1)$ *with* $(n-1) = (m-1)d$. *Then there is a monomorphism* $\lambda_{n,m} : F_{n,0} \to F_{m,0}$ *with* $g_{n,i}\lambda_{n,m} = g_{m,i}^d$ *for each* $i \geq 0$ *and there is a monomorphism* $\lambda'_{m,n} : F_{m,0} \to F_{n,0}$ *with* $g_{m,i}\lambda'_{m,n} = g_{n,di}$ *for each* $i \geq 0$.

*Proof.* We consider $\lambda_{n,m}$. If $T_{m,0}$ denotes conjugation by $g_{m,0}$, then we have $\sigma_1 T_{m,0} = \sigma_1^m$. The hypotheses on $m$ and $n$ now imply that $\sigma_1 T_{m,0}^d = \sigma_1^n$. But $T_{m,0}^d$ is just conjugation by $g_{m,0}^d$ and the images of $g_{m,0}^d$ under powers of $\sigma_1$ are just the elements $g_{m,i}^d$. Since $g_{m,0}^d$ and $g_{m,1}^d$ do not commute, the properties of $\lambda_{n,m}$ now follow from the previous lemma.

For $\lambda'_{m,n}$, let $\Sigma = \sigma_1^d$. Now

$$\Sigma T_{n,0} = \sigma_1^d T_{n,0} = \sigma_1^{d+n-1} = \sigma_1^{d+d(m-1)} = \sigma_1^{dm} = \Sigma^m.$$

The properties of $\lambda'_{m,n}$ follow from the previous lemma and the fact that $g_{n,0}$ and $g_{n,d} = g_{n,0}\Sigma$ do not commute. $\qquad\blacksquare$

Note that the hypothesis in the previous proposition is satisfied for all $n \geq 2$ if $m = 2$. Further, if we send $t_1$ to $g_{2,0}$ and, for each $i \geq 0$, $g_{n,i}$ to $g_{2,i+1}^{n-1}$, we get a monomorphism from $F_n$ into $F_{2,0}$. This implies that there are monomorphisms into $F_{2,0}$ from each of the groups in Definition 2.1.1 and Section 2.3.

2.6. **Applying the Thompson groups to automorphisms.** We are interested in $\mathrm{Out}_o(A_n(\mathbf{R}))$ and $\mathrm{Out}_o(B_n(\mathbf{R}))$ as defined just before Proposition 1.3.3. We will show that these groups can be obtained from the Thompson groups and their automorphisms and that to a large extent, the automorphism groups of the Thompson groups can be obtained from $\mathrm{Out}_o(A_n(\mathbf{R}))$ and $\mathrm{Out}_o(B_n(\mathbf{R}))$.

Because of Proposition 1.3.3, it is convenient to look at automorphisms that commute with certain translations. We will use the notation of the next paragraph extensively. When we do, $g$ will always be a translation.

Let $G$ be a group of self homeomorphisms of $\mathbf{R}$ to which the McCleary-Rubin theorem applies and identify $\mathrm{Aut}(G)$ with the normalizer of $G$ in the full group of self homeomorphisms of $\mathbf{R}$. If $g$ is a self homeomorphism of $\mathbf{R}$, with $g$ not necessarily in $G$, then let $\mathrm{Aut}(G, g)$ be the set of automorphisms of $G$ that fix $g$, and let $C(G, g)$ be the set of inner automorphisms of $G$ that fix $g$ (and so is the centralizer of $g$ in $G$). We refer to $\mathrm{Aut}(G, g)$ as the *automorphism group of the pair* $(G, g)$. We use $\mathrm{Aut}^0(G, g)$ and $C^0(G, g)$ to denote the subgroups of the respective groups that consist of the elements that fix 0.

The automorphisms of $B_n(\mathbf{R})$ are somewhat better behaved than those of $A_n(\mathbf{R})$, and so we start with $B_n(\mathbf{R})$. The following is a definition:

$$1 \to B_n(\mathbf{R}) \to \mathrm{Aut}_o(B_n(\mathbf{R})) \to \mathrm{Out}_o(B_n(\mathbf{R})) \to 1.$$

From Proposition 1.3.3, we know that every element of $\mathrm{Out}_o(B_n(\mathbf{R}))$ has a representative the commutes with $t_{n-1}$. Automorphisms (normalizers) that commute with $t_{n-1}$ preserve the defining property of $F_{n,\infty}$ within $B_n(\mathbf{R})$. Conversely, Lemma 2.2.7 says that automorphisms of $F_{n,\infty}$ are automorphisms of $B_n(\mathbf{R})$. Thus we obtain

$$1 \to \mathrm{Aut}_B(F_{n,\infty}) \to \mathrm{Aut}_o(F_{n,\infty}) \to \mathrm{Out}_o(B_n(\mathbf{R})) \to 1$$

where $\mathrm{Aut}_B(F_{n,\infty}) = B_n(\mathbf{R}) \cap \mathrm{Aut}_o(F_{n,\infty})$. If we push further the observation that elements of $\mathrm{Out}_o(B_n(\mathbf{R}))$ have representatives the commute with $t_{n-1}$, we get

$$1 \to C(B_n(\mathbf{R}), t_{n-1}) \to \mathrm{Aut}(F_{n,\infty}, t_{n-1}) \to \mathrm{Out}_o(B_n(\mathbf{R})) \to 1.$$

We now pass to a normal subgroup of finite index. We define $\mathrm{Aut}_\pi(F_{n,\infty}, t_{n-1})$ and $\mathrm{Out}_\pi(B_n(\mathbf{R}))$ to be the kernels of the homomorphism $\pi_n$ applied to $\mathrm{Aut}(F_{n,\infty}, t_{n-1})$ and $\mathrm{Out}_o(B_n(\mathbf{R}))$ repsectively. We will see in Theorem 4.4.3 that $(\mathrm{Out}_o(B_n(\mathbf{R})))\pi_n$ is all of $\Pi_{n-1}$ which will show that the index of these subgroups is $(n-1)!$. We have

$$1 \to C(B_n(\mathbf{R}), t_{n-1}) \to \mathrm{Aut}_\pi(F_{n,\infty}, t_{n-1}) \to \mathrm{Out}_\pi(B_n(\mathbf{R})) \to 1.$$

We can make a further reduction. Given $h \in \mathrm{Aut}_\pi(F_{n,\infty}, t_{n-1})$, we have $b\phi_n = (0h)\phi_n = 0(h\pi_n) = 0$. Thus $t_b$ is in $B_n(\mathbf{R})$, and $h' = ht_b^{-1}$ fixes 0 and represents the same element in $\mathrm{Out}_o(B_n(\mathbf{R}))$ as $h$. We now have

$$1 \to C^0(B_n(\mathbf{R}), t_{n-1}) \to \mathrm{Aut}_\pi^0(F_{n,\infty}, t_{n-1}) \to \mathrm{Out}_\pi(B_n(\mathbf{R})) \to 1.$$

The order in the above sequence of short exact sequences is largely forced. We cannot restrict to normalizers that fix 0 before we restrict to normalizers that have $\pi_n$ trivial. This is because without the triviality of $\pi_n$, the translation $t_b$ used above might not be in $B_n(\mathbf{R})$.

The automorphisms of $A_n(\mathbf{R})$ are better behaved in that they contain more translations, and the restriction to normalizers that fix 0 is easier to achieve. However, Proposition 1.3.3 does not apply to $\mathrm{Out}_o(A_n(\mathbf{R}))$. We will get around this by trivializing permutations once again, but here we are not dealing with a normal subgroup of $\mathrm{Out}_o(A_n(\mathbf{R}))$ since $\pi_n$ is not well defined on $\mathrm{Out}_o(A_n(\mathbf{R}))$.

We start with the defining sequence

$$1 \to A_n(\mathbf{R}) \to \mathrm{Aut}_o(A_n(\mathbf{R})) \to \mathrm{Out}_o(A_n(\mathbf{R})) \to 1.$$

Given $h \in \mathrm{Aut}_o(A_n(\mathbf{R}))$, we have $b = 0h \in \mathbf{Z}[\frac{1}{n}]$. Thus $t_b$ is in $A_n(\mathbf{R})$, and $h' = ht_b^{-1}$ fixes 0 and represents the same element in $\mathrm{Out}_o(A_n(\mathbf{R}))$ as $h$. We have

$$1 \to A_n^0(\mathbf{R}) \to \mathrm{Aut}_o^0(A_n(\mathbf{R})) \to \mathrm{Out}_o(A_n(\mathbf{R})) \to 1.$$

However, elements of $A_n(\mathbf{R})$ that fix 0 are in $B_n(\mathbf{R})$, so we get

$$1 \to B_n^0(\mathbf{R}) \to \mathrm{Aut}_o^0(A_n(\mathbf{R})) \to \mathrm{Out}_o(A_n(\mathbf{R})) \to 1.$$

With the kernel in the above sequence containing only elements with trivial $\pi_n$, we can pass to a subgroup. We let $\mathrm{Aut}_\pi^0(A_n(\mathbf{R}))$ be the kernel of $\pi_n$ restricted to $\mathrm{Aut}^0(A_n(\mathbf{R}))$. We let $\mathrm{Out}_\pi(A_n(\mathbf{R}))$ be the image of $\mathrm{Aut}_\pi^0(A_n(\mathbf{R}))$ in $\mathrm{Out}_o(A_n(\mathbf{R}))$. We will see in Theorem 6.1.8 that the indexes are not $(n-1)!$ and are not completely known. We get

$$1 \to B_n^0(\mathbf{R}) \to \mathrm{Aut}_\pi^0(A_n(\mathbf{R})) \to \mathrm{Out}_\pi(A_n(\mathbf{R})) \to 1.$$

We now observe that a minor change to the proof of [2, Lemma 3.7] gives the following. In the change, we can assume that the normalizer $h$ fixes 0 and our extra assumption says that $(1h)\phi_n = 1$.

**Proposition 2.6.1.** *Each element of* $\mathrm{Out}_\pi(A_n(\mathbf{R}))$ *has a representative that fixes* $t_1$. $\qquad\square$

Since automorphisms (normalizers) that commute with $t_1$ preserve the defining property of $F_n$ within $A_n(\mathbf{R})$ and automorphisms of $F_n$ are automorphisms of $A_n(\mathbf{R})$, we can now bring in a Thompson group:

$$1 \to \mathrm{Aut}_B^0(F_n) \to \mathrm{Aut}_\pi^0(F_n) \to \mathrm{Out}_\pi(A_n(\mathbf{R})) \to 1$$

where $\mathrm{Aut}_B^0(F_n) = B_n^0(\mathbf{R}) \cap \mathrm{Aut}(F_n)$. As with $B_n(\mathbf{R})$ above, we can push our use of Proposition 2.6.1 farther and get

$$1 \to C^0(B_n(\mathbf{R}), t_1) \to \mathrm{Aut}_\pi^0(F_n, t_1) \to \mathrm{Out}_\pi(A_n(\mathbf{R})) \to 1.$$

The following summarizes how we can relate automorphism groups of Thompson groups to $\mathrm{Out}(A_n(\mathbf{R}))$ and $\mathrm{Out}(B_n(\mathbf{R}))$. We have been somewhat more successful with $B_n(\mathbf{R})$ than with $A_n(\mathbf{R})$. More information about the kernels $\mathrm{Aut}_B(F_{n,\infty})$ and $\mathrm{Aut}_B^0(F_n)$ will be given in Proposition 6.1.3.

**Proposition 2.6.2.** *There are short exact sequences and isomorphisms*

$$1 \to \mathrm{Aut}_B(F_{n,\infty}) \to \mathrm{Aut}_o(F_{n,\infty}) \to \mathrm{Out}_o(B_n(\mathbf{R})) \to 1,$$

$$1 \to \mathrm{Aut}_B^0(F_n) \to \mathrm{Aut}_\pi^0(F_n) \to \mathrm{Out}_\pi(A_n(\mathbf{R})) \to 1,$$

$$\mathrm{Out}_\pi(B_n(\mathbf{R})) \simeq \mathrm{Aut}_\pi^0(F_{n,\infty}, t_{n-1}) \Big/ C^0(B_n(\mathbf{R}), t_{n-1}), \text{ and}$$

$$\mathrm{Out}_\pi(A_n(\mathbf{R})) \simeq \mathrm{Aut}_\pi^0(F_n, t_1) \Big/ C^0(B_n(\mathbf{R}), t_1).$$

## 3. PL automorphisms

### 3.1. Inner automorhpisms.

**Definition 3.1.1.** We concentrate on the denominator $C^0(B_n(\mathbf{R}), t_{n-1})$ in Proposition 2.6.2. The group $C^0(B_n(\mathbf{R}), t_{n-1})$ is clearly isomorphic to the group $B_n^0(S_{n-1})$ and by Corollary 2.3.1.1 to $F_{n,[0,n-1]}$. The group $F_{n,[0,n-1]}$ is a subgroup of $F_{n,0}$ so each $W \in F_{n,[0,n-1]}$ is a word in generators $g_i$, $i \geq 0$. Given $W \in F_{n,[0,n-1]}$ we can build an element of $C^0(B_n(\mathbf{R}), t_{n-1})$ from $W$ as the infinite composition (written as a product)

$$W\Phi_n = \prod_{k=-\infty}^\infty t_{n-1}^{-k} W t_{n-1}^k = \prod_{k=-\infty}^\infty W\sigma_{n-1}^k.$$

Note that the infinite composition makes sense (and the order is irrelevant) since the supports of $W\sigma_{n-1}^k$ and $W\sigma_{n-1}^j$ are disjoint if $k \neq j$.

**Lemma 3.1.2.** *Let $n \geq 2$.*

*(a) $\Phi_n : F_{n,[0,n-1]} \to C^0(B_n(\mathbf{R}), t_{n-1})$ is an isomorphism.*

*(b) For $0 \leq j < n-1$, $W \in F_{n,[0,n-1]}$ and $w = W\Phi_n$, we get $w^{-1} g_j w = W^{-1} g_j W W_1$ where $W_1 = W\sigma_{n-1}$.*

*Proof.* (a) is clear.

Let $W_k = W\sigma_{n-1}^k$. To see (b), note first that all $W_k$ commute with $g_j$ when $k < 0$. Note next, that as a word in the generators $g_i$, $W_k$ for $k \geq 2$ uses only those $g_i$ with $i \geq 2n-2$. These values of $i$ exceed $n-2$ by at least $n$, and for $k \geq 2$, we get $W_k^{-1} g_j = g_j W_{k+1}^{-1}$. Since $W_k$ and $W_l$ commute when $k \neq l$, (b) follows. $\quad\square$

### 3.2. General PL automorphisms.

We will show that $\mathrm{Aut}_{PL}(B_n(\mathbf{R}))$ has a structure similar to that of $B_n(\mathbf{R})$, but that a larger set of slopes is used. We start with some definitions needed to accomodate the larger structure.

For an integer $n \geq 2$ let $\langle\!\langle n \rangle\!\rangle$ denote the multiplicative subgroup of $\mathbf{Q}$ generated by the (prime) divisors of $n$. We have that $\langle\!\langle n \rangle\!\rangle$ is the group of multiplicative units of $\mathbf{Z}[\frac{1}{n}]$, and that $n$ is prime iff $\langle n \rangle = \langle\!\langle n \rangle\!\rangle$. We let $U_{n-1}$ denote the group of units in $\mathbf{Z}_{n-1}$ and we note that all divisors of $n$ are in $U_{n-1}$. We let $D_{n-1}$ denote the subgroup of $U_{n-1}$ generated by the divisors of $n$. The group $U_{n-1}$ acts on $\mathbf{Z}_{n-1}$ by multiplication and this action gives an isomorphism from $U_{n-1}$ to $\mathrm{Aut}(\mathbf{Z}_{n-1})$. It seems to be a non-trivial number theoretical problem to determine the structure of $D_{n-1}$.

Note that $\phi_n : \langle\langle n \rangle\rangle \to U_{n-1}$ is a homomorphism whose kernel contains $\langle n \rangle$ and whose image is $D_{n-1}$. Thus we get an induced epimorphism $\phi_n : \langle\langle n \rangle\rangle/\langle n \rangle \to D_{n-1}$.

If $A$ and $B$ are groups and there is a homomorphims $\eta : B \to \mathrm{Aut}(A)$, then we can form the semidirect product $A \rtimes B$ using the action of $B$ on $A$ given by $\eta$. We can think of $A \rtimes B$ as acting on $A$ affinely in that $(a, b) \in A \rtimes B$ takes $g \in A$ to $g^b a$. If $A$ is abelian, the result can be written $gb + a$. By sending the function $x \mapsto mx + a$ to the pair $(a, m)$, we can identify $\mathrm{Aff}(\mathbf{R})$, the affine functions on $\mathbf{R}$, with $\mathbf{R} \rtimes \mathbf{R}^*$ (where $\mathbf{R}^*$ is the multiplicative group of non-zero elements of $\mathbf{R}$) and we refer to $\mathbf{Z}_{n-1} \rtimes \mathrm{Aut}(\mathbf{Z}_{n-1}) = \mathbf{Z}_{n-1} \rtimes U_{n-1}$ as $\mathrm{Aff}(\mathbf{Z}_{n-1})$, the *group of affine transformations* of $\mathbf{Z}_{n-1}$.

We will use the following standard constructions that give new semi-direct products from $A \rtimes B$.

(i) If $H \subseteq A$ is invariant under the action of $G \subseteq B$ we get $H \rtimes G \subseteq A \rtimes B$ with the action of $H \rtimes G$ on $H$ given as a restriction of the original action on $A$.

(ii) If $N \triangleleft A$ is invariant under the action of $B$, then we get a surjection $A \rtimes B \to (A/N) \rtimes B$ with the action of $(A/N) \rtimes B$ on $A/N$ induced by the original action on $A$.

(iii) If $K \triangleleft B$ centralizes $A$ (is in the kernel of $\eta$), then we get a surjection $A \rtimes B \to A \rtimes (B/K)$ with the action of $A \rtimes (B/K)$ on $A$ induced by the original action on $A$.

Recall that $\Pi_n$ is the group of permutations of the set $\{0, \ldots, n-1\}$. We let $\Pi_n^0$ be the subgroup of $\Pi_n$ consisting of all elements that fix 0.

**Lemma 3.2.1.** *Let $n \geq 2$ be an integer and let $f : \mathbf{R} \to \mathbf{R}$ be a homeomorphism satisfying Properties 1,2, 4 and 5 of Definition 1.1.1. Then all the slopes of $f$ are units in $\mathbf{Z}[\frac{1}{n}]$, and therefore elements of $\langle\langle n \rangle\rangle$. If $G$ is a group of homeomorphisms $f : \mathbf{R} \to \mathbf{R}$ satisfying Properties 1,2, 4 and 5 of Definition 1.1.1 and*

*(3′) the slopes of $f$ are in a single coset of $\langle\langle n \rangle\rangle/\langle n \rangle$,*

*then:*

(a) *$\phi_n$ is constant on the set of slopes of any given $f \in G$,*

(b) *There are homomorphisms*

$$G \to \mathbf{Z}_{n-1} \rtimes \langle\langle n \rangle\rangle/\langle n \rangle \to \mathbf{Z}_{n-1} \rtimes D_{n-1} \subseteq \mathrm{Aff}(\mathbf{Z}_{n-1})$$

*where the first homomorphism is given by $f \mapsto (0f\phi_n, xf')$ for any real $x$.*

(c) *The composition is a homomorphism $\pi_n : G \to \mathrm{Aff}(Z_{n-1}) \subseteq \Pi_{n-1}$ giving the action of any $f \in G$ on the values of $x\phi_n$ for $x \in \mathbf{Z}[\frac{1}{n}]$.*

(d) *If $f \in G$ fixes 0, then $f\pi_n \in \Pi_{n-1}^0$ is the automorphism of $Z_{n-1}$ given by multiplication by $xf'\phi_n$ for any $x \in \mathbf{Z}[\frac{1}{n}]$.*

*Proof.* If $s$ is the slope of one affine segment of $f$, and two points are given in the domain of this segment which are in $\mathbf{Z}[\frac{1}{n}]$ at a distance $1/n^k$ of each other, then their images will be in $\mathbf{Z}[\frac{1}{n}]$ at a distance $s/n^k$. This shows that $s \in \mathbf{Z}[\frac{1}{n}]$. By Lemma 1.1.2, $f^{-1}$ also satisfies 1,2,4,5, so $1/s \in \mathbf{Z}[\frac{1}{n}]$ and $s$ is a unit in $\mathbf{Z}[\frac{1}{n}]$.

Now suppose that all $f \in G$ also satisfy 3′.
(a) Since $\phi_n$ is a ring homomorphism, and $n\phi_n = 1$, $\phi_n$ is constant on each coset in $\langle\langle n \rangle\rangle/\langle n \rangle$.

(b) An $f \in G$ is made up of pieces of functions in $\mathrm{Aff}(\mathbf{R})$ so we can define $[f]$ as the set of affine functions used by pieces of $f$. For $f, g \in G$ we have $[fg] \subseteq [f][g]$. Identifying $\mathrm{Aff}(\mathbf{R})$ with $\mathbf{R} \rtimes \mathbf{R}^*$ let us note that $[f] \subseteq \mathbf{Z}[\frac{1}{n}] \rtimes \langle\!\langle n \rangle\!\rangle$ which is a subgroup of $\mathbf{R} \rtimes \mathbf{R}^*$ by *(i)* above. Since $\Delta_n \subseteq \mathbf{Z}[\frac{1}{n}]$ is invariant under multiplication by $\langle\!\langle n \rangle\!\rangle$ and $\mathbf{Z}_{n-1} \simeq \mathbf{Z}[\frac{1}{n}]/\Delta_n$ is centralized by multiplication by $\langle n \rangle$, we get surjections

$$\mathbf{Z}[\tfrac{1}{n}] \rtimes \langle\!\langle n \rangle\!\rangle \to \mathbf{Z}_{n-1} \rtimes \langle\!\langle n \rangle\!\rangle \to \mathbf{Z}_{n-1} \rtimes \langle\!\langle n \rangle\!\rangle / \langle n \rangle$$

from *(ii)* and *(iii)* above. We are interested in $[f]^*$ the image of $[f]$ under the composition of these surjections.

Any two affine pieces of $f$ that are adjacent in the graph of $f$ give elements in $[f]$ that agree on some element of $\mathbf{Z}[\frac{1}{n}]$. For $h_1$ and $h_2$ in $[f]$ and an $x \in \mathbf{Z}[\frac{1}{n}]$ on which they agree, we have $xh'_1 + 0h_1 = xh'_2 + 0h_2$ which implies that $0h_2 - 0h_1 = x(h'_1 - h'_2)$. Applying $\phi_n$ and using (a) gives $0h_1\phi_n = 0h_2\phi_n$. This fact and (a) imply that $[f]^*$ is a single element of $\mathbf{Z}_{n-1} \rtimes \langle\!\langle n \rangle\!\rangle / \langle n \rangle$. Now $[fg]^* \subseteq [f]^*[g]^*$ gives the first homomorphism.

The kernel of $\phi_n : \langle\!\langle n \rangle\!\rangle / \langle n \rangle \to D_{n-1}$ centralizes $\mathbf{Z}_{n-1}$ and we get the second homomorphism from *(iii)* above.

(c) Consider again the set $[f]$ of affine actions on $\mathbf{R}$ used by pieces of an $f \in G$. Under the surjection $\mathbf{Z}[\frac{1}{n}] \rtimes \langle\!\langle n \rangle\!\rangle \to \mathbf{Z}_{n-1} \rtimes \langle\!\langle n \rangle\!\rangle$ induced by $\phi_n : \mathbf{Z}[\frac{1}{n}] \to \mathbf{Z}_{n-1}$, these actions are reduced to their actions on the values of $\phi_n$.

(d) This is immediate. $\qquad\blacksquare$

**Lemma 3.2.2.** *Let $n \geq 2$ be an integer and let $h : \mathbf{R} \to \mathbf{R}$ be an orientation preserving PL homeomorphism. The following are equivalent:*

  *(a) $h$ normalizes $A_n(\mathbf{R})$.*

  *(b) $h$ normalizes $B_n(\mathbf{R})$.*

  *(c) All of the following hold:*

   *(i) all the breaks of $h$ occur at elements of $\mathbf{Z}[\frac{1}{n}]$*

   *(ii) $(\mathbf{Z}[\frac{1}{n}])h = \mathbf{Z}[\frac{1}{n}]$*

   *(iii) all the slopes of $h$ are in a single coset of $\langle\!\langle n \rangle\!\rangle / \langle n \rangle$.*

*Proof.* (a) $\Rightarrow$ (b) Follows from Proposition 1.4.2.

(b) $\Rightarrow$ (c) Notice that since $N(B_n(\mathbf{R}))$ is a group, we actually have that both $h$ and $h^{-1}$ normalize $B_n(\mathbf{R})$, and it suffices to check each of the three conditions for one of $h$ and $h^{-1}$. For the first condition, let $x$ be a break point of $h^{-1}$, and let $y = xh^{-1}$. By Proposition 1.2.2 there is $f \in B_n(\mathbf{R})$ that maps a neighborhood of $y$ affinely into the domain of an affine segment of $h$. Since $f$ does not break at $y$ and $h$ does not break at $yf$, we know that $h^{-1}fh$ breaks at $x$. Since $h^{-1}fh \in B_n(\mathbf{R})$ we get $x \in \mathbf{Z}[\frac{1}{n}]$.

We get the second condition from Lemma 1.5.2, Proposition 1.3.3 and the fact that elements of $B_n(\mathbf{R})$ take $\mathbf{Z}[\frac{1}{n}]$ onto $\mathbf{Z}[\frac{1}{n}]$.

For the third condition, let $x, y \in \mathbf{R}$ be points where the derivative $h'$ is continuous. We want to show to show that $xh' \equiv yh' \pmod{\langle n \rangle}$, i.e., $(xh')^{-1}yh' \in \langle n \rangle$. Since $h$ is affine in neighborhoods of $x$ and $y$, we may assume that $x, y \in \mathbf{Z}[\frac{1}{n}]$ and $x\phi_n = y\phi_n$. By Proposition 1.2.2 there is $f \in B_n(\mathbf{R})$ such that $xf = y$, and moreover we may assume that $f'$ is continuous at $x$. Since $h$ normalizes $B_n(\mathbf{R})$, we have $h^{-1}fh \in B_n(\mathbf{R})$. If we let $z = xh$, then the chain rule gives $z(h^{-1}fh)' = (zh^{-1'})(zh^{-1}f')(zh^{-1}fh') = (xh')^{-1}(xf')(yh')$. Since the slopes of $f$ and $h^{-1}fh$ are in $\langle n \rangle$, we get $(xh')^{-1}(yh') \in \langle n \rangle$, as needed.

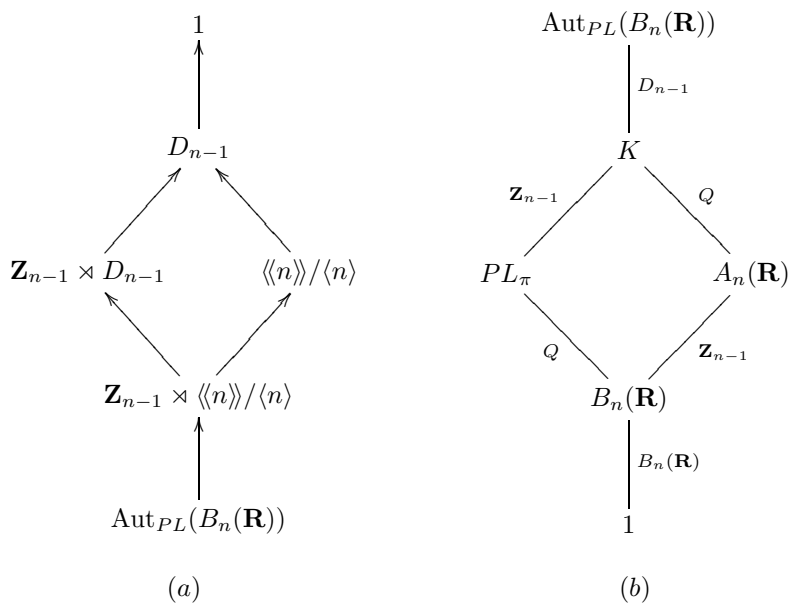

Figure 1. PL automorphisms

(c) $\Rightarrow$ (a) Let $f \in A_n(\mathbf{R})$. Since $h$ satisfies Conditions 1,2,4 and 5 of Definition 1.1.1, then so does $h^{-1}fh$. Since all the slopes of $h$ are in the same coset of $\langle\langle n \rangle\rangle/\langle n \rangle$, then the chain rule implies that the slopes of $h^{-1}fh$ are in $\langle n \rangle$ and $h^{-1}fh \in A_n(\mathbf{R})$. $\quad\square$

**Lemma 3.2.3.** *For each coset of $\langle\langle n \rangle\rangle/\langle n \rangle$ there is a homeomorphism $h : \mathbf{R} \to \mathbf{R}$ in $\mathrm{Aut}_{PL}(B_n(\mathbf{R}))$ with slopes in the given coset. Moreover, $h$ commutes with $t_{n-1}$ and fixes $0$.*

*Proof.* It suffices to prove the statement for cosets of the form $p\langle n \rangle$ where $p$ is a prime divisor of the composite $n$. Consider the map $h_1 : [0, n-1] \to [0, n-1]$ consisting of two affine segments—one from $(0,0)$ to $(\frac{n}{p}-1, n-p)$ with slope $p$ and the second from $(\frac{n}{p}-1, n-p)$ to $(n-1, n-1)$ with slope $p/n$. Since $p \mid n$, the break points are all in $\mathbf{Z} \subseteq \mathbf{Z}[\frac{1}{n}]$. Now $h = h_1\Phi_n$ has the required properties. $\quad\square$

**Proposition 3.2.4.** *There are homomorphisms that commute as shown in Figure 1(a), with the kernels of those compositions with domain $\mathrm{Aut}_{PL}(B_n(\mathbf{R}))$ forming corresponding subgroups and consecutive quotients as shown in Figure 1(b). The group $Q$ is defined only as the kernel of $\phi_n : \langle\langle n \rangle\rangle/\langle n \rangle \to D_{n-1}$. The group $K$ is the set of elements $f \in \mathrm{Aut}_{PL}(B_n(\mathbf{R}))$ with $xf'\phi_n = 1$ for all $x \in \mathbf{Z}[\frac{1}{n}]$, and $PL_\pi$ is the set of elements with $f\pi_n$ trivial. The diagram (b) is a lattice in that $PL_\pi \cap A_n(\mathbf{R}) = B_n(\mathbf{R})$ and $PL_\pi A_n(\mathbf{R}) = K$.*

*Proof.* Two of the homomorphisms in (a) come from Lemma 3.2.1. The remaining homomorphisms and the fact that the diagram commutes are straightforward. Checking that $K$, $PL_\pi$, $A_N(\mathbf{R})$ and $B_n(\mathbf{R})$ are the kernels of the corresponding homomorphisms is also straightforward from the definitions. That (b) forms a lattice follows from the fact that the corresponding diagram of kernels of homomorphisms with domain $\mathbf{Z}_{n-1} \rtimes \langle\langle n \rangle\rangle/\langle n \rangle$ forms a lattice. $\quad\blacksquare$

**Corollary 3.2.4.1.** *We have isomorphisms*

$$\mathrm{Out}_{PL}(A_n(\mathbf{R})) \simeq \langle\!\langle n \rangle\!\rangle / \langle n \rangle, \quad \text{and} \quad \mathrm{Out}_{PL}(B_n(\mathbf{R})) \simeq \mathbf{Z}_{n-1} \rtimes \langle\!\langle n \rangle\!\rangle / \langle n \rangle.$$

$\square$

From here we can determine the torsion of $\mathrm{Out}_{PL}(A_n(\mathbf{R}))$ and $\mathrm{Out}_{PL}(B_n(\mathbf{R}))$.

**Corollary 3.2.4.2.** *Let $n \geq 2$, let $k$ be the largest positive integer such that $n$ is a $k$-th power, and let $m$ be such that $n = m^k$. Identify $\mathrm{Out}_{PL}(A_n(\mathbf{R}))$ and $\mathrm{Out}_{PL}(B_n(\mathbf{R}))$ with the groups $\langle\!\langle n \rangle\!\rangle / \langle n \rangle$ and $\mathbf{Z}_{n-1} \rtimes \langle\!\langle n \rangle\!\rangle / \langle n \rangle$ as given in Corollary 3.2.4.1.*

(a) *The torsion elements of $\mathrm{Out}_{PL}(A_n(\mathbf{R}))$ form the cyclic subgroup $\langle m \rangle / \langle n \rangle$ of order $k$ in $\langle\!\langle n \rangle\!\rangle / \langle n \rangle$.*

(b) *The torsion elements of $\mathrm{Out}_{PL}(B_n(\mathbf{R}))$ form the subgroup $\mathbf{Z}_{n-1} \rtimes \langle m \rangle / \langle n \rangle$ in $\mathbf{Z}_{n-1} \rtimes \langle\!\langle n \rangle\!\rangle / \langle n \rangle$. The orders of the torsion elements are all the divisors of $\frac{k}{d}(m^d - 1)$ where $d$ ranges over all divisors of $k$.*

(c) *$\mathrm{Out}_{PL}(A_n(\mathbf{R}))$ and $\mathrm{Out}_{PL}(B_n(\mathbf{R}))$ have elements of infinite order if and only if $n$ is not a prime power.*

*Proof.* (a) Write $m = p_1^{\alpha_1} \cdots p_l^{\alpha_l}$ with $p_1, \ldots, p_l$ distinct primes and $\alpha_1, \ldots, \alpha_l > 0$. Let $s \in \langle\!\langle n \rangle\!\rangle$, $s = p_1^{\beta_1} \cdots p_l^{\beta_l}$ with $\beta_1, \ldots, \beta_l \in \mathbf{Z}$ be such that $s \langle n \rangle$ has finite order $u$. Then $s^u = n^v = m^{kv}$ for some $v \in \mathbf{Z}$, and $\gcd(u, v) = 1$. It then follows that $u\beta_i = kv\alpha_i$, and $u | k\alpha_i$. By the choice of $k$ we have $\gcd(\alpha_1, \ldots, \alpha_l) = 1$ so $u | k$, and $s = m^{\frac{k}{u}v}$. Therefore the torsion elements of $\langle\!\langle n \rangle\!\rangle / \langle n \rangle$ are precisely the powers of $m\langle n \rangle$, i.e., $\langle m \rangle / \langle n \rangle$.

(b) For $s \in \langle\!\langle n \rangle\!\rangle$ let us also use $s$ to denote its coset $s\langle n \rangle$ in $\langle\!\langle n \rangle\!\rangle / \langle n \rangle$. Since the action of $\langle\!\langle n \rangle\!\rangle / \langle n \rangle$ on $\mathbf{Z}_{n-1}$ is by multiplication, for any $a \in \mathbf{Z}_{n-1}$, $(a, s)^u = (a + sa + \cdots + s^{u-1}a, s^u)$. So $(a, s)$ has finite order if and only if $s\langle n \rangle$ has finite order in $\langle\!\langle n \rangle\!\rangle / \langle n \rangle$. Thus the torsion elements of $\mathbf{Z}_{n-1} \rtimes \langle\!\langle n \rangle\!\rangle / \langle n \rangle$ form the subgroup $\mathbf{Z}_{n-1} \rtimes \langle m \rangle / \langle n \rangle$.

For $d$ a divisor of $k$,

$$\left(1, m^d\right)^{\frac{k}{d}} = \left(\frac{m^k - 1}{m^d - 1}, m^k\right) = \left(\frac{n-1}{m^d - 1}, 1\right)$$

so the element $(1, m^d)$ has order $\frac{k}{d}(m^d - 1)$.

For any $0 \leq v < k$ if we let $d = \gcd(v, k)$, $v' = v/d$, $k' = k/d$, then for $a \in \mathbf{Z}_{n-1}$ we have

$$
\begin{aligned}
(a, m^v)^{\frac{k}{d}(m^d - 1)} &= \left(a \frac{(m^k)^{v'} - 1}{(m^d)^{v'} - 1}, m^{v'k}\right)^{m^d - 1} \\
&= \left(a \frac{1 + m^k + \cdots + m^{k(v'-1)}}{1 + m^d \cdots + m^{d(v'-1)}} \cdot \frac{m^k - 1}{m^d - 1}, 1\right)^{m^d - 1} \\
&= \left(a \frac{1 + (m^d)^{k'} + \cdots + (m^d)^{k'(v'-1)}}{1 + m^d \cdots + m^{d(v'-1)}}(n-1), 1\right) = (0, 1).
\end{aligned}
$$

In the last line, the denominator divides the numerator since a standard exercise using the factorization of $x^r - 1$ into cyclotomic polynomials [9, (2) on P. 191] shows that if $P(x) = 1 + x + \cdots x^{\nu-1}$ and $\gcd(\nu, k) = 1$, then $P(x)$ divides $P(x^k)$.

(c) Follows from (a) and (b) and the definition of $\langle\!\langle n \rangle\!\rangle$. $\square$

## 4. Relations among the automorphism groups

In this section we let the Thompson groups with elements whose slopes are in $\langle m \rangle$ interact with groups with elements whose slopes are in $\langle n \rangle$. By being careful, we keep all references to the generators $g_i$ confined to those with non-negative subscript—the generators of $F_{m,0}$ and $F_{n,0}$.

Because of Section 2.6 we are interested in the automorphism groups of certain pairs $(G, g)$. We will use homomorphisms between pairs to relate their automorphism groups. For groups $G$ and $H$ with $g \in G$ and $h \in H$, we say that $\theta : (G, g) \to (H, h)$ is a *homomorphism of the pairs* if $\theta : G \to H$ is a homomorphism and $g\theta = h$. Recall that the definition of the automorphism group of the pair $(G, g)$ does not require $g \in G$.

This section studies monomorphisms between the pairs $(F_{n,\infty}, t_{n-1})$ for various values of $n$ and uses this information to develop two techniques, lifting and rotating, that create "new" automorphisms of $B_n(\mathbf{R})$ and $A_n(\mathbf{R})$ from "old" ones. We also study how lifting and rotating interact with the homomorphisms $\pi_n$ to show that the $\pi_n$ are onto.

### 4.1. Lifts.

We are interested in automorphisms of the pair $(F_{n,\infty}, t_{n-1})$. We will obtain interesting automorphisms of the pair $(F_{n,\infty}, t_{n-1})$ from automorphisms of the pair $(F_{m,\infty}, t_{m-1})$, $m < n$, by the process that we will call *lifting* whose details are given in Theorem 4.1.6. We do this by showing that there is a monomorphism from $(F_{m,\infty}, t_{m-1})$ to $(F_{n,\infty}, t_{n-1})$ that is "extensible" in that every automorphism of the pair $(F_{m,\infty}, t_{m-1})$ extends to one of the pair $(F_{n,\infty}, t_{n-1})$.

**Definition 4.1.1.** The proof of the next lemma, and several others after, will use the following function. Let $2 \le m < n$ be integers. For every $j \in \mathbf{Z}$ there are unique $0 \le i \le (m-2)$ and $k \in \mathbf{Z}$ for which $j = i + k(m-1)$. If we set $j\zeta_{m,n} = i + k(n-1)$ with $i$, $j$ and $k$ as in the previous sentence, then we have a well defined, strictly increasing function from $\mathbf{Z}$ to $\mathbf{Z}$ that satisfies satisfies $(j+m-1)\zeta_{m,n} = j\zeta_{m,n}+n-1$, and that "commutes" with $\phi$ in that $(j\zeta_{m,n})\phi_n = j\phi_m$ for any $j \in \mathbf{Z}$.

**Proposition 4.1.2.** *For integers $2 \le m < n$, there is a unique monomorphism of group pairs $\tau_{m,n} : (F_{m,\infty}, t_{m-1}) \to (F_{n,\infty}, t_{n-1})$ so that $(g_{m,i})\tau_{m,n} = g_{n,i}$ for $0 \le i \le (m-2)$. For this monomorphism, $(g_{m,j})\tau_{m,n} = g_{n,j\zeta_{m,n}}$ for all $j$.*

*Proof.* With $j$ equal to the unique expression $i+k(m-1)$ with $0 \le i \le (m-2)$ and $k \in \mathbf{Z}$, we get that $g_{m,j}$ is equal to a unique $t_{m-1}^{-k}(g_{m,i})t_{m-1}^{k}$ for $0 \le i \le (m-2)$ and $k \in \mathbf{Z}$. From the requirement that $t_{m-1}\tau_{m,n} = t_{n-1}$, we get

$$(g_{m,j})\tau_{m,n} = (t_{m-1}^{-k}g_{m,i}t_{m-1}^{k})\tau_{m,n} = t_{n-1}^{-k}g_{n,i}t_{n-1}^{k} = g_{n,i+k(n-1)} = g_{n,j\zeta_{m,n}}$$

and the function $\tau_{m,n}$ is uniquely determined on its generators. There is thus at most one monomorphism with the properties given.

The defining relations of $F_{m,\infty}$ take the form $g_{m,i}^{-1}g_{m,j}g_{m,i} = (g_{m,j+n-1})$ whenever $i < j$. Now

$$\begin{aligned}
(g_{m,i}^{-1}\tau_{m,n})(g_{m,j}\tau_{m,n})(g_{m,i}\tau_{m,n}) &= g_{n,i\zeta_{m,n}}^{-1}g_{n,j\zeta_{m,n}}g_{n,i\zeta_{m,n}} \\
&= g_{n,j\zeta_{m,n}+n-1} \\
&= g_{n,(j+m-1)\zeta_{m,n}} \\
&= (g_{m,j+m-1})\tau_{m,n}
\end{aligned}$$

where the second equality holds because $i\zeta_{m,n} < j\zeta_{m,n}$. Thus $\tau_{m,n}$ preserves the relations of $F_{m,\infty}$ and extends to a homomorphism. It is a monomorphism by Proposition 2.1.5 since the elements $g_{n,i}$ do not commute. □

The next few lemmas will be needed to show that automorphisms of the image of the monomorphism of the previous lemma extend to automorphisms of the larger group. We refer to Definition 2.2.1 for what it means for a subdivision of $\mathbf{R}$ to be supported in a subset $A \subseteq \mathbf{R}$.

**Lemma 4.1.3.** *For $n \geq 2$, let $P = g_{i_1} g_{i_2} \cdots g_{i_s}$ be a positive word in semi-normal form in the generators of $F_{n,0}$ and let $k < l$ be integers. Then the subdivision $D$ associated with $P$ is supported in $[k,l]$ if and only if $k \leq i_1$ and for all $1 \leq j \leq s$ we have $i_j < l + (n-1)(j-1)$.*

Note the shift in role of $j$ from subscript to value in the last inequality.

*Proof.* The interval subdivided by $g_{i_1}$ is the interval $[i_1, i_1 + 1]$ so the conditions on $i_1$ are necessary and sufficient as far as $i_1$ goes. We need to show that the second condition is necessary and sufficient to guarantee that the interval subdivided by each remaining generator is contained somewhere within $[k,l]$. Since $P$ is in semi-normal form, $i_j \geq i_{j-1}$ for each $2 \leq j \leq s$. Let $[k_j, k_j + 1]$, $k_j \in \mathbf{Z}$, be the standard interval containing the interval subdivided by $g_{i_j}$. It is clear from the numbering scheme of Definition 2.2.3 that the interval subdivided by $g_{i_j}$ is either contained in $[k_{j-1}, k_{j-1} + 1]$ and $k_j = k_{j-1}$, or it is some $[k_j, k_j + 1]$ with $k_j \in \mathbf{Z}$ and $k_j > k_{j-1}$. Thus the interval subdivided by $g_{i_j}$ cannot be to the left of $k$.

In the numbering scheme of Definition 2.2.3, the interval $[k_j, k_j + 1]$ is assigned the number $k_j$ before any subdivisions are done and this number increases by $(n-1)$ each time an interval to the left of $[k_j, k_j + 1]$ is subdivided. If $k_j > k_{j-1}$, then all generators in $P$ before $g_{i_j}$ are associated to subdivisions of intervals to the left of $[k_j, k_j+1]$ and the number of $[k_j, k_j+1]$ is $k_j + (n-1)(j-1) = i_j$ when the subdivision associated to $g_{i_j}$ is done. If $k_j = k_{j-1}$, then when the subdivision associated to $g_{i_j}$ is done, the number of $[k_{j-1} + 1, k_{j-1} + 2]$ is $k_{j-1} + 1 + (n-1)(j-1) > i_j$ and we get $i_j \leq k_j + (n-1)(j-1)$. If all $k_j < l$, then all $i_j < l + (n-1)(j-1)$. For the converse, assume all $i_j < l + (n-1)(j-1)$. Clearly $k_1 = i_1 < l$. If $k_j = k_{j-1}$, we are done by induction, and if $k_j > k_{j-1}$, then $k_j = i_j - (n-1)(j-1) < l$. □

**Lemma 4.1.4.** *For $2 \leq m < n$, let $P = g_{i_1} g_{i_2} \cdots g_{i_s}$ be a positive word in semi-normal form in the generators of $F_{m,0}$ and let $k < l$ be integers. Assume the subdivision $D$ associated with $P$ is supported in $[k,l]$. Then the subdivision associated with $P\tau_{m,n}$ is supported in $[k',l']$ where $k' = k\zeta_{m,n}$ and $l' = (l-1)\zeta_{m,n} + 1$.*

*Proof.* From Lemma 4.1.3, we have $k \leq i_1$ so $k' = k\zeta_{m,n} \leq i_1\zeta_{m,n}$ and $k'$ satisfies the only inequality from Lemma 4.1.3 that it needs to satisfy. Now $i_j < l + (m-1)(j-1)$ can be rewritten $i_j \leq (l-1) + (m-1)(j-1)$ which implies $i_j\zeta_{m,n} \leq (l-1)\zeta_{m,n} + (n-1)(j-1)$ and we have $i_j\zeta_{m,n} < (l-1)\zeta_{m,n} + 1 + (n-1)(j-1)$ and we are done. □

In the next lemma, the support of a homeomorphism refers to its set of non-fixed points.

**Lemma 4.1.5.** *For $2 \leq m < n$, let $f$ be in $F_{m,0}$, let $k < l$ be integers, let $k' = k\zeta_{m,n}$ and let $l' = (l-1)\zeta_{m,n} + 1$. (a) If $f$ has support in $[k,l]$, then $f\tau_{m,n}$ has support in $[k',l']$. (b) If $f$ is the identity to the left of $k$ and is translation by $(m-1)$* *to the right of $l$, then $f\tau_{m,n}$ is the identity to the left of $k'$ and translation by $(n-1)$ to the right of $l'$.*

*Proof.* Let $f = PN^{-1}$ be in semi-normal form where $P$ and $N$ are positive words. There are allowable subdivisions $D$ and $E$ associated with $P$ and $N$ respectively. We can assume that $D$ subdivides no interval $[s, s+1]$, $s \in \mathbf{Z}$, that is outside $[k, l]$ since such a subdivision would have to be matched exactly by $E$ on $[s, s+1]$ in case (a) or $[s+m-1, s+m]$ in case (b) and the two matching subdivisions could be removed. Similarly, $E$ can be assumed to subdivide no such interval outside $[k, l]$ in case (a) or outside $[k, l+m-1]$ in case (b). We also note the the number of generators in $P$ equals the number in $N$ in case (a) and exceeds the number in $N$ by one in case (b).

By Lemma 4.1.4, $P' = P\tau_{m,n}$ and $N' = N\tau_{m,n}$ are supported in $[k', l']$ and $[k', l']$ respectively in case (a) or $[k', l']$ and $[k', l'+n-1]$ respectively in case (b). Further the number of generators in $P'$ equals the number of generators of $N'$ in case (a) and exceeds the number in $N'$ by one in case (b). Since $f\tau_{m,n} = P'(N')^{-1}$ is the identity near $-\infty$, the conclusions in (a) and (b) follow from the information gathered. $\qquad\square$

**Theorem 4.1.6.** *Let $2 \le m < n$ be integers. Let $\tau_{m,n}$ be as in Proposition 4.1.2. Then there is a monomorphism*

$$\Theta_{m,n} : \mathrm{Aut}^0(F_{m,\infty}, t_{m-1}) \to \mathrm{Aut}^0(F_{n,\infty}, t_{n-1})$$

*so that if $\theta' = \theta\Theta_{m,n}$, then (a) $g_{n,i}\theta' = (g_{m,i}\theta)\tau_{m,n}$ for $0 \le i \le (m-2)$ and $g_{n,j}\theta' = g_{n,j}$ for $(m-1) \le j \le (n-2)$, and (b) $(f\tau_{m,n})\theta' = (f\theta)\tau_{m,n}$ for all $f \in F_{m,\infty}$.*

*Proof.* Note that (a) specifies $\theta'$ on the generators $g_{n,i}$ of $F_{n,\infty}$ for $0 \le i \le (n-2)$. If $\theta'$ exists, then it is sufficiently specified by (a) since it is required to fix $t_{n-1}$ and $F_{n,\infty}$ is generated by $\{t_{n-1}, g_{n,0}, \ldots, g_{n,n-2}\}$. The specifications in (a) for $0 \le i \le (m-2)$ imply that the equality in (b) holds when $f = g_{m,i}$, $0 \le i \le (m-2)$, and when $f = t_{m-1}$. These are the generators of $F_{m,\infty}$ and (b) will hold when it is shown that $\theta'$ is a homomorphism.

We have to show two things. For a given $\theta \in \mathrm{Aut}^0(F_{m,\infty}, t_{m-1})$, we have to show that $\theta' = \theta\Theta_{m,n}$ is an automorphism of the pair $(F_{n,\infty}, t_{n-1})$. We also have to show that $\Theta_{m,n}$ is an injective homomorhism. That $\Theta_{m,n}$ is a homomorphism is clear from the definition. That $\Theta_{m,n}$ is a monomorphism follows from the fact that if $\theta$ does not fix some $g_{m,i}$, then $(g_{m,i}\theta)\tau_{m,n} \ne g_{n,i}$ since $\tau_{m,n}$ is a monomorphism.

To show that $\theta'$ is an automorphism, we have to show that it is an invertible homomorphism. To show that it is a homomorphism, we have to show that the images of the generators of $F_{n,\infty}$ satisfy the relations of $F_{n,\infty}$.

Since $\theta \in \mathrm{Aut}^0(F_{m,\infty}, t_{m-1})$, we know that $\theta$ is realized as conjugation by a self homeomorphism of $\mathbf{R}$ that commutes with $t_{m-1}$ and that fixes 0. We will abuse notation and use $\theta$ to denote this self homeomorphism of $\mathbf{R}$. Thus $0\theta = 0$ and $(m-1)\theta = (m-1)$.

Take an $i$ with $0 \le i \le (m-2)$. The generator $g_{m,i}$ is the identity to the left of $i$ and is equal to $t_{m-1}$ to the right of $i+1$. Thus $\theta^{-1}g_{m,i}\theta$ is the identity to the left of $i\theta$ and is equal to $t_{m-1}$ to the right of $(i+1)\theta$. However, $0 \le i \le (m-2)$ gives $0 \le i < i+1 \le (m-1)$ which leads to $0 \le i\theta < (i+1)\theta \le (m-1)$ and $\theta^{-1}g_{m,i}\theta$ is the identity to the left of 0 and equals $t_{m-1}$ to the right of $m-1$. It follows from Lemma 4.1.5 that $(g_{n,i})\theta'$ is the identity to the left of 0 and equals $t_{n-1}$ to the right of $m-1$.

We have $g_{n,i+n-1}\theta' = (t_{n-1}^{-1}g_{n,i}t_{n-1})\theta' = t_{n-1}^{-1}(g_{n,i}\theta')t_{n-1}$ since this is how $\theta'$ is defined.

Checking that $(g_{n,i}^{-1}\theta')(g_{n,j}\theta')(g_{n,i}\theta') = (g_{n,j+n-1}\theta') = t_{n-1}^{-1}(g_{n,j}\theta')t_{n-1}$ for all $i < j$ is straightforward by considering separately the various cases where $g_{n,i}$ and $g_{n,j}$ are or are not in the image of $\tau_{m,n}$.

This establishes that $\theta'$ is an endomorphism of $F_{n,\infty}$. Since $\Theta_{m,n}$ preserves compositions and inverses, $\theta'$ is invertible. $\qquad\square$

**Corollary 4.1.6.1.** *Given $2 \le m < n$, $f \in F_{m,\infty}$ and $\alpha$ and $\beta$ in $\mathrm{Aut}^0(f_{m,\infty}, t_{m-1})$ with $f\alpha = f\beta$, then $(f\tau_{m,n})(\alpha\Theta_{m,n}) = (f\tau_{m,n})(\beta\Theta_{m,n})$.* $\qquad\square$

### 4.2. Rotations.

**Lemma 4.2.1.** *For $n \ge 2$, let $\alpha \in Aut^0(F_{n,\infty}, t_{n-1})$ and let $h : \mathbf{R} \to \mathbf{R}$ be the homeomorphism fixing 0 and commuting with $t_{n-1}$ that realizes $\alpha$. Let $j \in \mathbf{Z}[\frac{1}{n}]$ and let $r = jh$. Then $k = t_j h t_r^{-1}$ is a homeomorphism of $\mathbf{R}$ that fixes 0 and commutes with $t_{n-1}$. Therefore it realizes a $\beta \in Aut^0(F_{n,\infty}, t_{n-1})$.*

*Proof.* Since $t_{n-1}$ commutes with $h$, $t_j$ and $t_r$, it commutes with $k$. Now, $0k = 0t_j h t_r^{-1} = jh t_r^{-1} = r t_r^{-1} = 0$. $\qquad\square$

In the setting of the lemma above we call $k$ (resp. $\beta$) the *$j$-step rotation* of $h$ (resp. $\alpha$). While the $j$-step rotation is a function from $Aut^0(F_{n,\infty}, t_{n-1})$ to itself, it is not in general a homomorphism. An elementary calculation shows that the $j'$-step rotation of the $j$-step rotation of $h$ is the $(j+j')$-step rotation of $h$. We are mainly interested in the case when $j$ is an integer. Note that the $(n-1)$-step rotation of $h$ is $h$.

**Proposition 4.2.2.** *For $n \ge 2$, let $\alpha \in Aut^0(F_{n,\infty}, t_{n-1})$ and let $h : \mathbf{R} \to \mathbf{R}$ be the homeomorphism fixing 0 and commuting with $t_{n-1}$ that realizes $\alpha$. Let $j$ be an integer, and $k$ be the $j$-step rotation of $h$.*

   (a) *$h$ normalizes $A_n(\mathbf{R})$ if and only if $k$ does.*

   (b) *If $h$ commutes with $t_1$ then $k = h$. Moreover, $h$ commutes with $t_1$ if and only if $h$ equals its one-step rotation.*

   (c) *$k\pi_n = \rho_{j\phi_n}(h\pi_n)\rho_{(j\phi_n)(h\pi_n)}^{-1}$ where $\rho_i$ denotes addition of $i$ (mod $n-1$) in the set $\{0, 1, \ldots, n-2\}$. In particular, if $h$ fixes the residue class of $j$ then $k\pi_n$ and $h\pi_n$ are conjugate.*

   (d) *$h$ fixes residue classes (i.e., $h\pi_n \equiv 1$) if and only if $k$ does.*

*Proof.* Since $h$ is the $(n-1-j)$-step rotation of $k$, in parts (a) and (d) it suffices to prove one direction.

(a) This follows because translations by elements of $\mathbf{Z}[\frac{1}{n}]$ normalize $A_n(\mathbf{R})$.

(b) If $h$ commutes with $t_1$ then $r = jh = 0t_1^j h = 0h t_1^j = 0t_1^j = j$ so $k = t_j h t_r^{-1} = t_1^j h t_1^{-j} = h$. For the converse, note that if $t_1 h t_{1h}^{-1} = h$, then $t_1^{n-1} h t_{1h}^{-n+1} = h$. However, $h$ commutes with $t_1^{n-1}$ and we get that $(t_1 t_{1h}^{-1})^{n-1}$ is the identity. Translations form a torsion free group, so $t_1 = t_{1h}$. Now we know that $h$ commutes with $t_1$.

(c) From Proposition 1.5.4 we have that $(jh)\phi_n = (j\phi_n)(h\pi_n)$, so

$$k\pi_n = (t_j\pi_n)(h\pi_n)(t_{jh}\pi_n)^{-1} = \rho_{j\phi_n}(h\pi_n)\rho_{(jh)\phi_n}^{-1} = \rho_{j\phi_n}(h\pi_n)\rho_{(j\phi_n)(h\pi_n)}^{-1}.$$

(d) Follows at once from (c). $\qquad\square$

4.3. **Symmetric lifts.** We introduce special lifts that will create automorphisms of $A_n$ from automorphisms of $A_2$ whenever $n \geq 2$. See Theorem 6.1.5.

**Theorem 4.3.1.** *Let $2 \leq m < n$ be integers for which $(m-1) \mid (n-1)$ and let $\tau_{m,n} : (F_{m,\infty}, t_1) \to (F_{n,\infty}, t_{n-1})$ be the monomorphism of Proposition 4.1.2. Then there is a monomorphism*

$$\Lambda_{m,n} : \mathrm{Aut}^0(F_{m,\infty}, t_{m-1}) \to \mathrm{Aut}^0(F_{n,\infty}, t_{m-1})$$

*(note the mix of $m$ and $n$ in the right hand group) taking $\theta \in \mathrm{Aut}^0(F_{m,\infty}, t_{m-1})$ to $\theta' \in \mathrm{Aut}^0(F_{n,\infty}, t_{m-1})$ so that $g_{n,i}\theta' = (g_{m,i}\theta)\tau_{m,n}$ for $0 \leq i \leq m-2$.*

*Proof.* The proof is similar to that of Theorem 4.1.6. With $\theta'$ defined on $g_{n,i}$, $0 \leq i \leq m-2$, we can set $g_{n,k}\theta' = (g_{n,j}\theta')\sigma_{m-1}^l$ where $0 \leq j \leq m-2$ and $k = j + l(m-1)$. Thus $\theta'$ is sufficiently specified. The main tasks are to show that $\theta'$ is an automorphism and that $\Lambda_{m,n}$ is an injective homomorphism. An argument identical to that in the proof of Theorem 4.1.6 shows that the $g_{n,k}\theta'$ satisfy the relations of $F_{n,\infty}$. That $\theta'$ is an automorphism follows as it does in Theorem 4.1.6. By the way $\theta'$ is defined, it commutes with $t_{m-1}$, and it is thus an automorphism of $(F_{n,\infty}, t_{m-1})$. That $\Lambda_{m,n}$ is a monomorphism is argued as it is for $\Theta_{m,n}$ in Theorem 4.1.6. $\qquad\blacksquare$

It might be surmised that if $2 < m < n$ and $(m-1) \mid (n-1)$, then a $\theta \in \mathrm{Aut}^0(F_{m,\infty}, t_{m-1})$ that is also a normalizer of $A_m(\mathbf{R})$ should lift to a $\theta' = \theta\Lambda_{m,n} \in \mathrm{Aut}^0(F_{n,\infty}, t_{m-1})$ that is also a normalizer of $A_n(\mathbf{R})$. We will see in Theorem 5.4.8 that this is not always the case.

4.4. **Lifting permutations.** Let $n \geq 2$. For $w \in F_{n,-\infty}$, $w \neq 1$, the leftmost break of $w$ is the largest $a \in \mathbf{R}$ such that $w$ is the identity on the interval $(-\infty, a]$. Also recall that $\Pi_n$ is the permutation group on $\{0, \ldots, n-1\}$ and that $\Pi_n^0$ is the subgroup of elements that fix 0.

**Lemma 4.4.1.** *For $n \geq 2$, let $w \in F_{n,-\infty}$, $w \neq 1$, and let $a$ be the leftmost break of $w$.*
*(a) If $h$ is an orientation preserving normalizer of $F_{n,\infty}$, then $h^{-1}wh \in F_{n,-\infty}$ and the leftmost break of $h^{-1}wh$ is $ah$.*
*(b) Let $w = PN^{-1}$ be a semi-normal form of $w$. Let $v = g_{i_0}g_{i_1}\ldots g_{i_{k-1}}$ be the largest common prefix of $P$ and $N$ as positive words in the generators. Let $g_{i_k}$ be the lowest indexed generator in the remaining suffixes. Then $a\phi_n = i_k\phi_n$.*

*Proof.* Item (a) follows because conjugation takes fixed points to fixed points. To work on (b), we note that $P' = v^{-1}P$ and $N' = v^{-1}N$ start with different generators, and that the lower subscript of the two starting generators is $i_k$. This implies that the semi-normal form $P'(N')^{-1}$ has its leftmost break at $i_k$. The leftmost break of $PN^{-1} = vP'(N')^{-1}v^{-1}$ is at $(i_k)v^{-1}$. Since $v$ is a word in elements with fixed points, we have that $v$ is in $B_n(\mathbf{R})$ and that $v$ preserves values of $\phi_n$. $\qquad\blacksquare$

**Proposition 4.4.2.** *Let $2 \leq m < n$ be integers. Let $\tau_{m,n}$ and $\Theta_{m,n}$ be as in Theorem 4.1.6. Let $\theta \in Aut^0(F_{m,\infty}, t_{m-1})$ and $\theta' = \theta\Theta_{m,n} \in Aut^0(F_{n,\infty}, t_{m-1})$. We will abuse notation and denote by $\theta$ (resp. $\theta'$) the self homeomorphism of $\mathbf{R}$ that fixes 0, commutes with $t_{m-1}$ (resp. $t_{n-1}$) and realizes the automorphism $\theta$ (resp. $\theta'$). Then the permutation $\theta'\pi_n$ is an extension to $\Pi_{n-1}^0$ of $\theta\pi_m$ in the sense that for $1 \leq i \leq m-2$, $i(\theta'\pi_n) = i(\theta\pi_m)$ and for $m-1 \leq i \leq n-2$, $i(\theta'\pi_n) = i$.*

*Proof.* Since $\theta'$ carries the fixed point set of $g_{n,i}$ to that of $g_{n,i}\theta'$, we have that $i\theta'$ is the leftmost break of $(g_{n,i})\theta'$. By the way $\Theta_{m,n}$ is defined for $m-1 \leq i \leq n-2$, $(g_{n,i})\theta' = g_{n,i}$ so its leftmost break is $i$ and $i(\theta'\pi_n) = i$ in these cases. For $1 \leq i \leq m-2$, write $(g_{m,i})\theta = PN^{-1}$ in semi-normal form, and let $v = g_{i_0}g_{i_1}\ldots g_{i_{k-1}}$ be the largest common prefix of $P$ and $N$ as positive words in the generators. Let $g_{i_k}$ be the lowest indexed generator in the remaining suffixes. Now $(P\tau_{m,n})(N\tau_{m,n})^{-1}$ is a semi-normal form for $(g_{n,i})\theta'$ where $v\tau_{m,n}$ is the largest common prefix of $P\tau_{m,n}$ and $N\tau_{m,n}$ and $g_{m,i_k}\tau_{m,n} = g_{n,i_k\zeta_{m,n}}$ is the lowest indexed generator in the remaining suffixes. Now

$$i(\theta'\pi_n) = (i_k\zeta_{m,n})\phi_n = i_k\phi_m = i(\theta\pi_m)$$

by the previous lemma and the remarks in Definition 4.1.1. $\qquad\blacksquare$

**Theorem 4.4.3.** *For any $n \geq 3$, the homomorphism $\pi_n : \mathrm{Aut}(F_{n,\infty}, t_{n-1}) \to \Pi_{n-1}$ is surjective. Further, the restriction $\pi_n : \mathrm{Aut}^0(F_{n,\infty}, t_{n-1}) \to \Pi_{n-1}^0$ is surjective.*

*Proof.* We start with the restricted homomorphism. There is nothing to show for $n = 3$ since $\Pi_2^0$ is trivial. For $n = 4$, Lemma 3.2.3 gives a self homeomorphism $h$ of $\mathbf{R}$ that fixes 0 and commutes with $t_3$ with all its slopes of the form $2 \cdot 4^k$, $k \in \mathbf{Z}$. Consider $\alpha \in Aut^0(F_{4,\infty}, t_3)$ the automorphism realized by $h$. By Lemma 3.2.1(d), $\alpha\pi_4$ is multiplication by 2 (mod 3), and therefore it is the transposition $\begin{pmatrix} 1 & 2 \end{pmatrix}$. For $n > 4$, let $\alpha' = \alpha\Theta_{4,n}$ be the lift of $\alpha$ from 4 to $n$. From Proposition 4.4.2 we have that $\alpha'\pi_n = \begin{pmatrix} 1 & 2 \end{pmatrix}$. From Proposition 4.2.2(c), the permutation of the 1-step rotation of $\alpha'$ is

$$\begin{aligned}
\rho_1(\alpha'\pi_n)\rho_2^{-1} =& \rho_1(\alpha'\pi_n)\rho_{n-3} \\
=& \begin{pmatrix} 0 & 1 & 2 & \ldots & n-3 & n-2 \\ 1 & 2 & 3 & \ldots & n-2 & 0 \end{pmatrix} \cdot \begin{pmatrix} 0 & 1 & 2 & 3 & \ldots & n-2 \\ 0 & 2 & 1 & 3 & \ldots & n-2 \end{pmatrix} \cdot \\
& \begin{pmatrix} 0 & 1 & 2 & 3 & \ldots & n-2 \\ n-3 & n-2 & 0 & 1 & \ldots & n-4 \end{pmatrix} \\
=& \begin{pmatrix} 0 & 1 & 2 & 3 & \ldots & n-2 \\ 0 & n-2 & 1 & 2 & \ldots & n-3 \end{pmatrix}
\end{aligned}$$

which is a cyclic rotation of the elements $\{1, 2, \ldots, n-2\}$. This rotation and the transposition $\begin{pmatrix} 1 & 2 \end{pmatrix}$ generate the whole group $\Pi_{n-1}^0$.

For the unrestricted homomorphism, $\sigma_1$ realized as conjugation by $t_1$ is in $\mathrm{Aut}(F_{n,\infty}, t_{n-1})$ and we get that the cyclic rotation of $\{0, \ldots, n-1\}$ is in the image. For $n = 3$, this is the only non-trivial element in $\Pi_2$. For $n > 3$, this rotation and the transposition $\begin{pmatrix} 1 & 2 \end{pmatrix}$ found above generate all of $\Pi_{n-1}$. This completes the proof. $\qquad\blacksquare$

## 5. Restrictions on automorphisms

In this section we extract geometric information from the circle, and transfer that information to the line to get "negative" results about automorphisms—that certain automorphisms are not PL and that certain automorphisms of $B_n(\mathbf{R})$ are not automorphisms of $A_n(\mathbf{R})$.

5.1. **Communication between line and circle.** Because of Section 2.6, we are interested in $\mathrm{Aut}^0(F_{n,\infty}, t_{n-1})$. If $\alpha \in \mathrm{Aut}^0(F_{n,\infty}, t_{n-1})$, then $\alpha$ is realized by conjugation by some homeomorphism $h : \mathbf{R} \to \mathbf{R}$ that commutes with $t_{n-1}$ and fixes 0. This tells us that $[0, n-1]h = [0, n-1]$. From the definition

$$xg_i = \begin{cases} x & x < i, \\ n(x - i) + i & i \le x \le i + 1, \\ x + n - 1 & x > i + 1. \end{cases}$$

we know that $g_i$ is the identity on $(-\infty, i]$, is $t_{n-1}$ on $[i + 1, \infty)$ and is $\widetilde{\mu_n} t_{n-1}^{-i}$ on $[i, i+1]$ where $\widetilde{\mu_n} : \mathbf{R} \to \mathbf{R}$ is defined by $x\widetilde{\mu_n} = nx$ and is the lift of $\mu_n : S_{n-1} \to S_{n-1}$ that fixes 0. Now

$$h^{-1}g_i h = \begin{cases} \text{identity} & \text{on } (-\infty, i]h, \\ (h^{-1}\widetilde{\mu_n}h)t_{n-1}^{-i} & \text{on } [i, i+1]h \\ t_{n-1} & \text{on } [i + 1, \infty)h. \end{cases}$$

since $h$ commutes with $t_{n-1}$. Recall $S_{n-1} = \mathbf{R}/(n-1)\mathbf{Z}$, and note that $h$ determines a well defined self homeomorphism $h_S$ of $S_{n-1}$ with

$$(x + (n-1)\mathbf{Z})h_S = (xh) + (n-1)\mathbf{Z}$$

and that $h$ is the lift of $h_S$ that fixes 0. If we refer to the restriction of $g_i$ to $[i, i+1]$ as the "active" part of $g_i$ and the restriction of $h^{-1}g_i h$ to $[i, i+1]h$ as the "active" part of $h^{-1}g_i h$, then we have that the active part of $g_i$ is a restriction to $[i, i+1]$ of that lift of $\mu_n$ that fixes $i$ and the active part of $h^{-1}g_i h$ is a restriction to $[i, i+1]h$ of that lift of $h_S^{-1}\mu_n h_S$ that fixes $ih$. For $0 \le i < n-1$, we have $ih \in [0, n-1]$. The active parts of the $h^{-1}g_i h$ determine the functions $h^{-1}g_i h$. Since $\{t_{n-1}, g_0 \ldots g_{n-2}\}$ is a generating set for $F_{n,\infty}$, the functions $h^{-1}g_i h$, $0 \le i < n-1$, determine the action of $\alpha \in \mathrm{Aut}^0(F_{n,\infty}, t_{n-1})$. Lastly, we note that the integers in $S_{n-1}$ are the fixed points of $\mu_n$ and their images under $h_S$ are the fixed points of $h_S^{-1}\mu_n h_S$. All of these observations lead to the next definition and lemma.

Let $\alpha \in \mathrm{Aut}^0(F_{n,\infty}, t_{n-1})$ be realized as conjugation by $h : \mathbf{R} \to \mathbf{R}$, and let $h_S : S_{n-1} \to S_{n-1}$ be determined by $h$. We call the degree $n$ function $h_S^{-1}\mu_n h_S$ the *characteristic function* $\chi_\alpha$ of $\alpha$. The remarks above have proven the following.

**Lemma 5.1.1.** *For $n \ge 2$, an $\alpha \in \mathrm{Aut}^0(F_{n,\infty}, t_{n-1})$ is determined by its characteristic function $\chi_\alpha$. Further, for $0 \le i < n-1$, the active part of $g_i\alpha$ is a lift of the restriction of $\chi_\alpha$ to the closed interval between the two consecutive fixed points $ih$ and $(i+1)h$ of $\chi_\alpha$.* $\qquad\square$

5.2. **The calculus of break values and a criterion for piecewise linearity.** If $f$ is PL, then we define the *break value* $xf^b$ of $f$ at $x$ to be $\log(xf'_+) - \log(xf'_-)$ where $xf'_+$ is the derivative from the right of $f$ at $x$ and $xf'_-$ is the derivative from the left of $f$ at $x$. We have $xf^b = 0$ if and only if $f$ is affine on some neighborhood of $x$.

**Lemma 5.2.1.** *For PL functions $f$, $g$ and $h$ with $h$ invertible, we have*

$$x(fg)^b = xf^b + (xf)g^b,$$
$$(xh)(h^{-1})^b = -xh^b, \text{ and}$$
$$(xh)(h^{-1}gh)^b = -xh^b + xg^b + (xg)h^b.$$

*Proof.* The first equation is an additive chain rule and follows from the definition. The second equation follows by applying the chain rule to $x(hh^{-1})^b$ and the third follows from the first two. $\qquad\square$

**Lemma 5.2.2.** *For $n \geq 2$, if $h$ is a normalizer of $B_n(S_{n-1})$, then $h^{-1}\mu_n h$ has slope $n$ on some neighborhood of all its fixed points and thus has break value 0 at all its fixed points.*

*Proof.* Conjugation by $h$ induces an isomorphism from the group of germs of $B_n(S_{n-1})$ with a fixed point $x \in \mathbf{Z}[\frac{1}{n}]$ to the group of germs with the fixed point $xh$. The groups involved are isomorphic to $\mathbf{Z} \times \mathbf{Z}$ with generators "slope $n$ to the right" and "slope $n$ to the left." Conjugation by $h$ preserves "to the right" and "to the left" in that behavior of the conjugate to one side is determined by behavior of the germ being conjugated on an appropriate side (depending on whether $h$ preserves orientation or not). The generator to one side is taken to the generator on the corresponding side or its inverse. But conjugation takes a repelling fixed point to a repelling fixed point. Thus "slope $n$ to the left" must be taken to "slope $n$ to the left" as opposed to "slope $\frac{1}{n}$ to the left." A similar statement applies to the right. Since $\mu_n$ has slope $n$ on both sides of all its fixed points, the result follows. $\qquad\square$

The above lemma and Lemma 5.1.1 give the following.

**Corollary 5.2.2.1.** *For $n \geq 2$, if $h$ is a normalizer of $F_{n,\infty}$ that commutes with $t_{n-1}$, then $h^{-1}g_i h$ has slope $n$ on some $[i, i+\epsilon_1)h$ and on some $(i+1-\epsilon_2, i+1]h$.* $\qquad\square$

**Lemma 5.2.3.** *For $n \geq 2$, if $h$ is a non-trivial PL normalizer of $B_n(S_{n-1})$ that fixes 0, then $h^{-1}\mu_n h$ has a non-zero break value in each open interval $(i, i+1)h$.*

*Proof.* Since $h$ is not the identity on $S_{n-1}$, it has a non-zero break value at some point $x$. Now $(i, i+1)$ hits each point of $S_{n-1}$ at least once under $\mu_n$. Let $y$ be the point in $(i, i+1)$ closest to $i$ so that $y\mu_n$ is a point at which $h$ has a non-zero break value. Note that $y$ is $y'\mu_n$ for a $y' \in (i, y)$ where the distance from $i$ to $y$ is $n$ times the distance from $i$ to $y'$. Thus $y'\mu_n = y$ is a point at which $h$ has zero break value. Now $\mu_n$ has zero break value at all points, so $(yh)(h^{-1}\mu_n h)^b = (y\mu_n)h^b$ which is not zero. $\qquad\square$

### 5.3. The unbent generator proposition and locally inner automorphisms.

All results in the rest of this section hinge on the next proposition and its consequences. This proposition will allow us to derive information by studying the effect of an automorphism on a single generator. Let $\alpha$ be in $Aut^0(F_{n,\infty}, t_{n-1})$. We say that $\alpha$ has an *unbent generator at $i$* if $g_i^\alpha$ is affine on the interval $[i, i+1]h$, where $h$ is the self homeomorphism of $\mathbf{R}$ that realizes $\alpha$.

**Proposition 5.3.1** (Unbent Generator). *For $n \geq 2$, let $\alpha \in Aut^0(F_{n,\infty}, t_{n-1})$ be piecewise linear. If for some $0 \leq i < n-1$, $\alpha$ has an unbent generator at $i$, then $\alpha$ is the identity.*

*Proof.* This follows from Lemma 5.1.1 and the contrapositive of Lemma 5.2.3. $\qquad\square$

The following definitions will help us concentrate on the effect that automorphisms have on single generators.

**Definition 5.3.2.** Let $\alpha$ be in $Aut^0(F_{n,\infty}, t_{n-1})$, $n \geq 2$.
(a) We say that $\alpha$ is *inner* if there is $w \in C^0(B_n(\mathbf{R}), t_{n-1})$ such that $g_i\alpha = w^{-1}g_i w$

for all $i$.

(b) We say that $\alpha$ is *inner at $i$ (via $w$)* if there is $w \in C^0(B_n(\mathbf{R}), t_{n-1})$ such that $g_i\alpha = w^{-1}g_iw$.

(d) We say that $\alpha$ is *locally inner* if it is inner at every $i \in \mathbf{Z}$.

Note that if $\alpha$ is inner at $i$ then it is inner at $i + k(n-1)$ for every $k$, and the same $w$ works for all $k$. Specifically

$$(g_{i+k(n-1)})\alpha = (g_i)\sigma_{n-1}^k\alpha = (g_i)\alpha\sigma_{n-1}^k = (w^{-1}g_iw)\sigma_{n-1}^{k} = w^{-1}g_{i+k(n-1)}w.$$

Note also that $\alpha$ is locally inner if it is inner at every $i \in \{0, \ldots, n-2\}$. Thus if $\alpha$ is locally inner, then its action is completely determined by an $(n-1)$-tuple $(w_0, \ldots, w_{n-2})$ of elements of $C^0(B_n(\mathbf{R}), t_{n-1})$, so that $g_i\alpha = w_i^{-1}g_iw_i$ for $i \in \{0, \ldots, n-2\}$. In turn, by Lemma 3.1.2(a) these $w_i$'s are determined by an $(n-1)$-tuple $(W_0, \ldots, W_{n-2})$ of elements of $F_{n,[0,n-1]}$ with $w_i = W_i\Phi_n$. Clearly, not every such $(n-1)$-tuple determines a locally inner automorphism.

We now define $\text{Aut}_{I_i}^0(F_{n,\infty}, t_{n-1})$ to be the set of elements in $\text{Aut}^0(F_{n,\infty}, t_{n-1})$ that are inner at $i$. We will see shortly that this is a group. We let $\text{Aut}_{g_i}^0(F_{n,\infty}, t_{n-1})$ be those elements in $\text{Aut}^0(F_{n,\infty}, t_{n-1})$ that fix $g_i$. We have $\text{Aut}_{g_i}^0(F_{n,\infty}, t_{n-1}) \subseteq \text{Aut}_{I_i}^0(F_{n,\infty}, t_{n-1})$ since an automorphism fixing $g_i$ is inner at $i$ via 1. We next let $\text{Aut}_{LI}^0(F_{n,\infty}, t_{n-1})$ be the set of locally inner automorphisms in $\text{Aut}^0(F_{n,\infty}, t_{n-1})$. This is the intersection of the $\text{Aut}_{I_i}^0(F_{n,\infty}, t_{n-1})$ as $i$ ranges over $\{0, \ldots, n-2\}$.

**Lemma 5.3.3.** *Let* $0 \leq i < n-1$ *where* $n \geq 2$.

(a) *If* $\alpha \in Aut^0(F_{n,\infty}, t_{n-1})$ *is inner at $i$ via both $w$ and $w'$ in $C^0(B_n(\mathbf{R}), t_{n-1})$, then $w = w'$.*

(b) *Let* $\alpha \in Aut^0(F_{n,\infty}, t_{n-1})$ *be inner at $i$. Then either $\alpha$ is inner or $\alpha$ is not PL.*

(c) $\text{Aut}_{I_i}^0(F_{n,\infty}, t_{n-1})$ *is a group and has $C^0(B_n(\mathbf{R}), t_{n-1})$ as a normal subgroup with complement $\text{Aut}_{g_i}^0(F_{n,\infty}, t_{n-1})$. Any element in $\text{Aut}_{I_i}^0(F_{n,\infty}, t_{n-1})$ and outside $C^0(B_n(\mathbf{R}), t_{n-1})$ is not PL.*

(d) $\langle\text{Aut}_{g_i}^0(F_{n,\infty}, t_{n-1})\rangle\pi_n = \langle\text{Aut}_{I_i}^0(F_{n,\infty}, t_{n-1})\rangle\pi_n = \Pi_n^{0,i,i+1}$ *the group of permutations that fix $\{0, i, i+1\}$.*

(e) $\text{Aut}_{LI}^0(F_{n,\infty}, t_{n-1})$ *is a subgroup of $\text{Aut}_\pi^0(F_{n,\infty}, t_{n-1})$. In general (for example, when $n$ is not a prime power) this subgroup is proper. Any element in $\text{Aut}_{LI}^0(F_{n,\infty}, t_{n-1})$ and not in $C^0(B_n(\mathbf{R}), t_{n-1})$ is not PL.*

*Proof.* (a) If $w \neq w'$, then conjugation by the PL homeomorphism $w'w^{-1}$ is a nontrivial element of $Aut^0(F_{n,\infty}, t_{n-1})$ that fixes $g_i$. This violates Proposition 5.3.1.

(b) If $\alpha$ is inner at $i$ via $w$, then $\alpha' \in Aut^0(F_{n,\infty}, t_{n-1})$ obtained by composing $\alpha$ with conjugation by $w^{-1}$ fixes $g_i$. By Proposition 5.3.1, $\alpha'$ is either trivial or not PL.

(c) To show closure under inversion, let $\alpha \in \text{Aut}^0(F_{n,\infty}, t_{n-1})$ be inner at $i$ via the element $w \in C^0(B_n(\mathbf{R}), t_{n-1})$. From $g_i\alpha = w^{-1}g_iw$, we get that $g_i = w(g_i\alpha)w^{-1}$ and $g_i\alpha^{-1} = (w\alpha^{-1})g_i(w^{-1}\alpha^{-1})$ so $\alpha^{-1}$ is inner at $i$ via $(w^{-1}\alpha^{-1})$. To show closure under composition, let $\alpha_1, \alpha_2 \in \text{Aut}^0(F_{n,\infty}, t_{n-1})$ be inner at $i$ via $w_1, w_2 \in C^0(B_n(\mathbf{R}), t_{n-1})$ respectively. Then

$$\begin{aligned}
g_i\alpha_1\alpha_2 &= (w_1^{-1}g_iw_1)\alpha_2\\
&= (w_1^{-1}\alpha_2)(g_i\alpha_2)(w_1\alpha_2)\\
&= (w_1^{-1}\alpha_2)(w_2^{-1}g_iw_2)(w_1\alpha_2)
\end{aligned}$$

This shows that $\alpha_1\alpha_2$ is inner at $i$ via $w_2(w_1\alpha_2)$ which is clearly in $C^0(B_n(\mathbf{R}), t_{n-1})$. The group $C^0(B_n(\mathbf{R}), t_{n-1})$ is normal in $\mathrm{Aut}^0(F_{n,\infty}, t_{n-1})$ and so is normal in $\mathrm{Aut}^0_{I_i}(F_{n,\infty}, t_{n-1})$. If $\alpha$ is inner at $i$ via $w \in C^0(B_n(\mathbf{R}), t_{n-1})$, then $\alpha'$ obtained by composing $\alpha$ with conjugation by $w^{-1}$ is in $\mathrm{Aut}^0_{g_i}(F_{n,\infty}, t_{n-1})$. Any element common to $C^0(B_n(\mathbf{R}), t_{n-1})$ and $\mathrm{Aut}^0_{g_i}(F_{n,\infty}, t_{n-1})$ is trivial by Proposition 5.3.1. The last provision follows from (b).

(d) We will use superscripts to indicate that certain points are fixed by automorphisms or permutations. The inclusions

$$\langle \mathrm{Aut}^0_{g_i}(F_{n,\infty}, t_{n-1})\rangle \pi_n \subseteq \langle \mathrm{Aut}^0_{I_i}(F_{n,\infty}, t_{n-1})\rangle \pi_n \subseteq \Pi_n^{0,i,i+1}$$

are obvious. For the inclusion $\Pi_n^{0,i,i+1} \subseteq \langle \mathrm{Aut}^0_{g_i}(F_{n,\infty}, t_{n-1})\rangle \pi_n$ consider $\sigma \in \Pi_n^{0,i,i+1}$ and let $\tau = \rho_{n-2-i}^{-1}\sigma\rho_{n-2-i} \in \Pi_n^{0,n-2-i,n-2}$. Now $\tau \in \Pi_{n-1}^{0,n-2-i}$ and by Theorem 4.4.3 there is an $\alpha \in \mathrm{Aut}^0(F_{n-1,\infty}, t_{n-2})$ such that $\alpha\pi_{n-1} = \tau$. Since $\tau$ fixes $(n-2-i)$, Proposition 1.2.2 gives a $W \in F_{0,[0,n-2]}$ taking $(n-2-i)\alpha$ to $(n-2-i)$. If $\alpha'$ is $\alpha$ composed with conjugation by $W\Phi_{n-1}$, then $\alpha' \in \mathrm{Aut}^{0,n-2-i}(F_{n-1,\infty}, t_{n-2})$ and $\alpha'\pi_{n-1} = \tau$. Let $\beta = \alpha'\Theta_{n-1,n}$. By Proposition 4.4.2, $\beta\pi_n = \tau$, and by Theorem 4.1.6, $\beta \in \mathrm{Aut}^{0,n-2}_{g_{n-2}}(F_{n,\infty}, t_{n-1})$. An easy argument using the notion of leftmost break shows that $\beta \in \mathrm{Aut}^{0,n-2-i,n-2}_{g_{n-2}}(F_{n,\infty}, t_{n-1})$. By Lemma 4.2.1, the $(n-2-i)$-step rotation of $\beta$ is $\gamma = t_{n-2-i}^{-1}\beta t_{n-2-i}$ so $\gamma$ is in $\mathrm{Aut}^{0,i,i+1}_{g_i}(F_{n,\infty}, t_{n-1})$. By Proposition 4.2.2, $\gamma\pi_n = \sigma$.

(e) When $n$ is not a prime power, the second sentence follows from Example 5.3.4 given after Corollary 5.3.3.1. The rest follows from (c) and (d). $\qquad\blacksquare$

From Lemma 5.3.3(a), we get a we get a well-defined injective function

$$\iota_n : \mathrm{Aut}^0_{LI}(F_{n,\infty}, t_{n-1}) \to F_{n,[0,n-1]}^{n-1}$$

so that if $\alpha \in \mathrm{Aut}^0_{LI}(F_{n,\infty}, t_{n-1})$ and $\alpha\iota_n = (W_0, \ldots, W_{n-2})$ then for each $i$ with $0 \le i \le n-2$, we have $g_i\alpha = w_i^{-1}g_iw_i$, where $w_i = W_i\Phi_n$. A further consequence of Lemma 5.3.3(a) is the following.

**Corollary 5.3.3.1.** *An element $\alpha \in Aut^0_{LI}(F_{n,\infty}, t_{n-1})$ is inner if and only if $\alpha\iota_n$ is in the diagonal of $F_{n,[0,n-1]}^{n-1}$.* $\qquad\blacksquare$

**Example 5.3.4.** If $n$ is not a prime power, then let $p$ be a prime divisor of $n$. For some $k$, $q = p^k$ will be congruent to 1 modulo $n-1$. By Lemma 3.2.3, there is a PL normalizer $h$ of $B_n(\mathbf{R})$ that fixes 0, that commutes with $t_{n-1}$ and that has slopes in $q\langle n\rangle$. The normalizer $h$ is not in $B_n(\mathbf{R})$ since $q$ is not in $\langle n\rangle$. However, it is easy to show that $h\pi_n$ is trivial. If for some $i$ there were a $w \in C^0(F_{n,\infty}, t_{n-1})$ so that $wh^{-1}g_{n,i}hw^{-1} = g_{n,i}$, then the PL normalizer $hw^{-1}$ would be non trivial and commute with $g_{n,i}$ contradicting Proposition 5.3.1. This shows that for our assumed $n$, the intersection of $\mathrm{Aut}_{PL}(F_{n,\infty}, t_{n-1})$ and $\mathrm{Aut}^0_{\pi}(F_{n,\infty}, t_{n-1})$ does not lie in $\mathrm{Aut}^0_{LI}(F_{n,\infty}, t_{n-1})$.

The next lemma shows how lifting affects an automorphism that is inner (at some $i$). We use $a^b$ to represent $b^{-1}ab$ or the action of $b$ on $a$ as appropriate.

**Lemma 5.3.5.** *Let $2 \le m < n$ be integers, and let $\tau_{m,n}$ and $\Theta_{m,n}$ be as in Theorem 4.1.6. Let $W$ be in $F_{m,[0,m-1]}$, let $W' = W\tau_{m,n}$, and let $w = W\Phi_m \in C^0(B_m(\mathbf{R}), t_{m-1})$ be regarded as an element of $\mathrm{Aut}^0(F_{m,\infty}, t_{m-1})$.*

(a) Then $W'$ is in $F_{n,[0,m-1]} \subseteq F_{n,[0,n-1]}$ and if $u = w\Theta_{m,n} \in \text{Aut}^0(F_{n,\infty}, t_{n-1})$ and $v = W'\Phi_n \in C^0(B_n(\mathbf{R}), t_{n-1}) \subseteq \text{Aut}^0(F_{n,\infty}, t_{n-1})$, then for all $i$ with $0 \leq i \leq m-2$, we have $(g_{n,i})^u = (g_{n,i})^v$.

(b) Now let $i$ be an integer with $0 \leq i \leq m-2$. Assume that $\alpha \in \text{Aut}^0(F_{m,\infty}, t_{m-1})$ is inner at $i$ via $w$. Let $\beta = \alpha\Theta_{m,n}$, and let $v = W'\Phi_n$. Then $(g_{n,i})^\beta = (g_{n,i})^v$ and $\beta$ is inner at $i$ via $v$.

*Proof.* We get (b) from (a) and Corollary 4.1.6.1 because $(g_{n,i})^\beta = (g_{n,i})^u = (g_{n,i})^v$ with $u$ as in (a). We thus concentrate on (a).

Let $W_1 = W\sigma_{m-1}$. That $W'$ is in $F_{n,[0,m-1]}$ follows from Lemma 4.1.5. Let

$$W_1' = W'\sigma_{n-1} = t_{n-1}^{-1}W't_{n-1} = (t_{m-1}^{-1}Wt_{n-1})\tau_{m,n} = (W\sigma_{m-1})\tau_{m,n} = W_1\tau_{m,n}.$$

Now two applications of Lemma 3.1.2(b) and one of Theorem 4.1.6 give

$$(g_{n,i})^v = (W')^{-1}g_{n,i}W'W_1' = (W^{-1}g_{m,i}WW_1)\tau_{m,n} = ((g_{m,i})^w)\tau_{m,n} = (g_{n,i})^u.$$

$\square$

**Corollary 5.3.5.1.** *Let $m < n$, let $\alpha$ be in $Aut^0(F_{m,\infty}, t_{m-1})$ and let $\beta = \alpha\Theta_{m,n} \in Aut^0(F_{n,\infty}, t_{n-1})$. If $\alpha$ is locally inner, then so is $\beta$.*

*Proof.* Lemma 5.3.5 works for each $i$ with $0 \leq i < m-1$. For $m-1 \leq i < n-1$, let $w_i = 1$. $\square$

Now that we know that locally inner automorphisms are closed under lifting, we can define another subgroup of $\text{Aut}^0(F_{n,\infty}, t_{n-1})$. As $n$ ranges over the integers greater than 1, the groups $\text{Aut}^0(F_{n,\infty}, t_{n-1})$ form a graded group that we denote by $\boldsymbol{\mathcal{A}}$. This graded group is closed under lifting. We let $\boldsymbol{\mathcal{LC}}$ be the smallest graded subgroup of $\boldsymbol{\mathcal{A}}$ that contains all the $C^0(B_n(\mathbf{R}), t_{n-1})$ and that is closed under lifts. The portion of $\boldsymbol{\mathcal{LC}}$ in $\text{Aut}^0(F_{n,\infty}, t_{n-1})$ will be denoted $LC^0(B_n(\mathbf{R}), t_{n-1})$. The next corollary follows from Corollary 5.3.5.1.

**Corollary 5.3.5.2.** *$LC^0(B_n(\mathbf{R}), t_{n-1})$ is a subgroup of $\text{Aut}_{LI}^0(F_{n,\infty}, t_{n-1})$.* $\square$

5.4. **Normalizers of $A_n$.** In this section, we obtain our main results from the machinery developed in Section 5.3. The results here will prove that certain normalizers of $A_n$ are not PL and that certain normalizers of $B_n$ are not normalzers of $A_n$.

Let $W \in F_{n,\infty}$ and $0 \leq j < n-1$. We say that $W$ *avoids (the residue class of) $j$* if $W$ can be written as a product of generators in Definition 2.1.1 and their inverses that uses no $g_i$ with $i$ congruent to $j$ modulo $n-1$. Let $w \in C^0(B_n(\mathbf{R}), t_{n-1})$. We say that $w$ *avoids (the residue class of) $j$* if we can write $w = W\Phi_n$ for some $W \in F_{n,[0,n-1]}$ such that $W$ avoids $j$.

Note that if $W$ avoids $j$ then it can be written in seminormal form without using generators whose subscripts are congruent to $j$ (mod $n-1$). This is due to the fact that the rewriting rules that lead to the seminormal form preserve the residue classes of the subscripts.

**Lemma 5.4.1.** *(a) $\{W \in F_{n,\infty} | W$ avoids $j\}$ is a subgroup of $F_{n,\infty}$.*
*(b) $\{w \in C^0(B_n(\mathbf{R}), t_{n-1}) | w$ avoids $j\}$ is a subgroup of $C^0(B_n(\mathbf{R}), t_{n-1})$.*
*(c) Let $m < n$, $W \in F_{m,\infty}$ and $\tau_{m,n} : (F_{m,\infty}, t_{m-1}) \to (F_{n,\infty}, t_{n-1})$ be the monomorphism of Proposition 4.1.2. Then $W\tau_{m,n}$ avoids every $j$ with $m-1 \leq j < n-1$.*

*Proof.* (a) It is immediate from the definition.
(b) This follows from (a) and Lemma 3.1.2(a).
(c) This is immediate from the definition of $\tau_{m,n}$ in Proposition 4.1.2. $\blacksquare$

**Lemma 5.4.2.** *Let $W \in F_{n,\infty}$, and $0 \leq j < n-1$. Then $W$ avoids $j$ if and only if $W$ can be represented by two allowable subdivisions in the sense of Lemma 2.2.2 such that for every $l \geq 0$, the intervals $\left[\frac{a}{n^l}, \frac{a+1}{n^l}\right]$, $a \equiv j \pmod{n-1}$, are not subdivided.*

*Proof.* Note that when an allowable subdivision is constructetd as in Definition 2.3.3, the number of an interval of the form $\left[\frac{a}{n^l}, \frac{a+1}{n^l}\right]$ is congruent to $a \pmod{n-1}$, so each generator $g_i$ in the seminormal form of $W$ will partition an interval $\left[\frac{a}{n^l}, \frac{a+1}{n^l}\right]$ with $a \equiv i \pmod{n-1}$. $\blacksquare$

**Corollary 5.4.2.1.** *Let $a < b < c$ be in $\mathbf{Z}$, $0 \leq j < n-1$, and $W_1, W_2 \in F_{n,\infty}$ with $\mathrm{supp}(W_1) \subseteq [a,b]$, $\mathrm{supp}(W_2) \subseteq [b,c]$. If $W_1 W_2$ avoids $j$, then so do $W_1$ and $W_2$.*

*Proof.* Write $W_1 W_2 = P N^{-1}$ in seminormal form so that the allowable subdivisions $D$ and $E$ corresponding to $P$ and $N$ respectively, do not subdivide intervals of the form $\left[\frac{a}{n^l}, \frac{a+1}{n^l}\right]$ with $a \equiv j \pmod{n-1}$. Combining the portions of $D$ and $E$ that lie in $(-\infty, b]$ with the standard subdivision on $[b, \infty)$ yields a pair of allowable subdivisions for $W_1$. Similarly, the portions of $D$ and $E$ in $[b, \infty)$ combined with the standard subdivision on $(-\infty, b]$ yield a pair of allowable subdivisions for $W_2$. Therefore $W_1$ and $W_2$ avoid $j$. $\blacksquare$

**Lemma 5.4.3.** *Let $b \in \mathbf{Z}$, $0 \leq j < n-1$, and $w \in C^0(B_n(\mathbf{R}), t_{n-1})$ be such that $w$ fixes $b$. Then $t_b w t_b^{-1} \in C^0(B_n(\mathbf{R}), t_{n-1})$ and $w$ avoids $j$ iff $t_b w t_b^{-1}$ avoids $j - b$.*

*Proof.* Since $w$ commutes with $t_{n-1}$ we may assume that $0 \leq b < n-1$. The fact that $t_b w t_b^{-1} \in C^0(B_n(\mathbf{R}), t_{n-1})$ is immediate. For the rest, it suffices to do one direction. Assume $w$ avoids $j$. Write $w = W \Phi_n$ with $W \in F_{n,[0,n-1]}$, so that $W$ avoids $j$. Then we have that $W$ fixes $b$. So write $W = W_1 W_2$ where $\mathrm{supp}(W_1) \subseteq [0,b]$ and $\mathrm{supp}(W_2) \subseteq [b, n-1]$. By Corollary 5.4.2.1 both $W_1$ and $W_2$ avoid $j$. Let $U_1 = t_b W_2 t_b^{-1}$ and $U_2 = t_{n-1-b}^{-1} W_1 t_{n-1-b}$, and $U = U_1 U_2$. Clearly both $U_1$ and $U_2$ avoid $j - b$ so by Lemma 5.4.1(a), $U$ avoids $j - b$. Since we have $t_b w t_b^{-1} = U \Phi_n$ the conclusion is that $t_b w t_b^{-1}$ avoids $j - b$. $\blacksquare$

**Lemma 5.4.4.** *For $n \geq 3$, let $\alpha$ be in $Aut^0(F_{n,\infty}, t_{n-1})$, and let $h$ be the self-homeomorphism of $\mathbf{R}$ that fixes $0$, commutes with $t_{n-1}$ and realizes $\alpha$. Let $0 \leq i < n-1$ and $\gamma_i, \gamma_{i+1} \in Aut_{PL}(B_n(\mathbf{R}), t_{n-1})$ be conjugation by PL homeomorphisms $w_i$ and $w_{i+1}$ such that $(g_i)\alpha = (g_i)\gamma_i$ and $(g_{i+1})\alpha = (g_{i+1})\gamma_{i+1}$. If $h$ normalizes $A_n(\mathbf{R})$, then $h^{-1} t_1 h = w_i^{-1} t_1 w_{i+1}$.*

*Proof.* Observe that $(g_i)\gamma_i \alpha^{-1} \sigma_1 \alpha \gamma_{i+1}^{-1} \sigma_1^{-1} = g_i$. Since $h$ normalizes $A_n(\mathbf{R})$, both $h^{-1} t_1 h$ and $w_i h^{-1} t_1 h w_{i+1}^{-1} t_1^{-1}$ are PL. Applying Proposition 5.3.1 to this PL normalizer of $B_n(\mathbf{R})$ which commutes with $t_{n-1}$ we get $w_i h^{-1} t_1 h w_{i+1}^{-1} t_1^{-1} = 1$, so $h^{-1} t_1 h = w_i^{-1} t_1 w_{i+1}$. $\blacksquare$

**Proposition 5.4.5** (Double unbent generator). *Let $\alpha \in Aut^0(F_{n,\infty}, t_{n-1})$, $n \geq 3$, be a normalizer of $A_n(\mathbf{R})$. If for some $0 \leq i < n-1$, $\alpha$ has an unbent generator at $i$ and $i+1$, then $\alpha$ is the identity.*

*Proof.* Let $h$ be the self-homeomorphism of $\mathbf{R}$ that fixes 0, commutes with $t_{n-1}$ and realizes $\alpha$. By Corollary 5.2.2.1, $h^{-1}g_i h$ has slope $n$ to the right of $ih$, and being affine on $[i, i+1]h$, forces this interval to have width 1. Similarly for $[i+1, i+2]h$, and the two intervals are consecutive. Let $b = ih - i = (i+1)h - (i+1) = (i+2)h - (i+2)$. Then $h^{-1}g_i h = t_b^{-1}g_i t_b$ and $h^{-1}g_{i+1}h = t_b^{-1}g_{i+1}t_b$. By Lemma 5.4.4, $h^{-1}t_1 h = t_1$, so $h$ commutes with $t_1$. It follows that $h^{-1}g_j h = t_b^{-1}g_j t_b$ for every $0 \leq j < n-1$, so $\alpha$ is conjugation by $t_b$. Since $h$ fixes 0, we must have $b = 0$ so $\alpha \equiv 1$. $\qquad\square$

**Lemma 5.4.6.** *For $n \geq 3$, let $\alpha \in Aut^0(F_{n-1,\infty}, t_{n-2})$ be locally inner. If $\beta = \alpha\Theta_{n-1,n}$ normalizes $A_n(\mathbf{R})$ then $\alpha$ is the identity.*

*Proof.* By Lemmas 5.4.1(c) and 5.3.5, for each $i$ there are $w_i \in C^0(B_n(\mathbf{R}), t_{n-1})$ such that $g_i\beta = w_i^{-1}g_i w_i$ and $w_i$ avoids $n-2$. Moreover $w_{n-2} = 1$ and $w_{i+n-1} = w_i$. It follows from Lemma 5.4.4 that $h^{-1}t_1 h = w_i^{-1}t_1 w_{i+1}$, so we get $h^{-1}t_1 h = t_1 w_0 = w_0^{-1}t_1 w_1 = \cdots = w_{n-4}^{-1}t_1 w_{n-3} = w_{n-3}^{-1}t_1$ so

$$(*) \qquad\qquad w_{i+1} = t_1^{-1}w_i t_1 w_0 \quad \text{for} \quad 0 \leq i \leq n-4,$$

$$(**) \qquad\qquad 1 = t_1^{-1}w_{n-3}t_1 w_0.$$

From $(*)$ and the fact that $w_0$ and $w_{i+1}$ fix 0, we get that $w_i$ fixes $n-2$ when $0 \leq i \leq n-4$. For $w_{n-3}$, use $(**)$ and we have that all $w_i$ fix $n-2$.

Now from $(**)$, Lemmas 5.4.1(b) and 5.4.3, and the fact that $w_{n-3}$ avoids $n-2$, we get that $w_0$ avoids 0. From this, $(*)$, Lemmas 5.4.1(b) and 5.4.3, and the fact that all $w_i$ avoid $n-2$, we get that all $w_i$ avoid 0. Once again, iterate the process to get that all $w_i$ avoid each of $0, 1, \ldots, n-2$. Therefore we must have all $w_i = 1$ and $\alpha \equiv 1$. $\qquad\square$

**Theorem 5.4.7.** *For $2 \leq m < n$, let $\alpha \in Aut^0(F_{m,\infty}, t_{m-1})$. For $m = n-1$ assume that $\alpha$ is locally inner. If $\beta = \alpha\Theta_{m,n}$ normalizes $A_n(\mathbf{R})$ then $\alpha$ is the identity.*

*Proof.* This follows from Lemma 5.4.6, when $m = n-1$, and from Proposition 5.4.5, when $m < n-1$. $\qquad\square$

The next application of Proposition 5.4.5 justifies the remark made after Theorem 4.3.1.

**Theorem 5.4.8.** *Assume $(m-1) \mid (n-1)$ with $m > 2$. Let $w \in C^0(B_m(\mathbf{R}), t_{m-1})$ represent a non-trivial automorphism $\alpha$ of $A_m(\mathbf{R})$ by conjugation. Let $\beta = w\Lambda_{m,n}$ as described in Theorem 4.3.1. Then $\beta$ does not normalize $A_n(\mathbf{R})$.*

*Proof.* Assume that $\beta$ normalizes $A_n(\mathbf{R})$. We have $w = W\Phi_m$ for some $W \in F_{m,[0,m-1]}$. Now if $W' = W\tau_{m,n}$, then by Lemma 4.1.5, we have $W' \in F_{n,[0,m-1]} \subseteq F_{n,[0,n-1]}$. If $w' = W'\Phi_n$, then Lemma 5.3.5 (which applies since $\beta$ and $\alpha\Theta_{m,n}$ agree on the relevant generators) gives $g_{n,i}\beta = (w')^{-1}g_{n,i}w'$ for $0 \leq i \leq (m-2)$. Now $w'$ avoids $j$ unless $0 \leq j\phi_n \leq (m-2)$. Thus if $w'$ were to commute with $t_{m-1}$, then $w'$ would be trivial. Since $\beta$ fixes $t_{m-1}$, the self homeomorphism $h$ of $\mathbf{R}$ for which $g_{n,j}\beta = h^{-1}g_{n,j}h$ for all $j$ must commute with $t_{m-1}$ as well. Thus $h \neq w'$. Now the composition $(w')^{-1}h$ is non-trivial and normalizes $A_n(\mathbf{R})$ since $w' \in A_n(\mathbf{R})$. However, conjugation by $(w')^{-1}h$ leaves each $g_{n,i}$ fixed with $0 \leq i \leq (m-2)$ which contradicts Proposition 5.4.5. $\qquad\square$

The next theorem will be used to show that there are many non-PL elements in $Aut(A_n(\mathbf{R}))$.

**Theorem 5.4.9.** *Let* $2 \leq m < n$ *be integers for which* $(m-1) \mid (n-1)$, *and consider* $\alpha \in \mathrm{Aut}^0_{LI}(F_{m,\infty}, t_{m-1})$. *Let* $\beta = \alpha \Lambda_{m,n}$ *as described in Theorem 4.3.1. If* $\alpha$ *is not trivial, then* $\beta$ *is not PL.*

*Proof.* From the construction in Theorem 4.3.1, we know that $\beta = \alpha \Theta_{m,n}$ on those $g_{n,i}$ with $0 \leq i \leq m-2$. By Lemma 5.3.5, $\beta$ is inner at $i$ via some $w_i$ in $C^0(B_n(\mathbf{R}), t_{n-1})$ for all $i$ with $0 \leq i \leq m-2$. For any $j$, $g_{n,j} = t_{m-1}^{-k} g_{n,i} t_{m-1}^k$ for some unique $k$ and $i$ with $0 \leq i \leq m-2$. An elementary calculation using the fact that $\beta$ commutes with $t_{m-1}$ shows that $\beta$ is inner at $j$ via $w_j$ where $w_j$ is equal to the conjugate of $w_i$ by $t_{m-1}^k$.

By Lemma 5.3.3(e), we are done if we show that $\beta$ is not inner. By Corollary 5.3.3.1, we are done if we show that $\beta$ is inner at two values via unequal elements of $C^0(B_n(\mathbf{R}), t_{n-1})$.

Consider an $i$ with $0 \leq i \leq m-2$. Since $\beta = \alpha \Theta_{m,n}$ on $g_{n,i}$, Lemma 5.4.1(c) says that $w_i$ avoids all $j$ with $m-1 \leq j < n-1$. Now $\beta$ is inner at $i+m-1$ by $w_i' = t_{m-1}^{-1} w_i t_{m-1}$. By Lemma 5.4.3, $w_i'$ avoids all $j$ with $2(m-1) \leq j < (n-1)+(m-1)$ modulo $(n-1)$. However, this includes all $j$ with $0 \leq j < m-1$ modulo $(n-1)$. If we were to have $w_i = w_i'$, we would have $w_i = 1$. This cannot happen for all $i$ with $0 \leq i \leq m-2$ since $\alpha$ is not trivial. This competes the proof. $\square$

**5.5. An example of torsion.** In the following, we leave to the reader the task of checking a large number of calculations.

**Theorem 5.5.1.** *For* $n \geq 5$ *there is an element of order* $n-2$ *in* $\mathrm{Out}_o(B_n(\mathbf{R}))$.

*Proof.* The example will be a rotation of a lift of an element in $\mathrm{Aut}^0(F_{4,\infty}, t_3)$. Consider the function fixing $t_3$ and sending

$$g_{4,0} \mapsto g_{4,0} g_{4,4} g_{4,2}^{-1},$$
$$g_{4,1} \mapsto g_{4,2} g_{4,3} g_{4,5}^{-1},$$
$$g_{4,2} \mapsto g_{4,2} g_{4,4} g_{4,2}^{-1}.$$

This extends to a function defined on all the $g_{4,i}$ that fixes $t_3$. Using the subdivisions of Section 2.2, it is easy to show that the the images of the $g_{4,i}$, $i \in \{0,1,2\}$, have their active parts confined to $[0,2]$, $[2, 2\frac{1}{2}]$ and $[2\frac{1}{2}, 3]$ in order, and that they satisfy the relations defining $F_{4,\infty}$. Thus the function extends to an endomorphism $\alpha$ of $F_{4,\infty}$ that commutes with $\sigma_3$. It is straightforward to check that the inverse of $\alpha$ is defined by

$$g_{4,0} \mapsto g_{4,0} g_{4,1} g_{4,3}^{-1},$$
$$g_{4,1} \mapsto g_{4,0} g_{4,2} g_{4,0}^{-1},$$
$$g_{4,2} \mapsto g_{4,1} g_{4,5} g_{4,3}^{-1}.$$

It is not needed here, but $\alpha$ is the automorphism used in Theorem 4.4.3 that is conjugation by a PL function with slopes in $2\langle 4 \rangle$ and thus $\alpha$ is not in $C^0(F_{4,\infty}, t_3)$.

If we lift $\alpha$ to $\beta_n$ by $\Theta_{4,n}$ for $n > 4$, then we get an automorphism defined by

$$g_{n,0} \mapsto g_{n,0} g_{n,n} g_{n,2}^{-1},$$
$$g_{n,1} \mapsto g_{n,2} g_{n,n-1} g_{n,n+1}^{-1},$$
$$g_{n,2} \mapsto g_{n,2} g_{n,n} g_{n,2}^{-1},$$
$$g_{n,i} \mapsto g_{n,i} \qquad \text{for } 3 \leq i \leq n-2.$$

If $h$ is the self homeomorphism of $\mathbf{R}$ that realizes $\beta_n$ by conjugation, then the 1-step rotation $\gamma_n$ of $\beta_n$ is realized as conjugation by $k = t_1 h t_{1h}^{-1}$. However, $1h = 2$ since it can be checked that the leftmost break of $g_{n,1}\beta_n$ is 2. Thus $\gamma_n$ is conjugation by $t_1 h t_2^{-1}$ and $\gamma_n = \sigma_1 \beta_n \sigma_2^{-1}$. This gives that $\gamma_n$ is determined by

$$
\begin{aligned}
g_{n,0} &\mapsto g_{n,0} g_{n,n-3} g_{n,n-1}^{-1}, \\
g_{n,1} &\mapsto g_{n,0} g_{n,n-2} g_{n,0}^{-1}, \\
g_{n,i} &\mapsto g_{n,i-1} \qquad \text{for } 2 \le i \le n-3, \\
g_{n,n-2} &\mapsto g_{n,n-3} g_{n,2n-3} g_{n,n-1}^{-1} \\
&= g_{n,n-2} g_{n,n-3} g_{n,n-1}^{-1}.
\end{aligned}
$$

where the last expression that is not in seminormal form will sometimes be convenient.

Let

$$
\begin{aligned}
u_{n,k} &= g_{n,0} g_{n,n-k} g_{n,n-1}^{-1}, \text{ and} \\
v_{n,k} &= g_{n,n-2} g_{n,n-k} g_{n,n-1}^{-1}.
\end{aligned}
$$

Note that $g_{n,0}\gamma_n = u_{n,3}$ and $g_{n,n-2}\gamma_n = v_{n,3}$. We claim that $u_{n,k}\gamma_n = u_{n,k+1}$ and $v_{n,k}\gamma_n = v_{n,k+1}$ for $3 \le k \le n-2$. In this range, $2 \le n-k \le n-3$ and

$$
\begin{aligned}
u_{n,k}\gamma_n &= g_{n,0} g_{n,n-3} g_{n,n-1}^{-1} g_{n,n-(k+1)} g_{n,2n-2} g_{n,2n-4}^{-1} g_{n,n-1}^{-1} \\
&= g_{n,0} g_{n,n-(k+1)} g_{n,n-1}^{-1} \\
&= u_{n,k+1}
\end{aligned}
$$

where the second equality holds because $n-(k+1)$ is less than both $n-1$ and $n-3$. The claim for $v_{n,k}\gamma_n$ has an almost identical proof.

Now $g_{n,n-3}\gamma_n^j = g_{n-3-j}$ for $1 \le j \le n-4$. Thus, for these values of $j$, $\gamma_n^j$ is non-trivial and satisfies the hypotheses of Proposition 5.3.1. For these values of $j$ we know that $\gamma_n^j$ is not represented as conjugation by a PL self homeomorphism of $\mathbf{R}$ and represents a non-trivial element in $\mathrm{Out}_o(B_n(\mathbf{R}))$. When we show that $\gamma_n^{n-2}$ represents a trivial element in $\mathrm{Out}_o(B_n(\mathbf{R}))$, we will know that $\gamma_n^{n-3}$ represents a non-trivial element in $\mathrm{Out}_o(B_n(\mathbf{R}))$ as well.

Now

$$
\begin{aligned}
g_{n,0}\gamma_n^{n-2} &= u_{n,n-1}\gamma_n \\
&= (g_{n,0} g_{n,1} g_{n,n-1}^{-1})\gamma_n \\
&= g_{n,0}^2 g_{n,n-2} g_{n,2n-2}^{-1} g_{n,0}^{-1}
\end{aligned}
$$

where the last equality follows from the definition of $\gamma_n$ and several applications of the relations of $F_{n,\infty}$. A similar calculation gives

$$
g_{n,n-2}\gamma_n^{n-2} = v_{n,n-1}\gamma_n = g_{n,0} g_{n,n-2} g_{n,3n-4} g_{n,2n-2}^{-1} g_{n,0}^{-1}.
$$

To help work with the other generators, we let $P = g_{n,n-2} g_{n,0}^{-1}$. An easy check shows $P\gamma_n = P$. Now $g_{n,1}\gamma_n = g_{n,0}P$ so for $1 \le i \le n-3$ we have

$$
\begin{aligned}
g_{n,i}\gamma_n^{n-2} &= (g_{n,0}P)\gamma_n^{n-2-i} = u_{n,n-i}P \\
&= g_{n,0} g_{n,i} g_{n,n-1}^{-1} g_{n,n-2} g_{n,0}^{-1} \\
&= g_{n,0} g_{n,i} g_{n,n-2} g_{n,2n-2}^{-1} g_{n,0}^{-1}.
\end{aligned}
$$

We will show that $\gamma_n^{n-2}$ is trivial in $\mathrm{Out}_o(B_n(\mathbf{R}))$ by showing that the action of $\gamma_n^{n-2}$ is realized as conjugation by an element of $C^0(F_{n,\infty}, t_{n-1})$. We claim that the action of $\gamma_n^{n-2}$ is conjugation by $w = P\Phi_n$. Using subdivisions, it is easy to check that $P$ has support in $[0, n-1]$ so that $P\Phi_n$ makes sense. By Lemma 3.1.2(b), $w^{-1}g_{n,i}w = P^{-1}g_{n,i}PP_1$ where $P_1 = P\sigma_{n-1} = g_{n,2n-3} g_{n,n-1}^{-1}$. It is now a straightforward exercise in using the relations of $F_{n,\infty}$ to check that for $0 \le i \le n-2$ we get $w^{-1}g_{n,i}w = g_{n,i}\gamma_n^{n-2}$. This completes the proof. ∎

The lemma above holds for $n = 4$ as well. The example is represented by the automorphism $\alpha$ in the proof. To verify this, one must actually prove the remark in the proof that the action of $\alpha$ is conjugation by a PL function with slopes in $2\langle 4 \rangle$. The order of this outer automorphism is 2. What is surprising is that the finiteness of the order persists in the particular combination of lift and rotation used. In fact, all $j$-step rotations of $\beta_n$, $n > 4$, by integers with $\phi_n(j) \ne 1$ and $\phi_n(j) \ne 2$ represent outer automorphisms of infinite order. This follows from Proposition 5.3.1 since each such rotation fixes $g_{n,0}$ and is non-trivial (and thus its positive powers are non-trivial). We do not know the status of the 2-step rotation of $\beta_n$, although numerical calculations hint that the order is infinite.

## 6. On the structure of $\mathrm{Aut}(A_n)$ and $\mathrm{Aut}(B_n)$

In this section we collect and comment on the main results of the paper. Section 6.1 gives general results about the automorphisms of $A_n$ and $B_n$ for a single value of $n$. Section 6.2 discusses local properties of automorphisms, i.e., the action on a single generator. Section 6.3 discusses the relations between automorphims of $A_m$ and $A_n$, and $B_m$ and $B_n$ for values $2 \le m < n$. The results are not arranged to show dependency. All results can be proven using material prior to this section, but in the order that they are presented here, some results in Section 6.1 depend on results quoted in Sections 6.2 and 6.3.

6.1. **Results for a single $n$.** The classes $[A_n]$ and $[B_n]$, $n \ge 2$, of homeomorphism groups on both the real line, and the circle are defined in Section 1.3. For $C = A$ or $C = B$, the class $[C_n]$ has a container group $C_n$. The McCleary-Rubin theorem applies to all groups in the classes, automorphisms can be viewed as normalizers, and the characteristic subgroups of $C_n$ are precisely the ones with the same automorphism group as $C_n$.

The next two lemmas are minor rephrasing of Lemmas 1.3.1, 1.3.2 and 2.2.7.

**Lemma 6.1.1.** *Let $X = \mathbf{R}$ or $X = S_r$, $C = A$ or $C = B$, and $r \in \mathbf{Z}[\frac{1}{n}]$ with $r \in \Delta_n$ for the case $C = B$. Then the following hold.*

(a) *The groups $C_n(\mathbf{R})$, $C_n(S_r)$, are in the class $[C_n]$.*
(b) *For $G$ in the class $[C_n]$, $N(G) \simeq \mathrm{Aut}(G)$.*
(c) *If $G$ is in $[C_n]$, we have $N(G) \subseteq N(C_n(X))$.*
(d) *$G \in [C_n]$ is characteristic in $C_n(X)$ if and only if $N(G) = N(C_n(X))$.*

Each of the two classes contains at least one of the Thompson groups.

**Lemma 6.1.2.** *The group $F_n$ is in $[A_n]$, and the group $F_{n,\infty}$, is in $[B_n]$.*

Moreover, in each case groups of automorphisms of the container and the corresponding Thompson group are closely related.

**Proposition 6.1.3.** *There are exact sequences and isomorphisms*

(1) $$1 \to \mathrm{Aut}_B(F_{n,\infty}) \to \mathrm{Aut}_o(F_{n,\infty}) \to \mathrm{Out}_o(B_n(\mathbf{R})) \to 1,$$

(2) $$1 \to \mathrm{Aut}_B^0(F_n) \to \mathrm{Aut}_\pi^0(F_n) \to \mathrm{Out}_\pi(A_n(\mathbf{R})) \to 1,$$

(3) $$\mathrm{Out}_\pi(B_n(\mathbf{R})) \simeq \mathrm{Aut}_\pi^0(F_{n,\infty}, t_{n-1}) \Big/ C^0(B_n(\mathbf{R}), t_{n-1}),$$

(4) $$\mathrm{Out}_\pi(A_n(\mathbf{R})) \simeq \mathrm{Aut}_\pi^0(F_n, t_1) \Big/ C^0(B_n(\mathbf{R}), t_1),$$

(5) $$1 \longrightarrow F_{n,\infty}^0 \longrightarrow \mathrm{Aut}_B^0(F_n) \overset{\theta_A}{\longrightarrow} T_{n,1} \times T_{n,1} \longrightarrow 1 \text{ , and}$$

(6) $$1 \longrightarrow F_{n,\infty} \longrightarrow \mathrm{Aut}_B(F_{n,\infty}) \overset{\theta_B}{\longrightarrow} T_{n,n-1} \times T_{n,n-1} \overset{d}{\longrightarrow} \mathbf{Z}_{n-1} \longrightarrow 1.$$

Note that the middle terms of (5) and (6) are the second terms of (2) and (1).

*Proof.* Items (1–4) are given in Proposition 2.6.2. In (5), the restriction that elements of $F_n$ must be integer translations near $\pm\infty$ forces an element $f$ of $\mathrm{Aut}_B^0(F_n)$ to satisfy $(x+1)f = xf + 1$ near $\pm\infty$. We thus take $f\theta_A$ to be the pair of germs of $f$ near $\pm\infty$ modulo integer translations. This gives a well defined pair of elements of $T_{n,1}$. If such a pair of elements is given, then they lift to a pair $(f_-, f_+)$ where each is in $\mathrm{Aut}(F_n)$. Since the lifts are adjustable by composing with translations by integers, we can arrange that $f_-\rho_n = f_+\rho_n = 0$ and that $xf_+ > x$ and $xf_- < x$ for all $x$. Then Proposition 1.2.2 lets us build an element $f$ of $A_n(\mathbf{R})$ that agrees with $f_\pm$ near $\pm\infty$ and that fixes 0. This shows that $\theta_A$ is onto. The kernel of $\theta_A$ consists of all elements of $B_n(\mathbf{R})$ that fix 0 and are integer translations near $\pm\infty$. Since residues are preserved in $B_n(\mathbf{R})$, the translations are by multiples of of $n-1$ and the kernel consists of all elements of $F_{n,\infty}$ that fix 0.

The arguments for (6) are similar to those for (5) except that we cannot adjust $f_-\rho_n$ and $f_+\rho_n$. In fact there is a well defined homomorphism $\rho_n : T_{n,n-1} \to \mathbf{Z}_{n-1}$ giving the amount of shift on residue classes induced by an element of $T_{n,n-1}$. The image of $\theta_B$ is the kernel of the map $d : T_{n,n-1} \times T_{n,n-1} \to \mathbf{Z}_{n-1}$ defined by $(f_1, f_2)d = f_1\rho_n - f_2\rho_n$. $\qquad\square$

The two groups $A_n(\mathbf{R})$ and $B_n(\mathbf{R})$ have the same PL normalizers which are completely determined by Lemma 3.2.2. From Lemma 3.2.2 and Corollary 3.2.4.1 we have:

**Theorem 6.1.4.** *Let $n \geq 2$.*

(a) $$\mathrm{Aut}_{PL}(A_n(\mathbf{R})) = \mathrm{Aut}_{PL}(B_n(\mathbf{R}))$$

(b) $$\mathrm{Out}_{PL}(A_n(\mathbf{R})) \simeq \langle\!\langle n \rangle\!\rangle / \langle n \rangle$$

(c) $$\mathrm{Out}_{PL}(B_n(\mathbf{R})) \simeq \mathbf{Z}_{n-1} \rtimes \langle\!\langle n \rangle\!\rangle / \langle n \rangle.$$

See Corollary 3.2.4.2 for the analysis of torsion that occurs in $\mathrm{Out}_{PL}(A_n(\mathbf{R}))$ and $\mathrm{Out}_{PL}(B_n(\mathbf{R}))$.

From Theorem 5.4.9, we get infinitely many non-PL automorphisms of $A_n$.

**Theorem 6.1.5.** *For $3 \leq n$, $\mathrm{Aut}(A_n(\mathbf{R}))$ contains an isomorphic copy $\widehat{F}$ of $F \simeq F_2 \simeq F_{2,\infty}$. The intersection of $\widehat{F}$ with $\mathrm{Aut}_{PL}(A_n(\mathbf{R}))$ is trivial.*

*Proof.* Theorem 4.3.1 says that $\Lambda_{2,n}$ injects $\mathrm{Aut}^0(F_{2,\infty}, t_1)$ into $\mathrm{Aut}^0(F_{n,\infty}, t_1)$. Since $F_{n,\infty}$ and $t_1$ generate $F_n$, $\Lambda_{2,n}$ injects $\mathrm{Aut}^0(F_{2,\infty}, t_1)$ into $\mathrm{Aut}^0(F_n, t_1)$. Since automorphisms of $F_n$ are automorphisms of $A_n(\mathbf{R})$, $\Lambda_{2,n}$ injects $\mathrm{Aut}^0(F_{2,\infty}, t_1)$ into $\mathrm{Aut}^0(A_n(\mathbf{R}), t_1)$. From Theorem 5.4.9, we get that the image of the non-trivial locally inner elements in $\mathrm{Aut}^0(F_{2,\infty}, t_1)$ are not PL. However, elements of $\mathrm{Aut}^0(F_{2,\infty}, t_1)$ are orientation preserving, and it follows from the main theorem of [2] that all elements of $\mathrm{Aut}^0(F_{2,\infty}, t_1)$ are in $C^0(B_2(\mathbf{R}), t_1)$ and are thus (locally) inner. Further $C^0(B_2(\mathbf{R}), t_1)$ is isomorphic to $F_{0,[0,1]}$ which is isomorphic to each of $F_{2,\infty}$, $F_2$ and Thompson's original group $F$. □

The two container groups and their automorphism groups present some unusual behavior in the sense that $A_n(\mathbf{R})$ and $B_n(\mathbf{R})$ "almost equal" whereas $\mathrm{Aut}(A_n(\mathbf{R}))$ and $\mathrm{Aut}(B_n(\mathbf{R}))$ are "very different". More precisely:

**Theorem 6.1.6.** *Let $n > 2$, $A = A_n(\mathbf{R})$, $B = B_n(\mathbf{R})$. The following properties hold:*

(a) *$B$ is characteristic in $A$ so that restriction to $B$ gives a homomorphism from $\mathrm{Aut}(A)$ to $\mathrm{Aut}(B)$.*

(b) *The homomorphism from $\mathrm{Aut}(A)$ to $\mathrm{Aut}(B)$ is a monomorphism so that we can regard $\mathrm{Aut}(A)$ as a subgroup of $\mathrm{Aut}(B)$.*

(c) *The index of $B$ in $A$ is finite.*

(d) *The index of $\mathrm{Aut}(A)$ in $\mathrm{Aut}(B)$ is infinite.*

*Proof.* The first three parts follow from Proposition 1.4.2 and Lemma 1.1.4. For part (d), Theorem 5.4.7 implies that $\mathrm{Aut}(A)$ has trivial intersection with the infinite groups $\langle \mathrm{Aut}^0(F_{m,\infty}, t_{m-1}) \rangle \Theta_{m,n}$ for $m \leq n - 2$ and $\langle \mathrm{Aut}^0_{LI}(F_{n-1,\infty}, t_{n-2}) \rangle \Theta_{n-1,n}$. □

The behavior in Theorem 6.1.6 although unusual, is found in other places. The following example is a variation of one supplied to us by F. Gross [15]. Note that while there are other groups having the properties imposed on $S$ in the hypotheses, the group $F_n^s$ of Proposition 2.4.2 will also do when $n > 2$.

**Proposition 6.1.7.** *Let $S$ be a simple group such that $\mathrm{Out}(S)$ is infinite, and $P$ a finite group acting transitively on the set $\{1, \ldots, n\}$, $n \geq 2$. Let $B = S^n$ and $A = S \wr P$, so that $A = B \rtimes P$. Then $A$ and $B$ satisfy the four properties of Theorem 6.1.6*

*Proof.* (a) This follows from from [21, Thm 9.2] since $S$ is not a dihedral group. However it can also be seen directly using the facts that $S$ is infinite simple and $P$ is finite as follows: Let $\alpha \in \mathrm{Aut}(A)$, and let $S_i$ denote the $i$-th factor of $B$. Since by construction $B$ is a normal subgroup of $A$ we have that $B\alpha \cap S_i$ is a normal subgroup of $S_i$, so by simplicity we have $B\alpha \cap S_i = S_i$ or 1. But $[S_i : B\alpha \cap S_i] = [(B\alpha)S_i : B\alpha] \leq [A : B\alpha] = n$ so we must have $B\alpha \cap S_i = S_i$, i.e., $S_i \subseteq B\alpha$. Therefore $B \subseteq B\alpha$. Now use $\alpha^{-1}$ to get equality.

(c) Is obvious.

In the rest of the proof, it will be convient at times to use $a^b$ to denote $b^{-1}ab$ when $a$ and $b$ are in the same group, and the action of $b$ on $a$ when $a \in G$ and $b \in \mathrm{Aut}(G)$. We let $\hat{} : \mathrm{Aut}(A) \to \mathrm{Aut}(B)$ be the monomorphism induced by restriction. By [23, Thm. 3.3.20] we know that $\mathrm{Aut}(B)$ is the wreath product $\mathrm{Aut}(S) \wr \Pi_n = (\mathrm{Aut}(S))^n \rtimes \Pi_n$ where $\Pi_n$ is acting on the set $\{1, \ldots, n\}$.

Observe that for $b \in B$ and $p \in P$, we have $bp = pb^p$, so for $\alpha \in \text{Aut}(A)$ if we write $p^\alpha = aq$ where $a \in B$ and $q \in P$, then applying $\alpha$ on both sides of $bp = pb^p$ we get $b^\alpha aq = aqb^{p\alpha} = ab^{p\alpha q^{-1}}q$, so $b^\alpha a = ab^{p\alpha q^{-1}}$, and therefore

$$(*) \qquad \widehat{a}\xi_a = \widehat{pa}q^{-1}$$

where $\xi_a$ is conjugation by $a$ in $B$.

(b) Let $\alpha \in \text{Aut}(A)$ be such that $\widehat{a} = 1$. From $(*)$ we get $\xi_a = pq^{-1}$, and since $\xi_a \in (\text{Aut}(S))^n$ and $pq^{-1} \in \Pi_n$, we must have $pq^{-1} = 1$ and $\xi_a = 1$ so $q = p$, and $a$ is in the center of $B$, which is trivial, so $p^\alpha = p$, and $\alpha = 1$.

(d) Let $\Delta_{\text{Aut}(S)}$ be the diagonal in $(\text{Aut}(S))^n$. We will show that when viewed as a subgroup of $\text{Aut}(B)$ we have $\text{Aut}(A) \subseteq (\Delta_{\text{Aut}(S)} \cdot (\text{Inn}(S))^n) \rtimes \Pi_n$, and therefore $[\text{Aut}(B) : \text{Aut}(A)] \geq |\text{Out}(S)|^{n-1}$. For $\alpha \in \text{Aut}(A)$, if we write $\widehat{a} = \delta\sigma$ where $\delta = (\delta_1, \ldots, \delta_n) \in (\text{Aut}(S))^n$ and $\sigma \in \Pi_n$, then from $(*)$ we get $\delta\sigma\xi_a = p\delta\sigma q^{-1}$ so $\delta\xi_a^{\sigma^{-1}}\sigma = \delta^{p^{-1}}p\sigma q^{-1}$, and since $\delta\xi_a^{\sigma^{-1}}, \delta^{p^{-1}} \in (\text{Aut}(S))^n$ and $\sigma, p\sigma q^{-1} \in \Pi_n$ we get $\delta\xi_a^{\sigma^{-1}} = \delta^{p^{-1}}$ so $\delta^{-1}\delta^{p^{-1}} \in (\text{Inn}(S))^n$. Finally, since the action of $P$ on $\{1, \ldots, n\}$ is transitive, $\delta_1, \delta_2, \ldots, \delta_n$ are all in the same coset of $\text{Out}(S)$, and we can write $(\delta_1, \delta_2, \ldots, \delta_n) = (\delta_1, \delta_1, \ldots, \delta_1) \cdot (1, \delta_1^{-1}\delta_2, \ldots, \delta_1^{-1}\delta_n) \in \Delta_{\text{Aut}(S)} \cdot (\text{Inn}(S))^n$. □

The "difference" between $\text{Aut}(A_n(\mathbf{R}))$ and $\text{Aut}(B_n(\mathbf{R}))$ is further illustrated by looking at the permutations in $\Pi_{n-1}$ that they support. Recall that $\text{Aff}(\mathbf{Z}_{n-1})$ denotes the group of affine transformations of $\mathbf{Z}_{n-1}$, i.e., $\mathbf{Z}_{n-1} \rtimes U_{n-1}$, and $D_{n-1}$ is the subgroup of $U_{n-1}$ generated by the divisors of $n$.

**Theorem 6.1.8.** *Let* $n > 2$.

(a) *The homomorphism* $\pi_n : \text{Aut}(B_n) \to \Pi_{n-1}$ *is surjective.*

(b) *The image of* $\text{Aut}_{PL}(A_n)$ *under* $\pi_n$ *equals* $\mathbf{Z}_{n-1} \rtimes D_{n-1}$.

(c) *The image of* $\text{Aut}(A_n)$ *under* $\pi_n$ *contains* $\mathbf{Z}_{n-1} \rtimes D_{n-1}$ *and is contained in* $\text{Aff}(\mathbf{Z}_{n-1})$.

*Proof.* (a) follows immediately from Theorem 4.4.3 and Lemmas 6.1.1 and 6.1.2.

(b) For one containment, Lemmas 3.2.1 and 3.2.2 imply that $\langle\text{Aut}_{PL}(A_n)\rangle\pi_n \subseteq \mathbf{Z}_{n-1} \rtimes D_{n-1}$. Lemmas 3.2.2 and 3.2.3 give the other inclusion.

(c) For the first containment, $\mathbf{Z}_{n-1} \rtimes D_{n-1} = \langle\text{Aut}_{PL}(A_n)\rangle\pi_n \subseteq \langle\text{Aut}(A_n)\rangle\pi_n$. Now for the second containment, given $\alpha \in \text{Aut}(A_n(\mathbf{R}))$, let $a = 0\alpha$, and $k$ be the self homeomorphism of $\mathbf{R}$ that realizes $\alpha$. Let $h = kt_a^{-1}$ and let $\beta$ be conjugation by $h$. Then $\beta \in \text{Aut}^0(A_n)$, and since $t_a \in A_n(\mathbf{R})$, then by Lemma 1.1.3, $t_a\pi_n$ is a power of the cyclic permutation $\sigma = \begin{pmatrix} 0 & 1 & 2 & \ldots & n-2 \end{pmatrix}$, which is in $\text{Aff}(\mathbf{Z}_{n-1})$. So it suffices to show that $\gamma = \beta\pi_n \in \text{Aff}(\mathbf{Z}_{n-1})$. Let $T = (t_1)\beta$. On one hand, since $T \in A_n(\mathbf{R})$ we have $T\pi_n$ is a power of $\sigma$. On the other hand, since $T = h^{-1}t_1 h$ we have that $T\pi_n$ is $\gamma^{-1}\sigma\gamma$ an $(n-1)$-cycle, so $T\pi_n = \sigma^j$ for some $j$ relatively prime to $(n-1)$. So $T\pi_n = \begin{pmatrix} 0 & j & 2j & \ldots & (n-2)j \end{pmatrix}$, and since $0\gamma = 0$, we must have

$$\gamma = \begin{pmatrix} 0 & 1 & 2 & \ldots & n-2 \\ 0 & j & 2j & \ldots & (n-2)j \end{pmatrix}$$

which is multiplication by $j$ modulo $n-1$, so $\gamma \in \text{Aff}(\mathbf{Z}_{n-1})$. □

We end this part by quoting two results about $\text{Aut}(B_n(\mathbf{R}))$ and $\text{Out}(B_n(\mathbf{R}))$. First, $\text{Aut}(B_n(\mathbf{R}))$ is complete (has trivial center and trivial outer automorphism group so that $\text{Aut}^2(B_n(\mathbf{R})) = \text{Aut}(B_n(\mathbf{R}))$), and second, $\text{Out}(B_n(\mathbf{R}))$ has torsion. These are Theorems 2.4.3 and 5.5.1 together with the comment following the proof of Theorem 5.5.1.

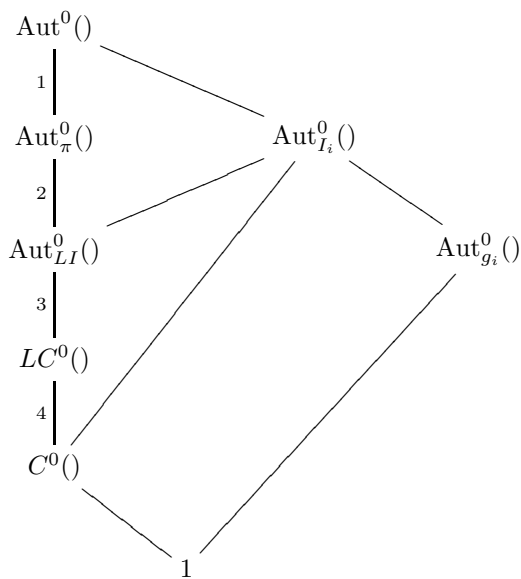

FIGURE 2. Subgroups of the automorphism groups

**Theorem 6.1.9.** *The automorphism group of* $\mathrm{Aut}(B_n(\mathbf{R}))$ *equals* $\mathrm{Aut}(B_n(\mathbf{R}))$.

**Theorem 6.1.10.** *For* $n \geq 4$ *there is an element of order* $n-2$ *in* $\mathrm{Out}(B_n(\mathbf{R}))$.

6.2. **Local properties of automorphisms.** Most of the analysis in Section 5 comes from studying the behavior of automorphisms on single generators. This leads to the concepts *unbent generator at* $i$, *inner*, *inner at* $i$ and *locally inner* in Definition 5.3.2.

The next two propositions tell us that, other than the identity, PL normalizers of $F_{n,\infty}$ can not have unbent generators, and normalizers of $F_{n,\infty}$ with two consecutive unbent generators can not normalize $A_n$. These results are Propositions 5.3.1 and 5.4.5.

**Proposition 6.2.1** (Unbent Generator). *For* $n \geq 2$, *let* $\alpha \in Aut^0(F_{n,\infty}, t_{n-1})$ *have an unbent generator at some* $i$. *If* $\alpha$ *is PL, then* $\alpha$ *is the identity.*

**Proposition 6.2.2** (Double unbent generator). *Let* $\alpha \in Aut^0(F_{n,\infty}, t_{n-1})$, $n \geq 3$, *be a normalizer of* $A_n(\mathbf{R})$. *If for some* $0 \leq i < n-1$, $\alpha$ *has an unbent generator at* $i$ *and* $i+1$, *then* $\alpha$ *is the identity.*

Section 5.3 introduces several subgroups of $\mathrm{Aut}^0(F_{n,\infty}, t_{n-1})$. In the subgroup diagram in Figure 2, we use empty parentheses to denote $(F_{n,\infty}, t_{n-1})$ after "Aut" and to denote $(B_n(\mathbf{R}), t_{n-1})$ after "$C^0$."

The next proposition summarizes what we know about the diagram in Figure 2.

**Proposition 6.2.3.** *In the diagram in Figure 2, containments 1, 2 and 4 are proper, the group* $C^0(B_n(\mathbf{R}))$ *is normal in all groups that contain it and the rectangle at the lower right is a lattice. There are no PL elements in* $\mathrm{Aut}^0_{LI}() - C^0()$, *but there are PL elements in* $\mathrm{Aut}^0_\pi() - \mathrm{Aut}^0_{LI}()$ *(at least when* $n$ *is not a prime power).*

*Proof.* That containments 1 and 4 hold is clear, and the other two follow from Lemma 5.3.3(e) and Corollary 5.3.5.2. Containment 1 is proper by Theorem 4.4.3, 2 is proper by Example 5.3.4 and 4 is proper by Proposition 5.3.1 since lifts are non-trivial and fix at least one generator. That the lower right rectangle is a lattice follows from Lemma 5.3.3(c). The existence and non-existence of PL elements is discussed in Example 5.3.4 and Lemma 5.3.3(e). □

We can say something about structure. For $0 \leq i < n-1$, let $\mathrm{Out}^0_{I_i}(B_n(\mathbf{R}))$ be the image of $\mathrm{Aut}^0_{I_i}(F_{n,\infty}, t_{n-1})$ under the map $\mathrm{Aut}(F_{n,\infty}, t_{n-1}) \to \mathrm{Out}_o(B_n(\mathbf{R}))$, and let $\mathrm{Out}^0_{LI}(B_n(\mathbf{R}))$ be the image of $\mathrm{Aut}^0_{LI}(F_{n,\infty}, t_{n-1})$. To denote the elements of $\mathrm{Aut}^0_{LI}(F_{n,\infty}, t_{n-1})$ that fix $g_i$, we will use $\mathrm{Aut}^0_{LI \cap g_i}(F_{n,\infty}, t_{n-1})$

**Proposition 6.2.4.** *Let*  $0 \leq i < n-1$.

*(a) The sequence*

$$1 \to C^0(B_n(\mathbf{R}), t_{n-1}) \to \mathrm{Aut}^0_{I_i}(F_{n,\infty}, t_{n-1}) \to \mathrm{Out}^0_{I_i}(B_n(\mathbf{R})) \to 1$$

*is split exact with section*  $\mathrm{Out}^0_{I_i}(B_n(\mathbf{R})) \xrightarrow{\sim} \mathrm{Aut}^0_{g_i}(F_{n,\infty}, t_{n-1})$ , *so*

$$\mathrm{Aut}^0_{I_i}(F_{n,\infty}, t_{n-1}) \simeq F_{n,[0,n-1]} \rtimes \mathrm{Aut}^0_{g_i}(F_{n,\infty}, t_{n-1})$$

*(b) The sequence*

$$1 \to C^0(B_n(\mathbf{R}), t_{n-1}) \to \mathrm{Aut}^0_{LI}(F_{n,\infty}, t_{n-1}) \to \mathrm{Out}^0_{LI}(B_n(\mathbf{R})) \to 1$$

*is split exact with a section*  $\mathrm{Out}^0_{LI}(B_n(\mathbf{R})) \xrightarrow{\sim} \mathrm{Aut}^0_{LI \cap g_i}(F_{n,\infty}, t_{n-1})$  *for each $i$. So*

$$\mathrm{Aut}^0_{LI}(F_{n,\infty}, t_{n-1}) \simeq F_{n,[0,n-1]} \rtimes \mathrm{Out}^0_{LI}(B_n(\mathbf{R}))$$

*(c)* $\mathrm{Out}^0_{I_i}(B_n(\mathbf{R}))$ *and* $\mathrm{Out}^0_{LI}(B_n(\mathbf{R}))$ *are torsion free.*

*Proof.* (a) follows from Lemma 5.3.3(c) and remarks in Definition 3.1.1.
(b) This follows immediately from (a), since conjugation by $w \in C^0(B_n(\mathbf{R}), t_{n-1})$ is locally inner.
(c) follows from the splittings in (a) and (b) because a non-trivial, orientation preserving, self homeomorphism of $\mathbf{R}$ has infinite order. □

6.3. **Results relating two values** $(m < n)$. For $2 \leq m < n$ there is an embedding $\tau_{m,n}$ of $F_{m,\infty}$ into $F_{n,\infty}$ taking $t_{m-1}$ to $t_{n-1}$, that extends to an embedding $\Theta_{m,n}$ of $\mathrm{Aut}^0(F_{m,\infty}, t_{m-1})$ into $\mathrm{Aut}^0(F_{n,\infty}, t_{n-1})$. We call the map $\Theta_{m,n}$ the *lift* from $m$ to $n$. The lift $\Theta_{m,n}$ also extends the embedding of $\Pi_m$ into $\Pi_n$. Moreover, locally inner automorphisms are lifted to locally inner automorphisms.

**Theorem 6.3.1.** *Let $2 \leq m < n$ be integers.*

*(a) There is a unique monomorphism of group pairs $\tau_{m,n} : (F_{m,\infty}, t_{m-1}) \to (F_{n,\infty}, t_{n-1})$ so that $(g_{m,i})\tau_{m,n} = g_{n,i}$ for $0 \leq i \leq (m-2)$.*

*(b) there is a monomorphism $\Theta_{m,n} : \mathrm{Aut}^0(F_{m,\infty}, t_{m-1}) \to \mathrm{Aut}^0(F_{n,\infty}, t_{n-1})$ so that if $\theta \in \mathrm{Aut}^0(F_{m,\infty}, t_{m-1})$ then $\theta' = \theta\Theta_{m,n}$ is the only map that makes the* *following diagram commute.*

$$
\begin{array}{ccccc}
(F_{m,\infty}, t_{m-1}) & \xrightarrow{\tau_{m,n}} & (F_{n,\infty}, t_{n-1}) & \longleftarrow & (\langle g_{m-1}, \ldots, g_{n-2}, t_{n-1} \rangle, t_{n-1}) \\
\downarrow{\scriptstyle \theta} & & \downarrow{\scriptstyle \theta'} & & \downarrow{\scriptstyle 1} \\
(F_{m,\infty}, t_{m-1}) & \xrightarrow{\tau_{m,n}} & (F_{n,\infty}, t_{n-1}) & \longleftarrow & (\langle g_{m-1}, \ldots, g_{n-2}, t_{n-1} \rangle, t_{n-1})
\end{array}
$$

*(c) The lift $\Theta_{m,n}$ preserves permutations, i.e., the following diagram commutes.*

$$
\begin{array}{ccc}
\mathrm{Aut}^0(F_{m,\infty}, t_{m-1}) & \xrightarrow{\Theta_{m,n}} & \mathrm{Aut}^0(F_{n,\infty}, t_{n-1}) \\
\downarrow{\scriptstyle \pi_m} & & \downarrow{\scriptstyle \pi_n} \\
\Pi_m & \lhook\joinrel\longrightarrow & \Pi_n
\end{array}
$$

*(d) The lift $\Theta_{m,n}$ maps $\mathrm{Aut}^0_{LI}(F_{m,\infty}, t_{m-1})$ to $\mathrm{Aut}^0_{LI}(F_{n,\infty}, t_{n-1})$. The map $\tau_{m,n}$ maps $F_{m,[0,m-1]}$ to $F_{n,[0,n-1]}$. If we define $\Upsilon_{m,n} : F_{m,[0,m-1]}^{m-1} \to F_{n,[0,n-1]}^{n-1}$ by $(W_0, \ldots, W_{m-2})\Upsilon_{m,n} = (W_0\tau_{m,n}, \ldots, W_{m-2}\tau_{m,n}, 1, \ldots, 1)$ then the following diagram commutes.*

$$
\begin{array}{ccc}
\mathrm{Aut}^0_{LI}(F_{m,\infty}, t_{m-1}) & \xrightarrow{\Theta_{m,n}} & \mathrm{Aut}^0_{LI}(F_{n,\infty}, t_{n-1}) \\
\downarrow{\scriptstyle \iota_m} & & \downarrow{\scriptstyle \iota_n} \\
F_{m,[0,m-1]}^{m-1} & \xrightarrow{\Upsilon_{m,n}} & F_{n,[0,n-1]}^{n-1}
\end{array}
$$

*Note that the maps $\iota_m, \iota_n$ are not homomorphisms.*

*Proof.* (a), (b), and (c) are just Proposition 4.1.2, Theorem 4.1.6, and Proposition 4.4.2 respectively.

(d) The statement about $\Theta_{m,n}$ is Corollary 5.3.5.1. The statement about $\tau_{m,n}$ follows from Lemma 4.1.5. Commutativity follows from Lemma 5.3.5. $\qquad\blacksquare$

In contrast with Theorem 6.3.1 part (b), standard lifts do not normalize $A_n(\mathbf{R})$, as we get from Theorem 5.4.7.

**Theorem 6.3.2.** *(a) For integer values $2 \leq m < n - 1$, there is an isomorphic copy of $Aut^0(F_{m,\infty}, t_{m-1})$ in $\mathrm{Aut}(B_n(\mathbf{R}))$ that has trivial intersection with $\mathrm{Aut}(A_n(\mathbf{R}))$. (b) For integer values $2 \leq m = n - 1$, there is an isomorphic copy of $Aut^0_{LI}(F_{m,\infty}, t_{m-1})$ in $\mathrm{Aut}(B_n(\mathbf{R}))$ that has trivial intersection with $\mathrm{Aut}(A_n(\mathbf{R}))$.*

However, when $m < n$ are such that $(m-1)|(n-1)$, we can do symmetric lifts from $m$ to $n$. For $\alpha \in \mathrm{Aut}^0(F_{m,\infty}, t_{m-1})$, the symmetric lift $\alpha\Lambda_{m,n}$ can be thought of as the composition of the standard lift $\alpha\Theta_{m,n}$ with its iterated $(m-1)$-step rotations. The lift $\alpha\Lambda_{m,n}$ normalizes $F_{n,\infty}$ and commutes with $t_{m-1}$, but in general it does not normalize $A_n(\mathbf{R})$. In the particular case $m = 2$ we get that symmetric lifts from $\mathrm{Aut}^0(F_{2,\infty}, t_1)$, do normalize $A_n(\mathbf{R})$. From Theorems 4.3.1 and 5.4.8 we get the following.

**Theorem 6.3.3.** *Let $2 \leq m < n$ be such that $(m-1)|(n-1)$.*

*(a) There is a monomorphism $\Lambda_n : Aut^0(F_{m,\infty}, t_{m-1}) \to \mathrm{Aut}^0(F_{n,\infty}, t_{m-1})$ so that if $\theta \in \mathrm{Aut}^0(F_{m,\infty}, t_{m-1})$, then $\theta' = \theta\Lambda_{m,n}$ is the only homomorphism in*

$\mathrm{Aut}^0(F_{n,\infty}, t_{m-1})$ *that makes the following diagram commute.*

$$
\begin{array}{ccc}
(F_{m,\infty}, t_{m-1}) & \xrightarrow{\tau_{m,n}} & (F_{n,\infty}, t_{n-1}) \\
\downarrow{\scriptstyle\theta} & & \downarrow{\scriptstyle\theta'} \\
(F_{m,\infty}, t_{m-1}) & \xrightarrow{\tau_{m,n}} & (F_{n,\infty}, t_{n-1})
\end{array}
$$

*In particular if $\theta \in \mathrm{Aut}^0(F_{2,\infty}, t_1)$ then $\theta\Lambda_{2,n}$ normalizes $A_n(\mathbf{R})$.*

(b) *If $m > 2$ and $\theta \in \mathrm{Aut}(A_m(\mathbf{R}))$ is conjugation by $1 \neq w \in C^0(B_m(\mathbf{R}), t_{m-1})$ then $\theta\Lambda_{m,n}$ does not normalize $A_n(\mathbf{R})$.*

## 7. Questions

We gather some questions raised by the material in the paper. Recall the graded groups $\mathcal{A}$ and $\mathcal{LC}$ defined in Section 5.3. We introduce two more graded subgroups of $\mathcal{A}$ to hang questions on. We let $\mathcal{LRC}$ be the smallest graded subgroup of $\mathcal{A}$ that contains all the $C^0(B_n(\mathbf{R}), t_{n-1})$ and that is closed under lifts and rotations. Note that $\mathcal{A}$ is closed under rotations. We let $\mathcal{LRPL}$ be the smallest graded subgroup of $\mathcal{A}$ that contains all the $\mathrm{Aut}^0_{PL}(F_{n,\infty}, t_{n-1})$ and that is closed under lifts and rotations.

**Question 7.1.** Is $\mathcal{LRPL} = \mathcal{A}$? This asks if we have found all automorphisms of $B_n(\mathbf{R})$.

**Question 7.2.** Is $\mathcal{LRC} = \mathcal{LRPL}$?

**Question 7.3.** Is $\mathcal{LC} = \mathcal{LRC}$?

Example 5.3.4 shows that $\mathcal{LC}$ is a proper subset of $\mathcal{LRPL}$, so at least one of Question 7.2 or 7.3 has answer no.

**Question 7.4.** Is containment 3 in the diagram above Proposition 6.2.3 proper?

To answer at least some of the questions above, it seems that techniques will have to be developed that go beyond the techniques of the current paper.

**Question 7.5.** Is $\mathrm{Aut}^0_{LI}(F_{n,\infty}, t_{n-1})$ preserved by rotations? Rotations preserve $\mathrm{Aut}^0(F_{n,\infty}, t_{n-1})$ by Lemma 4.2.1.

**Question 7.6.** Is there torsion in $\mathrm{Out}(B_3(\mathbf{R}))$?

**Question 7.7.** Is there torsion in $\mathrm{Out}_\pi(B_n(\mathbf{R}))$?

The example of torsion that we construct in Section 5.5 has a non-trivial permutation and only works for $n \geq 4$. Is there an example of torsion that does not permute the residues?

**Question 7.8.** Is there a bound on the torsion in $\mathrm{Out}(B_n(\mathbf{R}))$? Even if non-trivial permutations have to be involved, there are relevant permutations with order much higher than that of our example.

**Question 7.9.** What is the image of $\mathrm{Aut}(A_n)$ under $\pi_n$? See Theorem 6.1.8(c).

**Question 7.10.** Is $\mathrm{Aut}(A_n(\mathbf{R}))$ normal in $\mathrm{Aut}(B_n(\mathbf{R}))$? Is the automorphism tower of $A_n(\mathbf{R})$ contained in $\mathrm{Aut}(B_n(\mathbf{R}))$? Is $\mathrm{Aut}^2(A_n(\mathbf{R})) = \mathrm{Aut}(A_n(\mathbf{R}))$? These questions ask about increasingly strict conditions on $\mathrm{Aut}(A_n(\mathbf{R}))$. There is no particular reason to suspect that even the first question has an affirmative answer.

**Question 7.11.** For $2 < m < n$ with $(m-1) \mid (n-1)$, does image of $\Lambda_{m,n}$ defined in Theorem 4.3.1 have trivial intersection with $\mathrm{Aut}(A_n(\mathbf{R}))$? Theorem 5.4.8 shows that the image of the inner automorphisms of $A_m(\mathbf{R})$ under $\Lambda_{m,n}$ has trvial intersection with $\mathrm{Aut}(A_n(\mathbf{R}))$.

**Question 7.12.** Does the image of $\Lambda_{2,n}$ contain a representative of every element of $\mathrm{Out}_o(A_n(\mathbf{R}))$? This asks if we have found all automorphisms of $A_n(\mathbf{R})$. This breaks down into two parts. Is the image of $\Lambda_{2,n}$ equal to all of $\mathrm{Aut}^0(F_{n,\infty}, t_1)$ modulo $C^0(B_n(\mathbf{R}), t_{n-1})$? Is there are representative of every element of $\mathrm{Out}_o(A_n(\mathbf{R}))$ in $\mathrm{Aut}^0(F_{n,\infty}, t_1)$? There is no reason to expect either question to have an affirmative answer.

The authors have thought about the above questions at least to some extent, and have not pursued such obvious questions as: what is the full structure of the automorphism and outer automorphism groups, do they have free subgroups, are any of them finitely generated or finitely presented, etc.